\documentclass[12pt]{article}
\usepackage{amsmath,amssymb}
\usepackage{amsthm,amscd}
\usepackage{graphicx,tikz}
\usepackage{enumerate}
\usepackage{titlesec}
\usepackage{epstopdf}
\usepackage{caption}
\usepackage{xcolor}
\usepackage[sans]{dsfont}
\numberwithin{equation}{section}
\newtheorem{Thm}{Theorem}[section]
\newtheorem*{Thm*}{Theorem}
\newtheorem{Prop}[Thm]{Proposition}
\newtheorem*{Prop*}{Proposition}
\newtheorem{Lem}[Thm]{Lemma}

\newtheorem{Cor}[Thm]{Corollary}
\newtheorem*{Cor*}{Corollary}
\newtheorem{Fact}[Thm]{Fact}

\newtheorem{Exa}[Thm]{Example}

\theoremstyle{remark}

\theoremstyle{definition}

\newtheorem{Def}[Thm]{Definition}

\newtheorem*{Def*}{Definition}

\numberwithin{equation}{section}

\newcommand{\g}[1]{{\mathfrak #1}}
\newcommand{\m}[1]{\mathbb{ #1}}
\newcommand{\mc}[1]{\mathcal{ #1}}

     \def\ol{\overline}    
   
\def\al{\alpha}       \def\be{\beta}        \def\ga{\gamma}
\def\de{\delta}       \def\eps{\varepsilon}  \def\ze{\zeta}
\def\th{\theta}       %\def\im{\imath}
\def\ka{\kappa}       \def\la{\lambda}      
\def\si{\sigma}                
\def\ph{\varphi}               
\def\om{\omega}              \def\De{\Delta}
\def\La{\Lambda}             \def\Ph{\Phi}

\theoremstyle{definition}

\theoremstyle{remark}

\newtheorem{Rmq}[Thm]{Remark}

\numberwithin{equation}{section}
\newfont{\goth}{eufm10 at 12pt}
\newfont{\gots}{eufm8 at 9pt}

\def\bt{\begin{Thm}}
\def\et{\end{Thm}}
\def\br{\begin{Rmq}}
\def\er{\end{Rmq}}

\def\bc{\begin{Cor}}
\def\ec{\end{Cor}}
\def\bp{\begin{Prop}}
\def\ep{\end{Prop}}
\def\bl{\begin{Lem}}
\def\el{\end{Lem}}
\def\bd{\begin{Def}}
\def\ed{\end{Def}}
\def\bq{\begin{quotation}}
\def\eq{\end{quotation}}
\def\bfa{\begin{Fact}}
\def\efa{\end{Fact}}
\def\bexa{\begin{Exa}}
\def\eexa{\end{Exa}}

\def\ra{\rightarrow}
\def\vs{\vspace{1em}}

\begin{document}
%\mbox{ }\vspace*{-5em}
\title{
Convolution and square 
in abelian groups II
}
\author{Yves Benoist
}
\date{}

\maketitle
%\vspace{-2em}
%\centerline{\footnotesize \mbox{Preliminary version } }

\begin{abstract}
\noindent 
A critical value on an abelian group $G$ of odd order $d$ is a value $\la$ such that 
the  functional equa\-tion $f\!\star\! f (2\,t) =\la f(t)^2$ on $G$  has a nonzero solution $f$. 
We construct many  critical  values by using
abelian varieties with complex multiplication.
\end{abstract}
\renewcommand{\thefootnote}{\fnsymbol{footnote}} 
\renewcommand{\thefootnote}{\arabic{footnote}} 
{\footnotesize \tableofcontents}\nopagebreak
%\newpage
%10
\section{Critical values}

%11
\subsection{Introduction}
\label{secintexp}

This paper, is the sequel of \cite{CSAGI} in which we introduced 
the notion of $d$-critical values or critical values on the cyclic group $\m Z/d\m Z$, where $d$ is an odd integer. Before recalling the precise definition of critical values  in Section \ref{secconsqu}, I would like to sum up the content of this first paper.
In \cite{CSAGI} we  reported some striking numerical experiments
and  explained that the value $\la:=\sqrt{a} +i\sqrt{b}$ is a $d$-critical value
on $\m Z/d\m Z$
when $a+b=d$ and $a\equiv \frac{(d+1)^2}{4}$ mod $4$. It was surprising that Jacobi theta functions and elliptic curves with complex multiplication were needed to prove that these apparently very simple values
are critical. 
\vs

\noindent
\begin{minipage}[t]{.52\textwidth}
\hspace*{2em}In the experimental lists of critical values given in \cite[Section 1.5]{CSAGI}, there still  remained 
intriguing $d$-critical values
that could not be  explained by 
the technics of \cite{CSAGI} relying on
Jacobi theta functions and elliptic curves. 
\end{minipage}\!\!
\hfill\!\!
\noindent
\raisebox{-1.5ex}{\fbox{
\begin{minipage}[t]{.4\textwidth}
\noindent
\centerline{$d=11$ $\la=1\!+\!\sqrt{5}\! +\! i\sqrt{\! 5\!-\!2 \sqrt{5}}$} \\
\centerline{$d=15$ $\la=1\!+\!\sqrt{5}\! +\! i\sqrt{\! 9\!-\!2 \sqrt{5}}$} \\
\centerline{$d=15$ $\la\!=\! 2\sqrt{\! 2\!-\!\sqrt{3}}\! +\!2i\! +\!i\sqrt{3}   $} \\
\centerline{$d=15$\; $\la=(\!\sqrt{3}\!+\! i\sqrt{2})\,(\!\sqrt{2}\!+\! i)$} \\
\centerline{$d\!=\!17$\; $\la\!=\!1\!+\! 2\sqrt{2}\! +\! 2i\sqrt{\!2\!-\!\sqrt{2}}$}  
\end{minipage}
}}
\vs 

The aim of the present paper is to explain 
a general construction of critical values $\la$
by using torsion points on  abelian varieties
and Riemann theta functions 
(Theorem \ref{thmmaiabe}). 
This construction gives all the known $d$-critical values with $d\leq 17$.  
As an output, focusing on abelian varieties with complex multiplication, we will obtain new explicit  critical values. 

%12
\subsection{Constructing critical values}
\label{secconsqu}

Before going on our introduction, we need to recall the definition of critical value from \cite{CSAGI}.
We let $G$ be a finite abelian group of odd order $d$,
but we keep in mind that the case where $G$ is the cyclic group $G=\m Z/d\m Z$
is very important. A nonzero function $f: G\ra \m C$ solution of the functional equation
\begin{equation*}
\label{eqntante2} \textstyle
\sum_{\ell\in G}f(k\!+\!\ell)\,f(k\!-\!\ell)
\: =\; \la \, f(k)^2
 \;\;\; \mbox{\rm for all $k$ in $G$,}
\end{equation*}
where $\la\in \m C$ is a parameter,
will be called a ``$\la$-critical
function on $G$'',
and a value $\la$ for which such a function $f$ 
exists will be called a ``critical value on $G$'', or a ``$d$-critical value''
when $G=\m Z/d\m Z$. 
This equation expresses 
a proportionality condition between the ``convolution square'' of $f$ and 
its ``multiplication square''.

Our main result  in this paper
is a construction of explicit critical values.
Beyond this precise list of critical values, 
the main interest and surprise is the fact that 
this naive question is related to abelian varieties with complex multiplication and to modular forms on the Siegel upper half-space.

We will explain the construction of these  critical values
from various points of view. 
First from the point of view of Riemann theta functions (Theorem \ref{thmmaithe}).
Then from the point of view of abelian varieties (Theorem \ref{thmmaiabe}).
And finally we will apply this construction from the point of view of number fields with complex multiplication
(Theorem \ref{thmmaimul}).

Eventhough the consequences of these theorems,
as Propositions \ref{proparacri}, \ref{prorarbrc1},
and \ref{profibonacci},  are  elementary, they do not have an
elementary proof.
Indeed,
our construction relies on the
Riemann theta functions $\theta(z,\tau)$ for special values of the Riemann matrix $\tau$,
and the key point in the proof relies  
on  modularity properties of these theta functions on the Siegel upper half-space
which is due to Siegel
(Lemma  \ref{lemtrafor}) and on its refinements due  to Igusa (Corollary \ref{cortrafor}) 
and to Stark (Lemma \ref{lemstasty}). 
It also relies on a construction of principally polarized CM abelian varieties 
due to Taniyama-Shimura (Fact \ref{facshitan}).
\vs

%13
\subsection{A few concrete examples of critical values} 
\label{secfewexa}

To give a flavor of the output of our method we just give here three concrete families of critical values
that  will be obtained, respectively, as part  of Proposition 
\ref{proparacri}, Proposition \ref{prorarbrc1}  
and Corollary \ref{corfibonacci}.

\begin{Cor*}
For $j=1,2$, let $d_j\!=\!a_j\!+\!b_j$  be positive integers with $d_1\wedge d_2=1$ and
$a_j\! -\! \frac{(d_j+1)^2}{4}\!\equiv\! 2$ {\rm mod}~$4$ .
Then  
$\la\!=\!(\sqrt{a_1}\!+\! i\sqrt{b_1})(\sqrt{a_2}\!+\! i\sqrt{b_2})
$ 
is $d$-critical.
\end{Cor*}

\begin{Cor*}
Let $d=a+b+c$ be positive integers with  
 $b^2>4ac$. \\
Assume either that $ b\!\equiv\! c\!\equiv\!  0$ {\rm mod} $4$
or that $ b\!\equiv\! c\!\equiv\!  a$ {\rm mod} $4$.\\ 
If $a\!\equiv 1$ {\rm mod} $4$,
then 
$
\la\!=\!\sqrt{a}\!+\!\sqrt{c}\! +\! i\sqrt{b\!-\! 2 \sqrt{ac}}\; 
$ 
is  $d$-critical.\\
If $a\!\equiv  3$ {\rm mod} $4$,
then 
$
\la\!=\! \sqrt{b\!-\! 2 \sqrt{ac}}\! +\! i\sqrt{a}\!+\! i\sqrt{c}\; 
$ 
is  $d$-critical.
\end{Cor*}

\begin{Cor*}
Let $n\geq 5$  be a prime number  and $L_n$ be the $n^{th}$ Lucas number. Then  
$\la_n=\prod\limits_{1\leq k\leq (n-1)/2}
(1+2i \sin(\frac{k\pi}{n}))$ or $-\la_n$
is  $L_n$-critical.
\end{Cor*}

%In the last corollary the critical value is determined only up to a $\pm $ sign. 
%Computing the sign would require extra work. 

%14
\subsection{Properties of critical values}

It is important to keep in mind a few properties, proven in \cite[Section 2.1]{CSAGI}.

\bp
\label{procrival} Let $\la$ be a critical value on  an abelian group $G$ of odd order $d$.\;
$(i)$ All the Galois conjugates of $\la$ are also critical values on $G$.\\ 
$(ii)$ One has $|\la|\leq d$ with equality if and only if $\la=d$.\\
$(iii)$
The ratio $d/\la$ is also a critical value on $G$.\\
$(iv)$  $\la$ is an algebraic integer such that  $\la\equiv 1$ mod $2$.\\
$(v)$ In particular, there exist only finitely many critical values on $G$.
\ep

\noindent
The condition ``$\la\!\equiv\! 1$ mod $2$'' means that the ratio $\frac{\la -1}{2}$ is an algebraic integer.

Since there are only finitely many critical values $\la$ on $G$, 
it will be surprising to see that, for all the critical values
$\la$ we will construct in this paper, the variety of $\la$-critical functions $f$ on $G$ is higher dimensional.

%15
\subsection{Strategy, organization and main results}

This paper is independent of \cite{CSAGI}, but it might be helpful for the reader to  be familiar with the proof given in  \cite{CSAGI}
that relies on Jacobi theta functions, and on elliptic curves with complex multiplication. The proof here will follow the same lines
replacing the Jacobi theta functions by the Riemann theta functions, 
the elliptic curves by abelian varieties, and the imaginary quadratic fields by CM number fields.
\vs 

{\bf In  Chapter \ref{secthefun}}, we explain 
how one can use the Riemann theta functions 
for constructing critical fonctions 
on any finite abelian group $G$ of odd order (Theorem \ref{thmmaithe}). 
Such an abelian group $G$ can be seen as a quotient 
${\bf d}^{-1}\m Z^g/\m Z^g$ for an integral matrix
${\bf d}$  with $\det( {\bf d})$ odd.
The  Riemann theta functions $\theta_\tau(z)$ are $\m Z^g$-periodic
functions on $\m C^g$ parame\-trized by a  matrix $\tau$
that lives in the Siegel upper half space $\mc H_g$,
i.e. $\tau$ is a complex symmetric matrix whose imaginary part is positive definite.
The key point is a condition on the Riemann matrix $\tau$
(Lemma \ref{lemratequ})
that ensures that the restriction 
of any translate of the function  $\theta_\tau(z)$ to the group $G$ is critical.

In this higher dimensional case, this condition, which involves $2^g-1$ 
equations while there are only $g(g+1)/2$ parameters,
seems difficult to satisfy.
However a nice fact due to Igusa (lemma \ref{cortrafor})
helps us to construct solutions of these equations: 
it is sufficient (and conjecturally also necessary) that 
there exists an element $\si$ of {\it the integral symplectic theta group 
${\rm Sp}^{\theta,2}_{g,\m Z}$ of level $2$} such that
$\si\,\tau={}^t{\bf d}\tau{\bf d}$. This 
group ${\rm Sp}^{\theta,2}_{g,\m Z}$ is a normal finite index 
subgroup of the integral symplectic group ${\rm Sp}(g,\m Z)$.
It acts naturally on the Siegel upper half space $\mc H_g$.

\begin{Thm*}\!{\bf \ref{thmmaithe}}
 Let $\tau\in \mc H_g$ and 
${\bf d}\in \mc M(g,\m Z)$ with $\det( {\bf d})$ odd.\\
Assume that there exists 
$
\si\!=\!\mbox{\scriptsize 
$\left(\!\begin{array}{cc} \al&\be\\  
\ga&\de\end{array}\!\right)$} \in {\rm Sp}_{g,\m Z}^{\theta,2}
$  
such that $\si\tau={}^{t}{\bf d}\tau{\bf d}$.
Then the function 
$\theta_\tau$ restricted to the group $G={\bf d}^{-1}\m Z^g/\m Z^g$
is $\la$-critical for  a critical value $\la:=\ka \, {\rm det}_{\m C}(\ga\tau+\de)^{1/2}|G|$
with $\ka^8=1$.
\end{Thm*}

{\bf In  Chapter \ref{secsymgro}}, we reinterpret our construction of critical values
in terms of the integral symplectic  theta subgroup of level $2$
(Corollary \ref{cormaisym}).
We also explain (Lemma \ref{lemthestr}) how to construct elements of 
the integral symplectic  theta group of level $2$
starting from an element of the rational symplectic  theta subgroup ${\rm Sp}^{\theta,2}_{g,\m Q}$  of level $2$. This relies on the {\it symplectic adapted basis theorem} (Proposition \ref{prohsidsi}).
Finally we explain how to construct easily
elements of the rational symplectic  theta subgroup  ${\rm Sp}^{\theta,2}_{g,\m Q}$ of level $2$
thanks to the Cayley transform (Lemma \ref{lemcaythe}).
\vs 

{\bf In Chapter \ref{secabevar}}, we are dealing 
with a principally polarized abelian variety $(A=\m C^g/\La,\om)$,
with its hermitian structure $H$ on its Lie algebra $\m C^g$
and with its integral symplectic form $\om={\rm Im}(H)$ on its lattice $\La$.
We reinterpret our general construction of critical functions and critical values
from the point of view of abelian varieties
in Theorem \ref{thmmaiabe} which is the main theorem of this paper. Note that the Riemann theta functions do not occur in the statement of this theorem, but only in its proof.

\begin{Thm*}\!{\bf \ref{thmmaiabe} }
Let $(A=\m C^g/\La,\om)$ be a principally polarized abelian variety, 
$\nu$ be  a unitary $\m Q$-endomorphism of $A$
preserving a theta structure of level $2$,\; $T_\nu$  its tangent map, 
$
G_\nu:=\La/(\La\cap \nu\La)
$
and $d_\nu:=|G_\nu|$.
Then there exists a critical value
$\la_\nu=\ka_\nu\, d_\nu^{1/2}\,{\rm det}_{\m C}(T_\nu)^{1/2}$
on  the group
$G_\nu$ with $\ka_\nu^4=1$.
\end{Thm*}
\noindent
The square $\ka_\nu^2=\pm 1$ can be calculated from the condition $\la_\nu\equiv 1$ mod $2$.

The finite abelian group $G_\nu$ depends not only on $\nu$ but also on $\La$. It might be cyclic even when $g>1$.
\vs

{\bf In  Chapter \ref{seccommul}}, we are dealing with a CM number field $K$
of degree $2g$, with a  CM type $\Phi:K\rightarrow \m C^g$ and 
with a lattice $\La\subset K$, and we assume 
that  the abelian variety is the CM abelian variety  
$A=\m C^g/\Phi(\La)$.
We specialize our general construction of critical values to that case and  express it from the point of view of CM number fields (Theorem \ref{thmmaimul}):\\
$\star$ 
The symplectic form on $\La$ is given by 
a nonzero imaginary element $t_0$ of $K$ thanks to the simple formula 
$\om(x,x')={\rm Tr}_{K/\m Q}(\frac{x\ol{x'}}{t_0})$.
A key point is to choose $t_0$ so that this symplectic form 
is integral with determinant $1$ on $\La$.\\
$\star$
The unitary $\m Q$-endomorphisms $\nu$  are nothing 
but  elements $\nu\in K$ of absolute value $1$. 
On can construct such  $\nu$ that preserve a theta structure of level $2$
thanks to the Cayley transform
(Lemma \ref{lemimauni}). \\
$\star$
The ratio $\la_\nu^2/d_\nu$ 
is, up to sign, equal to the {\it reflex norm} 
$N_\Phi(\nu)$.
In particular the critical value $\la_\nu$ is a $d_\nu$-Weil number.
\vs 

{\bf In  Chapter \ref{seclisexa}},  we show on examples how to 
compute, at least up to sign, explicit critical values by using 
Theorem \ref{thmmaimul}.

In Section \ref{secimaqua} we  discuss the case where $K$ is an imaginary quadratic field
and hence $A$ is a CM elliptic curve. 
This case is the one we studied  in \cite{CSAGI}.

In Section \ref{secproima}, we discuss the case where
$K$ is a product of two imaginary fields or, equivalently, 
$A$ is isogenous to the product of two elliptic curves with complex multiplication.
This gives the  $1^{\rm st}$ corollary of Section \ref{secfewexa}.

In Section \ref{secquafie}, we discuss the case where $K$ is a quartic CM fields 
and hence $A$ is a CM abelian surface.
This gives the $2^{\rm nd}$ corollary of Section \ref{secfewexa}.\vs

{\bf In Chapter \ref{secabemok}}, we come back to the general CM abelian varieties, but we  specialize our 
Theorem \ref{thmmaimul} 
to the case where the lattice $\La$ is a fractional ideal $\g m$ of $K$, or equivalently, to the case 
where the abelian variety $A=\m C^g/\Phi(\La)$
has multiplication by $\mc O_K$ (Theorem \ref{thmmaimul2}).
This case is particularly nice because of the following three reasons:\\
$\star$
A theorem of Taniyama-Shimura relying on class field theory tells us exactly, 
for which CM type $\Phi$ and for which ideal $\g m$
this abelian variety $A$ is principally polarized (Fact \ref{facshitan}).\\
$\star$ The group $G_\nu$ does not depend on the ideal $\g m$.\\
$\star$ 
Exemples of unitary elements $\nu\in K$ that preserve a theta structure of level $2$ are 
$\nu=\mu/\ol{\mu}$ with $\mu=1+s-\ol{s}$ where $s\in K$ has deno\-minator prime to $2$ and where the norm $N_{K/\m Q}(\mu)\in \m Q$ has an odd numerator
(Lemma \ref{lemitkl2}).  

We quote below  the part of Corollary \ref{cormaimul2} where the extension $K/(K\cap \m R)$ is {\it ramified} i.e. ``ramified at a finite place''.
 
\begin{Cor*}{\bf  \ref{cormaimul2}.B} 
Let $K$ be a CM field such that $K/(K\cap \m R)$ is ramified.
Let $s\in \mc O_K$, $\mu:=1+s-\ol{s}$ 
with $N_{K/\m Q}(\mu)$ odd.
Then for all CM types $\Phi$ of $K$, there exists
a critical value $\la_\mu=\ka_\mu N_\Phi(\mu)$ on $G_\mu:=\mc O_K/\mu\mc O_K$  with $\ka_\mu ^4=1$.
\end{Cor*}
Note that, the principally polarized CM abelian varieties and their theta structures of level $2$ do not occur in the statement of this corollary, but only in its proof.

In Section \ref{seccycfie}, we focus on a very important case, when $K$ is the cyclotomic field $K_n=\m Q[\zeta_n]$. In this case, 
the conclusion of this corollary is always true even when 
$K$ is unramified over its maximal real subfield $K\cap \m R$.
%Indeed, for all CM type $\Phi$ of $K_n$ 
%one can always find an ideal $\g m$ of $K_n$ 
%such that the abelian variety $\m C^g/\Phi(\g m)$ 
%is principally polarized (Corollary \ref{corcycfie}). 

In  Section \ref{secfibluc},
we take $s=\zeta_n$.
This gives the $3^{\rm rd}$  corollary of Section \ref{secfewexa}.\vs 

{\bf In  Chapter \ref{secsigcri}},  we come back
to the examples of Chapter \ref{seclisexa}, and 
explain how to remove the remaining ambiguity on the sign of the critical values.
The key point is a precise formula for the theta cocycle 
$j(\si,\tau)$ 
in Lemma \ref{lemstasty}
which is due to Stark, combined with tricky explicit computations. This theta cocycle is the one that shows up in the transformation formula 
for the Riemann theta functions.\vs 

I would like to thank  E. Ullmo and R. Salvati Manni for useful comments on this project.

%20
\section{Theta functions}
\label{secthefun}

The aim of this chapter is to explain our general construction of critical values from the point of view of 
the Riemann theta functions (Theorem \ref{thmmaiabe}).

%21
\subsection{Riemann matrices and symplectic group}
\label{secriemat}

We begin by preliminary classical notation and definition
(see \cite{Deb99} or \cite{Beau13}).
\vs 

A Riemann matrix $\tau$ is a complex symmetric matrix  whose imaginary part is positive definite.
For $g\geq 1$, let $\mc H_g$ be the Siegel upper half-space which is the space 
of Riemann matrices of size $g$,
$$
\mc H_g\; =\;
\{ \tau\in \mc M(g,\m C)\mid\; {}^t\tau=\tau ,\;\; {\rm Im}\, \tau >0\}.
$$

Let ${\rm Sp}(g,\m R):=\{ \si\in {\rm GL}(2g,\m R)\mid \; {}^t\si J\si =J\}
$,
where 
$
J=\mbox{\scriptsize 
$\left(\!\begin{array}{cc} {\bf 0}&\mathds{1}_g\\  
-\mathds{1}_g&{\bf 0}\end{array}\!\right)$} 
$,
be the real symplectic group. This group is the stabilizer of the symplectic form 
$\om $ on $\m R^{2g}$  given by 
$$
\om(x,y)={}^tx\, J\, y\, ,
$$
and seen as a group of
$2$ by $2$ block matrices of size $g$, it is given by
\begin{align}
\label{eqnspgr}
\nonumber
{\rm Sp}(g,\m R)&=\{
\si=\mbox{\scriptsize 
$\left(\!\begin{array}{cc} \al&\be\\  
\ga&\de\end{array}\!\right)$} \mid \; 
{}^t\al\ga={}^t\ga\al,\;\;
{}^t\be\de={}^t\de\be,\;\;
{}^t\al\de-{}^t\ga\be=\mathds{1}_g
\},\\
&= 
\{
\si=\mbox{\scriptsize 
$\left(\!\begin{array}{cc} \al&\be\\  
\ga&\de\end{array}\!\right)$} \mid \; 
\si^{-1}=\mbox{\scriptsize 
$\left(\!\begin{array}{cc} {}^t\de&-{}^t\be\\  
-{}^t\ga&{}^t\al\end{array}\!\right)$}\;\}.
\end{align} 
The group ${\rm Sp}(g,\m R)$ acts transitively on the Siegel 
upper half-space $\mc H_g$,
$$
\si \tau:=(\al \tau +\be)(\ga \tau +\de)^{-1}.
$$
One cannot confuse this notation $\si\tau$ with the product of matrices since $\si$ has size $2g$ while $\tau$ has size $g$.
The stabilizer of the element $\tau_0=i\mathds{1}_g\in \mc H_g$ is 
the unitary group $U(g,\m R)$, 
so that $\mc H_g\simeq {\rm Sp}(g,\m R)/U(g,\m R)$.

Let ${\rm Sp}(g,\m Z):={\rm GL}(2g,\m Z)\cap {\rm Sp}(g,\m R)$
be the integral symplectic group. 
For $\ell\geq 1$, the following subgroups of ${\rm Sp}(g,\m Z)$ play 
an important role in the theory of theta functions.
The first one is the {\it  integral congruence symplectic group ${\rm Sp}_{g,\m Z}^{\ell}$ 
of level $\ell$}.
$$
{\rm Sp}_{g,\m Z}^{\ell}:= \{ \si\in {\rm Sp}(g,\m Z)\mid \si \equiv\mathds{1}_{2g}\;{\rm mod}\; \ell\}.
$$
The second one is the {\it integral symplectic theta group ${\rm Sp}_{g,\m Z}^{\theta,\ell}$ of level $\ell$.}
\begin{equation*}
\label{eqnthestr}
{\rm Sp}_{g,\m Z}^{\theta,\ell}:= 
\{
\si\!=\!\mbox{\scriptsize 
$\left(\!\!\begin{array}{cc} \al\!&\be\\  
\ga\!&\de\end{array}\!\!\right)$} \in {\rm Sp}_{g,\m Z}^{\ell} \mid \,
({}^t\al\ga)_0\equiv ({}^t\be\de)_0\equiv 0\; {\rm mod}\; 2\ell\}, 
\end{equation*}
where for a $g\times g$ symmetric matrix $S$, the notation $S_0$ means the diagonal of $S$. 
This group is sometimes called the {\it Igusa group of level $\ell$} as in 
\cite[p.10]{Freitag91}
We will discuss in more details this symplectic theta group in Section \ref{secthesub}.
When $\ell=1$ we just write ${\rm Sp}_{g,\m Z}^{\theta}$ for ${\rm Sp}_{g,\m Z}^{\theta,1}$
\vs

In this paper we will mainly need these groups ${\rm Sp}_{g,\m Z}^{\ell}$ 
and ${\rm Sp}_{g,\m Z}^{\theta,\ell}$ for $\ell=2$.
Indeed  when $\ell=2$, the modular variety $X^{\theta,\ell}_g$ of {\it theta structures of level $\ell$} defined by
$X^{\theta,\ell}_g:={\rm Sp}_{g,\m Z}^{\theta,\ell}\backslash \mc H_g$
 will play an important role in this paper.

%22
\subsection{Critical values and theta functions}
\label{seccrithe}

We now recall the definition of the Riemann theta function:
$$\textstyle
\theta_\tau (z)=\theta(z,\tau):=\sum\limits_{m\in \m Z^{^g}}e^{i\pi{}^tm\tau m}e^{2i\pi {}^tmz},
\;\;\mbox{\rm for $z\in \m C^g$ and $\tau\in \mc H_g$.}
$$
This function is a holomorphic function of $z$ 
which is $\m Z^g$-periodic.  One has
$\theta_\tau(z+q)=\theta_\tau(z)$ for all $q$ in $\m Z^g$.
We can now explain our construction of $\la$-critical functions.
The construction involves 
a matrix ${\bf d}$  with  integer coefficients
and $\det({\bf d})\neq 0$,
and its 
associate group $G_{\bf d}:={\bf d}^{-1}\m Z^g/\m Z^g$
whose order $|G_{\bf d}|$ is equal to $|{\rm det}(\bf d)|$. 
Very often, we will choose ${\bf d}={\rm diag}(d_1,...,d_g)$ where each coefficient is positive and divides 
the next one: $d_1|d_2|\cdots|d_g$, 
Note that any finite abelian group is isomorphic to a unique group 
$G_{\bf d}$ with such a diagonal matrix ${\bf d}$.

\bd
\label{defrapirb}  Let $\tau\in \mc H_g$ and 
${\bf d}\in \mc M(g,\m Z)$ with $\det({\bf d})$ odd.\\
We will say that the function 
$\theta_{\tau}$ is $(\la,{\bf d})$-critical 
if, for all $z$ in $\m C^g$, the function
$
f_{z,\tau}:\ell\mapsto \theta(z\!+\!\ell,\tau)
$ 
is $\la$-critical on the group $G_{\bf d}:={\bf d}^{-1}\m Z^g/\m Z^g$.
\ed

This means that, for all $z$ in $\m C^g$,
$$\textstyle
\sum\limits_{\ell\,\in\,G_{\bf d}} \theta(z+\ell,\tau)\,\theta(z-\ell,\tau)
\;=\;\la\,\theta(z,\tau)^2.
$$

\br 
\label{remevencri}
In particular the function $f_{0,\tau}: \ell\mapsto \theta(\ell,\tau)$ is 
a $\la$-critical function on $G_{\bf d}$ which is even,
that is  $f_{0,\tau}(-\ell)=f_{0,\tau}(\ell)$ for all $\ell$ in $G_{\bf d}$.
\er

Note, when ${\bf d}={\rm diag}(d_1,...,d_g)$ as above,  that the group $G_{\bf d}$ has order $|G_{\bf d}|=d_1\cdots d_g$, 
and that this group  $G_{\bf d}$ is cyclic of order $d$ if and only if 
$1=d_1=\cdots =d_{g-1}<d_g=d$.

Here is our construction of critical values seen from the point of view of theta functions

\bt
\label{thmmaithe}
 Let $\tau\in \mc H_g$ and 
${\bf d}\in \mc M(g,\m Z)$ with $\det({\bf d})$ odd.\\
Assume that there exists 
$
\si\!=\!\mbox{\scriptsize 
$\left(\!\begin{array}{cc} \al&\be\\  
\ga&\de\end{array}\!\right)$} \in {\rm Sp}_{g,\m Z}^{\theta,2}
$ 
such that $\si\tau={}^t{\bf d}\tau{\bf d}$.\\
$a)$ Then there exists $\la\in \m C$ such that the function 
$\theta_\tau$ is $(\la,{\bf d})$-critical.\\
$b)$ One has $\la=\ka\, {\rm det}_{\m C}(\ga\tau+\de)^{1/2}|G_{\bf d}|$,
where $\ka^8=1$.
\et

Note that replacing $\si$ by $-\si$ does not change the assumptions
but changes the sign of the determinant when $g$ is odd.

In $a)$ the converse is true 
for $g\leq 3$ and is expected to be true for all $g$:
if there exists $\la$ such that 
the function $\theta_\tau$ is $(\la,{\bf d})$-critical
then there should exist $\si\in {\rm Sp}_{g,\m Z}^{\theta,2} $ 
such that $\si\tau={}^t{\bf d}\tau{\bf d}$. This will follow from 
Remark \ref{remph2phl}. 

A more precise formula for $\la$ will be given as Formula \eqref{eqnlajsitau}.  
Notice that it is easy to determine the $8^{\rm th}$ root of unity $\ka$ up to sign
without using \eqref{eqnlajsitau} by using instead
Proposition \ref{procrival} which says that $\la$ is an algebraic integer satisfying $\la\equiv 1$ mod $2$. Indeed, the only $8^{\rm th}$ roots of unity which are equal to $1$ mod $2$ are $\pm 1$.

%23
\subsection{Theta functions with characteristic}
\label{secthecar}

For the proof of Theorem \ref{thmmaithe} 
we will need to introduce the Riemann theta functions with characteristic (see \cite{BiLa04}). 
We will also need three classical formulas satisfied by these functions,
the ``addition formula'', the ``isogeny formula'',  and the ``transformation formula''. We will only need special cases of these 
formulas that we state below.
\vs 

The {\it theta functions with characteristic} $a$, $b$ in $\m C^g$,
are defined by, for $z\in \m C^g$ and $\tau\in \mc H_g$,
\begin{eqnarray*}
\theta
\mbox{\scriptsize 
$\left[\!\!\!
\begin{array}{c} a \\b\end{array}\!\!\!
\right]$} 
(z,\tau)
&:=&\textstyle
\sum\limits_{m\in \m Z^{^g}}e^{i\pi{}^t(m+a)\tau (m+a)}e^{2i\pi {}^t(m+a)(z+b)}.
\end{eqnarray*}
Note that these functions satisfy the following periodicity 
when translating the characteristic 
by elements $m$, $n$ in $\m Z^{g}$, 
\begin{eqnarray}
\label{eqnthetra}
\theta
\mbox{\scriptsize 
$\left[\!\!\!
\begin{array}{c} a+m \\b+n\end{array}\!\!\!
\right]$} 
(z,\tau)
&=&\textstyle
e^{2i\pi{}^ta n}\;
\theta
\mbox{\scriptsize 
$\left[\!\!\!
\begin{array}{c} a \\b\end{array}\!\!\!
\right]$} 
(z,\tau).
\end{eqnarray}
In this paper, we will mainly use the following special cases
of theta functions with characteristics. 
For $\xi\in \m Z^g/2\m Z^g$, seen as a subset of $\m Z^g$, we define
\begin{eqnarray}
\label{eqnthexiz}
\theta_{[\xi]}(z,\tau)=
\theta
\mbox{\scriptsize 
$\left[\!\!\!
\begin{array}{c} \xi/2 \\0\end{array}\!\!\!
\right]$} 
(2z,2\tau)
&:=&\textstyle
\sum\limits_{m\in \xi}e^{i\pi{}^tm\frac{\tau}{2} m}e^{2i\pi {}^tmz}.
\end{eqnarray}

Note that one has the equalities:
\begin{eqnarray*}
\label{eqnt0t1t2}\textstyle
\theta_{[0]}(z,\tau)=\theta(2z,2\tau)
&\;{\rm and}\;&\textstyle
\sum\limits_{\xi\in \,\m Z^{^g}\!/2\m Z^{^g}}\theta_{[\xi]}(z,\tau)=\theta(z,\tau/2).
\end{eqnarray*}

Here is the addition formula that we need.

\bl 
\label{lemaddfor}
For all $z,w$ in $\m C^g$, $\tau\in \mc H_g$, one has
\begin{eqnarray}
\label{eqnaddfor}
\theta(z+w,\tau)\,\theta(z-w,\tau)&=&
\textstyle\sum\limits_{\xi\in \m Z^{^g}/2\m Z^{^g}}
\theta_{[\xi]}(w,\tau )\, \theta_{[\xi]}(z,\tau ).
\end{eqnarray}
\el 

\begin{proof} 
Just write the left-hand side $LHS$ as a double sum over $m$, $n$ in $\m Z^g$
and split this double sum according to the class 
$\xi\in \m Z^g/2\m Z^g$ in which $m\!-\! n$ lives, 
and note that one has the equivalence: \;\; 
$m-n\in \xi\Longleftrightarrow m+n\in \xi$.
Use then a change of variable $p:=m-n$ and $q:=m+n$, 
this gives
\begin{eqnarray*}
LHS&=&\textstyle
\sum\limits_{\xi\in \m Z^{^g}/2\m Z^{^g}}\sum\limits_{p \in \xi}\sum\limits_{q \in \xi}
e^{i\pi{}^tp\frac{\tau}{2} p}e^{i\pi{}^tq\frac{\tau}{2} q}
e^{2i\pi {}^tpw}e^{2i\pi {}^tqz}\\
&=&\textstyle
\sum\limits_{\xi\in \m Z^{^g}/2\m Z^{^g}} \theta_{[\xi]}(w,\tau )\, \theta_{[\xi]}(z,\tau ),
\end{eqnarray*}
as required. \end{proof}

The second  formula is a simple but useful isogeny formula.

\bl 
\label{lemisofor}  
 Let $\tau\in \mc H_g$ and 
${\bf d}\in \mc M(g,\m Z)$ with $\det({\bf d})$ odd.\\ 
Set $G_{\bf d}:={\bf d}^{-1}\m Z^g/\m Z^g$. Then for all $\xi\in \m Z^g/2\m Z^g$, one has
\begin{eqnarray*}
\label{eqnisofor}\textstyle
\sum\limits_{\ell \in G_{\bf d}}\theta_{[\xi]}(\ell ,\tau)
&=&|G_{\bf d}|\;\theta_{[\xi]}(0,{}^t{\bf d}\tau{\bf d}).
\end{eqnarray*}
\el

\begin{proof} Just write the left-hand side $LHS$ as a double sum over $m$ in $\m Z^g$
and $\ell $ in $G_{\bf d}$ and notice that 
$\sum_{\ell \in G_{\bf d}}e^{2i\pi {}^tm\ell }$ is equal to 
the order
$|G_{\bf d}|$ of the group $G_{\bf d}$ 
when $m$ belongs to ${\bf d}\m Z^g$
and is equal to $0$ otherwise.
Hence 
\begin{eqnarray*}
LHS
&=&
\textstyle
|G_{\bf d}|
\sum\limits_{m \in {\bf d}\m Z^{^g}\cap \xi}
e^{i\pi{}^tm\frac{\tau}{2} m}\\
&=& |G_{\bf d}|\;\theta_{[\xi]}(0,{}^t{\bf d}\tau{\bf d}).
\end{eqnarray*}
In the last equality we used ${\rm det}({\bf d})$  odd by writing $m={\bf d}p$ with $p\in \xi$.
\end{proof}

%24
\subsection{The theta cocycle}
\label{secthecoc}

The last formula is a transformation formula for the 
theta functions with characteristic.
It deals with an element 
$\si
=\mbox{
\scriptsize 
	$\left(\!
	\begin{array}{cc} \al&\be   \\
	\ga&\de
	\end{array}\!
	\right)$} \in {\rm Sp}(g,\m Z).
$
This formula is particularly simple when $\si$ belongs to the theta group
and when it is expressed with the {\it modified} theta function
\begin{eqnarray}
\label{eqnmodthe}
\widetilde{\theta}
\mbox{\scriptsize 
$\left[\!\!\!
\begin{array}{c} a \\b\end{array}\!\!\!
\right]$} 
(z,\tau)
&=& e^{-i\pi {}^t\! a(z+b)}\;
\theta
\mbox{\scriptsize 
$\left[\!\!\!
\begin{array}{c} a \\b\end{array}\!\!\!
\right]$} 
(z,\tau).
\end{eqnarray}
Note that there is no {\it modification} when $z=b=0$.

\bl 
\label{lemtrafor}  Let $\tau\in \mc H_g$  and  $\si\in {\rm Sp}_{g,\m Z}^{\theta}$. Then, for $a$, $b$ in 
$\m C^g$, one has
\begin{eqnarray}
\label{eqntrafor}
\widetilde{\theta}
\mbox{\scriptsize 
$\left[\!\!\!
\begin{array}{c} \de a \!-\!\ga b\\-\be a\! +\! \al b\end{array}\!\!\!
\right]$} 
(0,\si\tau)
&=& j(\si,\tau)\;\,
\widetilde{\theta}
\mbox{\scriptsize 
$\left[\!\!\!
\begin{array}{c} a \\ b\end{array}\!\!\!
\right]$} 
(0,\tau),
\;\;{\rm where}\\
\label{eqntraforbis}
j(\si,\tau)
&=&
\ka(\si)\;
{\rm det}_{\m C}(\ga\tau+\de)^{\frac12}
\end{eqnarray}
\el
In this formula, $j(\si,\tau)$ is a cocycle on ${\rm Sp}_{g,\m Z}^{\theta}\times \mc H_g$ called the {\it theta cocycle} which is analytic in $\tau$:
one has $j(\si_1\si_2,\tau)=j(\si_1,\si_2\tau)\,j(\si_2,\tau)$.
The constant $\ka(\si)$ is a eigth root of unity,
$\ka(\si)^8=1$,
that depends only on $\si$ for a continuous choice of 
the square root ${\rm det}_{\m C}(\ga\tau+\de)^{\frac12}$
of the complex number 
${\rm det}_{\m C}(\ga\tau+\de)$.
The precise value of $j(\si,\tau)$
will  be explained in Section \ref{secsigtrafor}.

\begin{proof} This is \cite[Th. 5.7]{Freitag91} or
\cite[Section 8.6 p.231]{BiLa04}.
The textbooks  \cite{Eichler66}, \cite{Mum83},  \cite{Kempf91}, or \cite{Poli03} also discuss this transformation formula. We recall the strategy of proof. 
One proves a more involved transformation
formula, \cite[Prop. 5.6, 5.7]{Freitag91}, 
for 
$
\theta\mbox{\scriptsize 
$\left[\!\!\!
\begin{array}{c} a \\b\end{array}\!\!\!
\right]$}
$ 
valid for all $\si$ in ${\rm Sp}(g,\m Z)$, by checking it on generators
of ${\rm Sp}(g,\m Z)$. 
The first generators are translations by
an integral symmetric matrix $\be$,
\begin{equation}
\theta
\mbox{\scriptsize 
$\left[\!\!\!
\begin{array}{c} a \\
-\be a\!+\!b\!+\!\be_0/2 \end{array}\!\!\!
\right]$} 
(0,\tau+\be)
= e^{i\pi {}^t\! a(-\be a +\be_0)}\;
\theta
\mbox{\scriptsize 
$\left[\!\!\!
\begin{array}{c} a \\ b\end{array}\!\!\!
\right]$} 
(0,\tau),
\end{equation}
where  $\be_0$ is the diagonal of $\be$ seen as an element of $\m Z^g$.

The formula for the second generator is the Poisson formula,
\begin{equation}
\theta
\mbox{\scriptsize 
$\left[\!\!\!
\begin{array}{c} -b \\ a\end{array}\!\!\!
\right]$} 
(0,-\tau^{-1})
= 
{\rm det}_{\m C}(-i\tau)^{\frac12}\;
e^{-2i\pi {}^t\! a b}\;
\theta
\mbox{\scriptsize 
$\left[\!\!\!
\begin{array}{c}  a \\b\end{array}\!\!\!
\right]$} 
(0,\tau),
\end{equation}
where the square root is
defined by holomorphic continuation in $\tau$ 
with the constraint that when $\tau=i\mathds{1}$ it is equal to $1$. 
One uses then the fact that the map $(\si,\tau)\mapsto \det_{\m C}(\si\tau+\de)$ 
is a cocycle on ${\rm Sp}(g,\m Z)\times \mc H_g$.  
\end{proof}

%\br In the  book \cite[App. of Chap. 1 p. 41]{Eichler66} and many  articles the translation formula is given for a function which is called $\theta(\tau,u,v)$. The difference is only in the notation since$\theta(\tau, u,v)=\widetilde{\theta}\mbox{\scriptsize $\left[\!\!\!\begin{array}{c} v \\ -u\end{array}\!\!\!\right]$}(0,\tau) $.\er 

%25
\subsection{The condition on theta constant}
\label{secconthe}

The first step in the proof of Theorem \ref{thmmaithe} is 
the following criterion on $\la, \tau, {\bf d}$ which ensures that the function $\theta_\tau$ 
is $(\la,{\bf d})$-critical. This criterion is a  relation between  
``theta constants'', i.e. theta functions evaluated at $z=0$.

\bl
\label{lemratequ}
Let $\tau\in \mc H_g$, $\la\in \m C$ and ${\bf d}\in {\mc M}(g,\m Z)$
with ${\rm det}({\bf d})$ odd.
The function $\theta_\tau$ is $(\la,{\bf d})$-critical
if and only if 
the ratios 
$$ 
|G_{\bf d}|\; \frac{\theta_{[\xi]}(0,{}^t{\bf d}\tau{\bf d})}{\theta_{[\xi]}(0,\tau)}
$$
do not depend on $\xi\in \m Z^g/2\m Z^g$ 
and are equal to $\la$.
\el

It might happen that for some $\xi$, 
the denominator $\theta_{[\xi]}(0,\tau)$ is zero. In this case, the condition means that 
$\theta_{[\xi]}(0,{}^t{\bf d}\tau{\bf d})$ has to be zero too.

\begin{proof} For $w$ in $\m C^g$ we introduce the function on $\m C^g$
$$z\mapsto F_w(z)=F_w(z,\tau):=\theta(z+w,\tau)\,\theta(z-w,\tau).$$
We want to know when the two functions 
$\sum_{\ell\in G_{\bf d}} F_{\ell}$ and $F_0=\theta^2$ are proportional. 
The key point  in the proof is that all these functions 
$F_w$ live in the same finite dimensional vector space 
and that this vector space has a very convenient basis:
$(\theta_{[\xi]})_{\xi\in \,\m Z^{^g}\!/2\m Z^{^g}}$. 
We only have to express that the coefficients of our two functions in this basis are proportional. 
These coefficients are given by the following calculation in which we apply successively the addition formula and the isogeny formula,
\begin{eqnarray*}\textstyle
\sum\limits_{\ell\in G_{\bf d}} F_{\ell }(z,\tau )
&=& \textstyle
\sum\limits_{\ell\in G_{\bf d}} \;
\sum\limits_{\xi\in \m Z^{^g}\!/2\m Z^{^g}}
\theta_{[\xi]}(\ell ,\tau)\; \theta_{[\xi]}(z,\tau)\\
&=& \textstyle
|G_{\bf d}|\;\sum\limits_{\xi\in \m Z^{^g}\!/2\m Z^{^g}}
\theta_{[\xi]}(0,{}^t{\bf d}\tau{\bf d})\; \theta_{[\xi]}(z,\tau) 
 \;\;\; {\rm and}
\end{eqnarray*}
\begin{eqnarray*}
\theta(z,\tau)^2&=&\textstyle
\sum\limits_{\xi\in \m Z^{^g}\!/2\m Z^{^g}}
\theta_{[\xi]}(0,\tau)\; \theta_{[\xi]}(z,\tau) .\hspace*{3em}
\end{eqnarray*}
These two functions are proportional with  proportionality factor $\la$
if and only if one has,
\begin{eqnarray}
\label{eqnlat0t1} 
\la
= 
|G_{\bf d}|\; \frac{\theta_{[\xi]}(0,{}^t{\bf d}\tau{\bf d})}{\theta_{[\xi]}(0,\tau)},
&\;\;& 
\mbox{\rm for all $\xi$ in $\m Z^g/2\m Z^g$}.
\end{eqnarray}
This is the criterion we were looking for.
\end{proof}

\br 
Note that, for every $\tau$ in $\mc H_g$ and  $z$ in $\m C^g$, 
there exists $\xi$ in $\m Z^{^g}\!/2\m Z^{^g}$ such that
$\theta_{[\xi]}(z,\tau)\neq 0$, see for instance \cite[Sections 3.4-3.8]{Beau13}.
\er

%26
\subsection{The moduli variety of theta structures}
\label{secembmod}

In order to exploit the criterion \eqref{eqnlat0t1},
the following corollary of Lemma \ref{lemtrafor} will be very useful.
It is due to Igusa (see \cite[Lemma 9.2 p.239]{BiLa04}).
\bc
\label{cortrafor}  
When $\si\in{\rm Sp}^{\theta,2}_{g,\m Z}$ and $\tau\in\mc H_g$, for all $\xi\in \m Z^g/2\m Z^g$, one has
\begin{eqnarray}
\label{eqntrafor2}
\frac{\theta_{[\xi]}(0,\si\tau)}{\theta_{[0]}(0,\si\tau)}&=& \frac{\theta_{[\xi]}(0,\tau)}{\theta_{[0]}(0,\tau)}\, .
\end{eqnarray}
\ec

\begin{proof} It will be useful to recall the proof of this corollary.\\
Introduce  
$\si'
:=\mbox{
\scriptsize 
$\left(\!
\begin{array}{cc} \al&\! 2\be\!   \\
\!\ga/2\!&\de
\end{array}\!
\right)$}
$ 
so that 
$
\si'(2\tau)=2\si\tau .
$
Since the matrix $\si$ is in ${\rm Sp}_{g,\m Z}^{\theta,2}$, the matrix 
$\si'$ is in ${\rm Sp}_{g,\m Z}^\theta$. 
We claim that, for all $\xi\in \m Z^g/2\m Z^g$, 
\begin{eqnarray}
\label{eqnthxisi}
\theta_{[\xi]}(0,\si\tau)
&=&
j(\si',2\tau)\;
\theta_{[\xi]}(0,\tau).
\end{eqnarray}
Indeed, we  compute remembering that,
by assumption,
the matrices $(\de -1)/2$, $\be/2$, the vector $\xi$ and 
the scalar ${}^t\xi{}^t\de\be\xi/4$  are all integral,
\begin{align*}
\theta_{[\xi]}(0,\si\tau)
&=
\theta
\mbox{\scriptsize 
$\left[\!\!\!
\begin{array}{c} \xi/2 \\ 0\end{array}\!\!\!
\right]$} 
(0,\si'(2\tau))
\hspace{4.2em}
\mbox{\rm by Definition \eqref{eqnthexiz},}\\
&=
\theta
\mbox{\scriptsize 
$\left[\!\!\!
\begin{array}{c} \de\xi/2 \\ -\be\xi\end{array}\!\!\!
\right]$} 
(0,\si'(2\tau))\hspace{4em}
\mbox{\rm by Property \eqref{eqnthetra},}\\
&=
\widetilde{\theta}
\mbox{\scriptsize 
$\left[\!\!\!
\begin{array}{c} \de\xi/2 \\ -\be\xi\end{array}\!\!\!
\right]$} 
(0,\si'(2\tau))
\hspace{4em}
\mbox{\rm by Definition \eqref{eqnmodthe}.}\\
\end{align*}
We now apply 
the transformation formula in Lemma \ref{lemtrafor}
to the pair $(\si',2\tau)$,
\begin{eqnarray*}
\theta_{[\xi]}(0,\si\tau)
& =& 
j(\si',2\tau)\;
\widetilde{\theta}
\mbox{\scriptsize 
$\left[\!\!\!
\begin{array}{c} \xi/2 \\ 0\end{array}\!\!\!
\right]$} 
(0,2\tau)\\
& =& 
j(\si',2\tau)\;
\theta
\mbox{\scriptsize 
$\left[\!\!\!
\begin{array}{c} \xi/2 \\ 0\end{array}\!\!\!
\right]$} 
(0,2\tau)\\
&=&
j(\si',2\tau)\;
\theta_{[\xi]}(0,\tau).
\end{eqnarray*}
This proves that the ratio
$\displaystyle\frac{\theta_{[\xi]}(0,\si\tau)}{\theta_{[\xi]}(0,\tau)}$
does not depend on $\xi$ as required.
\end{proof}

\noindent
Corollary \ref{cortrafor}  tells us that 
the map $\Phi_2$ given in homogeneous coordinates by
\begin{eqnarray*}
\Phi_2: \; \mc H_g
&\longrightarrow&
\m P(\m C^{2^g})\\
\tau&\mapsto & [\cdots ,\theta_{[\xi]}(0,\tau),\cdots]
\end{eqnarray*}
induces a well defined holomorphic map
$$
\ph_2: X^{\theta,2}_g\longrightarrow
\m P(\m C^{2^g}).
$$

\br 
\label{remph2phl}
These maps $\ph_2$ and their analogs 
$
\ph_\ell: X_g^{\theta,\ell}\longrightarrow
\m P(\m C^{\ell^g}),
$ 
for $\ell\geq 2$, have a long history, as the quasi-projective realizations of the
moduli varieties of theta structures of level $\ell$. 
We will not need here the precise definition of these structures. 
But here are some comments that  relate our computation to the existing litterature.

For $\ell\geq 3$, 
according to successive works of Igusa, Mumford and Salvati Manni,
these maps $\ph_\ell$ are proven to be  embeddings,
see \cite{Igusa72},   \cite{Mum74}, 
and \cite[Section 8.10]{BiLa04}.

For $\ell=2$, the situation is more delicate and has been studied in detail by Salvati Manni:\\ 
$\star$ He proves that the map $\ph_2$ is  generically injective, see  \cite[Prop. 1]{Salva94a}.\\
$\star$ He conjectures that the map $\ph_2$ is injective,
see  \cite[Theorem 3]{Salva94a} where a tentative proof is given. See also
\cite[Theorem 3.6]{Grush09} and 
\cite[Section 2]{OuraSalva08} for more comments on this question.\\
$\star$ For $g\leq 3$, the map $\ph_2$ is indeed injective, see \cite{Runge98}.\\
$\star$ For $g\geq 4$, the map $\ph_2$ is not a biholomorphism.
\er

\begin{proof}[Proof of Theorem \ref{thmmaithe}]
$(a)$ 
Our assumptions and Corollary \ref{cortrafor} tell us that
$$
\Phi_2({}^t{\bf d}\tau{\bf d})
=
\Phi_2(\tau).
$$
This equality is nothing but the criterion of Lemma \ref{lemratequ}. 
Therefore the function $\th_\tau$ is $(\la,{\bf d})$-critical for a critical value $\la$.

$(b)$ To compute the critical value $\la$, 
we use again Lemma \ref{lemratequ}  combined 
with Corollary \ref{cortrafor}, and more precisely with Formulas \eqref{eqnthxisi}
and \eqref{eqntraforbis}. We obtain, for all $\xi$ in $\m Z^g/2\m Z^g$,
\begin{eqnarray}
\label{eqnlajisitau}
\nonumber
\la&=&|G_{\bf d}|\;\frac{\theta_{[\xi]}(0,{}^t{\bf d} \tau{\bf d})}{\theta_{[\xi]}(0,\tau)}=
|G_{\bf d}|\;
\frac{\theta_{[\xi]}(0,\si\tau)}{\theta_{[\xi]}(0,\tau)} ,\\
\label{eqnlaepde}
\la &=&
j(\si',2\tau)\,|G_{\bf d}|\;
\; =\;
\ka(\si')\;{\rm det}_{\m C}(\ga\tau+\de)^{1/2}\,|G_{\bf d}|
\; , 
\end{eqnarray}
where the matrix 
$\si':=\mbox{\scriptsize 
$\left(\!
\begin{array}{cc} \al&\! 2\be\!   \\ \!\ga/2\!&\de\end{array}\!
\right)$}$ 
belongs to ${\rm Sp}_{g,\m Z}^\theta$ and $\ka(\si')^8=1.$
\end{proof}

\br 
\label{remgaoltau}  One has $|\la|=|G_{\bf d}|^{1/2}$ and Equation \eqref{eqnlaepde} 
can be written as
\begin{eqnarray}
\label{eqnlajisitaubis}
\la=\ka(\si')\,
{\rm det}_{\m C}(\ga\ol{\tau}+\de)^{-1/2}.
\end{eqnarray}

Indeed, let $\La_\tau$ be the lattice $\La_\tau:=\tau\m Z^g\oplus \m Z^g$ of $\m C^g$.
Since $\si\tau={}^t{\bf d}\tau {\bf d}$, 
the complex matrix
$M:={}^t(\ga\tau+\de)^{-1}\in {\rm GL}(\m C^g)$ satisfies the 
equality between $g\!\times\! 2g$ complex matrices:
$$
M \,(\tau\;\, {\bf 1}_g)\, {}^t\si 
=
({}^t{\bf d}\tau {\bf d}\;\;{\bf 1}_g).
$$ 
This equality implies that 
$$
M(\La_\tau)=\La_{\,{}^t\!{\bf d}\tau{\bf d}}.
$$ 
Comparing the covolume of these lattices, one gets
the equality 
$$
|det_{\m C}(M)|=|\det({\bf d})|=|G_{\bf d}|,
$$
which implies \eqref{eqnlajisitaubis}. 
\er

%30
\section{The symplectic group}
\label{secsymgro}

In this chapter we first give in Corollary \ref{cormaisym} a reformulation of theorem \ref{thmmaithe}
with a more precise formula for the critical value.
We then explain various tools for studying the symplectic group
like the symplectic adapted basis theorem
and the Cayley transform that will be useful in the following chapters.

%31
\subsection{Critical values and symplectic transformation}
\label{seccrisym}

We recall that a matrix $h\in {\rm Sp}(g,\m R)$ is said to be 
{\it elliptic} if it is diagonalizable over $\m C$ with all eigenvalues of modulus $1$. It is equivalent to say that $h$ has a fixed point in $\mc H_g$.

\bc
\label{cormaisym}
 Let $\tau\in \mc H_g$, let
${\bf d}\in \mc M(g,\m Z)$ with $\det({\bf d})$ odd,
let\\ $
\si\!=\!\mbox{\scriptsize 
$\left(\!\begin{array}{cc} \al&\be\\  
\ga&\de\end{array}\!\right)$} \in {\rm Sp}_{g,\m Z}^{\theta,2}
$ 
such that the sympletic matrix $h:=
\mbox{\scriptsize 
$\left(\!\begin{array}{cc} {}^t{\bf d}^{-1}\al\!&\!{}^t{\bf d}^{-1}\be\\  
{\bf d}\ga&{\bf d}\de\end{array}\!\right)$}
$ is elliptic, let $\tau$ be a fixed point of $h$ in $\mc H_g$
and let $\si':=\!\mbox{\scriptsize 
$\left(\!\begin{array}{cc} \al&\!2\be\!\\  
\!\ga/2\!&\de\end{array}\!\right)$}$.
Then 
\begin{equation}
\label{eqnlajsitau}
\la:= j(\si',2\tau)\, |\det({\bf d})|
\end{equation} is a critical value on the 
group ${\bf d}^{-1}\m Z^g/\m Z^g$.
\ec

\begin{proof}[Proof of Corollary \ref{cormaisym}]
This is a direct corollary of Theorem \ref{thmmaithe}
with formula \eqref{eqnlajisitau}. Indeed the condition $h\tau=\tau$
is equivalent to $\si\tau={}^t{\bf d}\tau{\bf d}$.
A $\la$-critical function can be chosen to be the restriction 
of the $\m Z^g$-periodic function $\th_\tau$ 
to the finite group 
$G_{\bf d}={\bf d}^{-1}\m Z^g/\m Z^g$.
\end{proof}

%32
\subsection{Theta subgroup of level $\ell$}
\label{secthesub}

In order to use efficiently 
Theorem \ref{thmmaithe}, we
recall equivalent definitions 
for the theta subgroup ${\rm Sp}_{g,\m Z}^{\theta,2}$ (see \cite[p.177-182]{DolOrt88} for more on this topic).

\bfa 
$(a)$ For $\ell\geq 1$, the group ${\rm Sp}_{g,\m Z}^{\ell}$ is a normal subgroup of  ${\rm Sp}(g,\m Z)$.\\
$(b)$ When $\ell$ is even,  one has the equalities
\begin{eqnarray*}
\label{eqnthestr2}
{\rm Sp}_{g,\m Z}^{\theta,\ell}&=& 
\{
\si\!=\!\mbox{\scriptsize 
$\left(\!\!\begin{array}{cc} \al\!&\be\\  
\ga\!&\de\end{array}\!\!\right)$} \in {\rm Sp}_{g,\m Z}^{\ell} \mid \,
\be_0\equiv \ga_0\equiv 0\; {\rm mod}\; 2\ell\}, \\
&=&
\{
\si\in {\rm Sp}_{g,\m Z}^{\ell}\mid \om(\si x,x)\equiv 0\;{\rm mod}\; 2\ell\;
\mbox{\rm for all $x$ in $\m Z^{2g}$}\}
\end{eqnarray*}  
$(c)$ When $\ell$ is even,  the group  ${\rm Sp}_{g,\m Z}^{\theta,\ell}$ is also a normal subgroup of ${\rm Sp}(g,\m Z)$.
\efa

\begin{proof}
This is classical, see \cite[Lemma 4 p.177]{Igusa72}. Here are some details.

$(a)$ ${\rm Sp}_{g,\m Z}^{\ell}$ is the kernel of the projection 
$\pi_\ell:{\rm Sp}(g,\m Z)\ra {\rm Sp}(g,\m Z/\ell\m Z)$. 

$(b)$ It is enough to check that these three sets have same
image in the quotient ${\rm Sp}(g,\m Z/2\ell\m Z)$. 
Since $\pi_{2\ell}$ is onto, and $\ell$ is even, the group
\begin{eqnarray*}
\label{eqngagell}
\pi_{2\ell}({\rm Sp}_{g,\m Z}^{\ell})&\simeq& 
\{
\si\!=\!\mbox{\scriptsize 
$\left(\!\!\begin{array}{cc} \mathds{1}+\ell a\!&\ell b\\  
\ell c\!&\mathds{1}+\ell d\end{array}\!\!\right)$}  \mid \,
b\equiv {}^tb\; ,\; c\equiv{}^tc\;{\rm and}\; d\equiv{}^ta\;  {\rm mod}\; 2\}
\end{eqnarray*}  
is abelian and is  a $\m F_2$-vector space $V_\ell\simeq \m F_2^{(2g+1)g}$. 
The image in $V_\ell$ of each of these three sets is the $\m F_2$-vector subspace 
$V^\theta_\ell$ of codimension 
$2g$ in $V_\ell$ given by the equation $b_0\equiv c_0\equiv 0 \;{\rm mod}\; 2$,
where, as before, $b_0$ and $c_0$ are the diagonals of the symmetric matrices $b$ and $c$.

$(c)$ The quotient ${\rm Sp}(g,\m Z)/{\rm Sp}_{g,\m Z}^{2}$ 
is the symplectic group ${\rm Sp}(g,\m F_2)$ 
over the finite field $\m F_2$, i.e. the stabilizer of the non-degenerate symmetric bilinear form $\om$ seen  on $\m F_2^{2g}$. One easily checks that the
action by conjugation of ${\rm Sp}(g,\m Z)$ on the  quotient $V_\ell\simeq {\rm Sp}_{g,\m Z}^{\ell}/{\rm Sp}_{g,\m Z}^{2\ell}$
preserves the vector subspace 
$V^\theta_\ell$. Hence the group ${\rm Sp}_{g,\m Z}^{\theta,\ell}$ is normal in ${\rm Sp}(g,\m Z)$.
\end{proof}

\br 
This also proves that the quotient $\widetilde{\rm Sp}(g,\m F_2):={\rm Sp}(g,\m Z)/{\rm Sp}_{g,\m Z}^{\theta,2}$ 
is an extension 
of ${\rm Sp}(g,\m F_2)$ by $\m F_2^{2g}$:
\begin{equation}
\label{eqnsptgf2}
1 \longrightarrow \m F_2^{2g} \longrightarrow 
\widetilde{\rm Sp}(g,\m F_2)\longrightarrow {\rm Sp}(g,\m F_2)\longrightarrow 1.
\end{equation}
\er

%33
\subsection{The symplectic adapted basis}
\label{secconcri}

In this section we discuss the structure of the {\it rational symplectic group}
${\rm Sp}(g,\m Q):={\rm GL}(2g,\m Q)\cap {\rm Sp}(g,\m R)$, and its relation with the  {\it integral symplectic group}
${\rm Sp}(g,\m Z)$. We also introduce the rational symplectic theta group 
${\rm Sp}^{\theta,2}_{g,\m Q}$ of level $2$.
\vs

The following proposition is a variation of the classical ``adapted basis theorem'' 
which takes into account the existence of a symplectic form.

\bp
\label{prohsidsi}
Let $h\in {\rm Sp}(g,\m Q)$. Then there exists $\si_1$ and $\si_2$ in 
${\rm Sp}(g,\m Z)$ and a diagonal matrix ${\bf d}={\rm diag}(d_1,\ldots,d_g)$ 
with $d_1|d_2|\ldots |d_g$ integral and
\begin{eqnarray*}
h&=&\si_1\;
\mbox{\small
$\left(\!
\begin{array}{cc} 
\!{ }^t{\bf d}^{-1}\!&{\bf 0}\\  
{\bf 0}&{\bf d}
\end{array}
\!\right)$}
\;\si_2.
\end{eqnarray*}
\ep

A proof of this proposition is given in \cite{RSG}. 
\vs

For $\ell\geq 1$, let $\m Z_{(\ell)}$ be the ring of rational numbers with denominator prime to $\ell$. 
We introduce
the {\it  rational congruence symplectic group of level $\ell$}
$$
{\rm Sp}_{g,\m Q}^{\ell}:= \{ h\in {\rm Sp}(g,\m Z_{(\ell)})\mid \si \equiv\mathds{1}_{2g}\;{\rm mod}\; \ell\},
$$
and the {\it rational symplectic theta group  of level $2$}
\begin{equation*}
\label{eqnthestrq}
{\rm Sp}_{g,\m Q}^{\theta,2}:= 
\{
h\!=\!\mbox{\scriptsize 
$\left(\!\!\begin{array}{cc} \al\!&\be\\  
\ga\!&\de\end{array}\!\!\right)$} \in {\rm Sp}_{g,\m Q}^{2} \mid \,
({}^t\al\ga)_0\equiv ({}^t\be\de)_0\equiv 0\; {\rm mod}\; 4\}, 
\end{equation*}
We will say that $h\in {\rm Sp}(g,\m Q)$  {\it preserves a theta structure 
of level $2$} if it belongs to ${\rm Sp}_{g,\m Q}^{\theta,2}$.
As in the integral case, the group ${\rm Sp}^{\theta,2}_{g,\m Q}$ is a normal subgroup of the group $ {\rm Sp}(g,\m Z_{(2)})$, and one has the inclusions
\begin{equation*}
\label{eqngagq2}
{\rm Sp}_{g,\m Q}^4
\subset 
{\rm Sp}^{\theta,2}_{g,\m Q}
\subset
{\rm Sp}_{g,\m Q}^{2}
\subset
{\rm Sp}(g,\m Z_{(2)}).
\end{equation*}
Indeed the reduction modulo $4$ of the group ${\rm Sp}^{\theta,2}_{g,\m Q}$ is the group  $\widetilde{\rm Sp}(g,\m F_2)$ which is a normal subgroup of the group 
${\rm Sp}(g,\m Z/4\m Z)\simeq {\rm Sp}(g,\m Z_{(2)})/{\rm Sp}_{g,\m Q}^4$.

\bl
\label{lemthestr}
Let $h\in {\rm Sp}(g,\m Z_{(2)})$ 
and write 
$ h=\si_1\;
\mbox{\small\scriptsize 
$\left(\!
\begin{array}{cc} 
\!{}^t{\bf d}^{-1}\!&{\bf 0}\\  
{\bf 0}&{\bf d}
\end{array}
\!\right)$}
\;\si_2
$
with ${\bf d}$ in $\mc M(g,\m Z)$ and both $\si_1$ and $\si_2$ in ${\rm Sp}(g,\m Z)$.
Then the following are equivalent:\\
$\star$ $h$ preserves a theta structure of level $2$, i.e. $h\in {\rm Sp}_{g,\m Q}^{\theta,2}$.\\
$\star$ ${\rm det}({\bf d})$ is odd and $\si_1\si_2\in {\rm Sp}_{g,\m Z}^{\theta,2}$.\\
$\star$ ${\rm det}({\bf d})$ is odd and $\si_2\si_1\in {\rm Sp}_{g,\m Z}^{\theta,2}$.
\el

\begin{proof}[Proof of Lemma \ref{lemthestr}] 
We  first note 
the equivalence 
$$h\in {\rm Sp}_{g,\m Q}^{2}\;
\;\Longleftrightarrow
\;\;{\rm det}({\bf d})
\; \mbox{\rm is odd\,  and}\; \si_2\si_1\in {\rm Sp}_{g,\m Q}^2.
$$
One conclude by noticing that ${\rm Sp}^{\theta,2}_{g,\m Q}$ is a normal subgroup of  $ {\rm Sp}(g,\m Z_{(2)})$.
\end{proof}

%34
\subsection{The Cayley transform}
\label{seccaytra}

The Cayley transform is a convenient tool that allows us to construct 
symplectic transformations $\si$ with rational coefficients
(Lemma \ref{lemcayley}),
and to recognize those that preserve a theta structure of level $2$ (Lemma \ref{lemcaythe}).
\vs 

Let $\g s\g p(g,\m Q)$ be the ``Lie algebra'' of ${\rm Sp}(g,\m Q)$,
that is,
\begin{eqnarray*}
\label{eqnspgq0}
\g s\g p(g,\m Q)
&=&\{
X\! \in\! \mc M(2g,\m Q) \mid \,
\om(Xx,y)+\om(x,Xy)=0\;\mbox{\rm  for all $x,y$ in $\m Q^{2g}$}\}.\\
&=&
\{
X\!=\!\mbox{\scriptsize 
$\left(\!\!\begin{array}{cc}  a\!& b\\  
 c\!& d\end{array}\!\!\right)$} \in \mc M(2g,\m Q) \mid \,
b= {}^tb\; ,\; c={}^tc\;,\; d=-{}^ta\;  \},\\
&=& 
\{
X\!=\!\mbox{\scriptsize 
$\left(\!\!\begin{array}{cc}  a\!& b\\  
 c\!& d\end{array}\!\!\right)$} \in \mc M(2g,\m Q) \mid \,
X=-X^*\;  \},\;\; \mbox{where}
\end{eqnarray*}
$$X^*\!=\! J\, { }^tXJ^{-1}=\mbox{\scriptsize 
$\left(\!\!\begin{array}{cc} {}^t d\!& - {}^t b\\  
 - {}^t c\!&  {}^t a\end{array}\!\!\right)$}.
$$
The Cayley transform will give a bijection between the following two 
Zariski open subsets 
\begin{eqnarray*}
\label{eqnspgqs}
\g s\g p^\bullet_{g,\m Q}&:=& 
\{
X\in \g s\g p(g,\m Q) \mid \,
{\rm det}(\mathds{1}-X)\neq 0   \}\;{\rm and}\\
{\rm Sp}^\bullet_{g,\m Q}&:=&\{ h\in {\rm Sp}(g,\m Q) \mid \,
{\rm det}(\mathds{1}+h)\neq 0 \}.
\end{eqnarray*}  

Here is the construction of the Cayley transform.

\bl
\label{lemcayley} 
The map $X\mapsto C(X)=h:=(\mathds{1}+X)(\mathds{1}-X)^{-1}$ is a bijection from 
$\g s\g p^\bullet_{g,\m Q}$ to ${\rm Sp}^\bullet_{g,\m Q}$ with inverse map  $h\mapsto X:=-(\mathds{1}-h)(\mathds{1}+h)^{-1}$.
\el

\begin{proof} 
The equality $h=C(X)$ can be rewritten as 
$
\frac12(\mathds{1}+h)(\mathds{1}-X)=\mathds{1}
$.

It implies that both matrices $\mathds{1}-X$ and $\mathds{1}+h$ are invertible. 

For such 
a pair $(X,h)$ in $\mc M(2g,\m Q)$, one has the equivalences:\\
$X\in \g s\g p(g,\m Q)$ $\Longleftrightarrow$  
$\om(Xx,y)+\om(x,Xy)=0$ for all $x,y$ in $\m Q^{2g}$ $\Longleftrightarrow$\\ 
$\om((\mathds{1}+X)x,(\mathds{1}+X)y)=\om((\mathds{1}-X)x,(\mathds{1}-X)y)$ for all $x,y$ in $\m Q^{2g}$
 $\Longleftrightarrow$\\ 
$\om(hx',hy')=\om(x',y')$ for all $x',y'$ in $\m Q^{2g}$
 $\Longleftrightarrow$
 $h\in {\rm Sp}(g,\m Q)$.
\end{proof}
Let\;\; $\g s\g p_{g,\m Q}^2=
\g s\g p(g,\m Z_{(2)})$ and
\begin{align*}
\g s\g p^{\theta,2}_{g,\m Q}&=
\{
X\!=\!\mbox{\scriptsize 
$\left(\!\!\begin{array}{cc}  a\!& b\\  
 c\!& d\end{array}\!\!\right)$} \in \g s\g p_{g,\m Q}^2 \mid \,
b_0\equiv c_0\equiv 0 \;{\rm mod}\; 2 \},\\
&=
\{
X\! \in\!{\g s\g p}_{g,\m Q}^2 \mid \,
\om(Xx,x)\equiv 0\;{\rm mod}\;2\; \mbox{\rm  for all $x$ in $\m Z_{(2)}^{2g}$}\}.
\end{align*}
where again $b_0$ and $c_0$ are the diagonals of the symmetric matrices $b$ and $c$.

The equivalence between these two definitions follows from the fact that, for $X$ in ${\g s\g p}_{g,\m Q}^2$, 
the $\m Z_{(2)}$-valued bilinear form $\om(Xx,y)$ on  $\m Z_{(2)}^{2g}$ is symmetric.

By construction one has the inclusions
$$
2\g s\g p_{g,\m Q}^2\subset 
\g s\g p^{\theta,2}_{g,\m Q}\subset\g s\g p_{g,\m Q}^2.
$$
Here is  a useful  interpretation of this intermediate 
$\m Z_{(2)}$-module
$\g s\g p^{\theta,2}_{g,\m Q}$  using the involution $X\mapsto X^*$.

\bl
\label{lemspgq2}
One has $\g s\g p^{\theta,2}_{g,\m Q}=
\{X=Y-Y^*\mid Y\in \mc M(2g,\m Z_{(2)})\}$
\el

\begin{proof} It is possible to write an element
$X\!=\!\mbox{\scriptsize 
$\left(\!\!\begin{array}{cc}  a\!& b\\  
 c\!& d\end{array}\!\!\right)$} \in \g s\g p_{g,\m Q}^2$
 as a sum 
 $X\!=\!\mbox{\scriptsize 
 $\left(\!\!\begin{array}{cc}  \al-{}^t\de\!& \be+{}^t\be\\  
 \ga+{}^t\ga\!& \de-{}^t\al\end{array}\!\!\right)$}$ with $\al$, $\be$, $\ga$, $\de$ in $\mc M(g,\m Z_{[2]})$  if and only if the diagonals of $b$ and $c$ are even.
\end{proof}
    
The Cayley transform will induce a bijection between the  two subsets 
\begin{eqnarray*}
\label{eqnspgq3}
\g s\g p^{\bullet,2}_{g,\m Q}&:=& 
\{
X\in \g s\g p_{g,\m Q}^2 \mid \,
{\rm det}(\mathds{1}-X)\in \m Z^*_{(2)}   \},\\
{\rm Sp}^{\bullet,2}_{g,\m Q}&:=&\{ h\in {\rm Sp}_{g,\m Q}^2 \mid \,
{\rm det}(\tfrac12(\mathds{1}+h))\in \m Z^*_{(2)} \}.
\end{eqnarray*}  
For a matrix $M\in\mc M(2g,\m Z_{(2)})$, i.e. a matrix with odd denominator, 
the condition ${\rm det}(M)\in \m Z^*_{(2)} $,  means that the determinant is invertible in the ring $\m Z_{(2)}$ and hence that the inverse $M^{-1}$
exists and belongs to $\mc M(2g,\m Z_{(2)})$, i.e. it also has odd denominator.

\bl
\label{lemcaythe} 
$(a)$ The Cayley transform  $X\mapsto C(X)=h=(\mathds{1}+X)(\mathds{1}-X)^{-1}$ 
induces a bijection from 
$\g s\g p^{\bullet,2}_{g,\m Q}$ onto ${\rm Sp}^{\bullet,2}_{g,\m Q}$.\\
$(b)$ For $X\in \g s\g p^{\bullet,2}_{g,\m Q}$ and $h=C(X)$, one has the equivalence:
\begin{eqnarray*}
\label{eqnsptspt}
X\in \g s\g p^{\theta,2}_{g,\m Q} 
&\Longleftrightarrow&
h\in {\rm Sp}^{\theta,2}_{g,\m Q}\, .
\end{eqnarray*}
\el

\begin{proof} 
$(a)$ Recall that the equality $h=C(X)$ means 
$
\frac12(\mathds{1}+h)(\mathds{1}-X)=\mathds{1}.
$
For such 
a pair $(X,h)$ in $\mc M(2g,\m Q)$, one has the equivalence\\
${\rm det}(\mathds{1}-X)\in \m Z^*_{(2)}$ $\Longleftrightarrow$ 
${\rm det}(\frac12(\mathds{1}+h))\in \m Z^*_{(2)}$ and, in that case,\\ 
$X\in \mc M(2g,\m Z_{(2)})$ $\Longleftrightarrow$  $\frac12(\mathds{1}+h)\in \mc M(2g,\m Z_{(2)})$.\\
Hence one has the equivalence:\;
$X\in \g s\g p^{\bullet,2}_{g,\m Q}$ $\Longleftrightarrow$  $h\in {\rm Sp}^{\bullet,2}_{g,\m Q}$.

$(b)$ In that case, one has the equivalences:\\ 
$X\in \g s\g p^{\theta,2}_{g,\m Q}$  $\Longleftrightarrow$ 
$\om(Xx,x)\equiv 0\;{\rm mod}\; 2$ for all $x$ in $\m Z_{(2)}^{2g}$ $\Longleftrightarrow$\\ 
$\om((\mathds{1}+X)x,(\mathds{1}-X)x)\equiv 0\;{\rm mod}\; 4$ for all $x$ in $\m Z_{(2)}^{2g}$
 $\Longleftrightarrow$\\ 
$\om(hx',x')\equiv 0\;{\rm mod}\; 4$ for all $x'$ in $\m Z_{(2)}^{2g}$
 $\Longleftrightarrow$
 $h\in {\rm Sp}^{\theta,2}_{g,\m Q}$. 
\end{proof}

%40
\section{Abelian varieties}
\label{secabevar}

The aim of this chapter is to interpret our general construction of critical values
from the point of view of abelian varieties (Theorem \ref{thmmaiabe}).

%41
\subsection{Principally polarized abelian varieties}
\label{secpolvar}

In this section and the next one we fix our choice of notation
and definition.\vs

Let $(A=V/\La,\om)$ be a {\it polarized abelian variety}.
This means that $A$ is a complex torus, with $V=\m C^g$, that $\La$ 
is a lattice in $V$, and that $\om:V\times V\ra \m R$
is a {\it real symplectic form } on $V$ satisfying the following two conditions:\\
$(i)$ the symplectic form $\om$ takes integral values on $\La\times\La$, and\\
$(ii)$ $\om$ is the imaginary part ${\rm Im}(H)$ of a positive hermitian form $H$ on $V$.

Assumption $(ii)$ implies that $\om(iv_1,iv_2)=\om(v_1,v_2)$ for all $v_1$, $v_2$ in $V$ and 
the hermitian form $H$ can be recovered as 
$
H(v_1,v_2)=\om(iv_1,v_2)+i\,\om(v_1,v_2)$.
Our convention is that $H(v_1,v_2)$ is linear in $v_1$ and antilinear in $v_2$.

We will always assume that the  
{\it polarization is  principal},
i.e. that the restriction of $\om$ to $\La\times\La$ has determinant $1$.
At first glance, this assumption looks harmless for us since every
polarized abelian variety is isogenous to a principally polarized abelian variety. 
The problem is that changing $\La$ 
might change the group $G_\nu$ 
in Theorem \ref{thmmaimul} and it 
will be a delicate issue to choose $\La$ 
so that $G_\nu$ is cyclic.

When the polarization is principal, there exists a symplectic basis\\
$(f_1,\ldots f_g,e_1,\ldots, e_g)$ of the lattice $\La$,
i.e. a $\m Z$-basis such that
$$
\om(e_j,e_k)
\;=\;
\om(f_j,f_k)
\;=\;
\om(f_j,e_k)-\de_{jk}
\;=\; 0
\;\;\; \;
\mbox{\rm for all $j,k$}.
$$
The family $(e_1,\ldots,e_g)$ is then a basis of $\m C^g$.
We denote by $\tau$ the $g\times g$ matrix given by
$
(f_1\ldots,f_g)=(e_1,\ldots,e_g)\tau .
$
This matrix $\tau$ 
is a symmetric complex matrix with positive definite 
imaginary part, that is 
$\tau$ belongs to $\mc H_g$. Moreover this map 
$(A,\om)
\longrightarrow \tau$ gives a bijection
\begin{eqnarray}
\label{eqnriethm}
\left\{\begin{array}{c}\!\!\mbox{principally polarized}\!\!\\
\mbox{abelian varieties}
\end{array}\right\}
&\longleftrightarrow&
{\rm Sp}(g,\m Z)\backslash \mc H_g.
\end{eqnarray}

More precisely, for $\tau$ in $\mc H_g$, we introduce the lattice 
$\La_\tau:=\tau\m Z^g\oplus\m Z^g$
of $\m C^g$,
the quotient torus $A_\tau:=\m C^g/\La_\tau$,  
the hermitian form  $H_\tau$ on $\m C^g$ whose matrix is 
$({\rm Im}\tau)^{-1}$ in the canonical basis $(e_1,\ldots,e_g)$ 
and  the imaginary part $\om_\tau$ of $H_\tau$.
The pair $(A_\tau,\om_\tau)$ is then a principally polarized abelian variety, and the map $\tau\mapsto (A_\tau,\om_\tau)$
is the inverse map of \eqref{eqnriethm}.

%42
\subsection{Unitary $\m Q$-endomorphisms}
\label{secisogen}

Let $(A=V/\La,\om)$ be a principally polarized abelian variety.

We denote by ${\rm End}(A)$ the ring of endomorphisms 
$\mu:A\ra A$, and by 
${\rm End}_{\m Q}(A):={\rm End}(A)\otimes_{\m Z}\m Q$
the $\m Q$-algebra of $\m Q$-endomorphisms $\nu$ of $A$.
An {\it isogeny} is an endomorphism of $A$ which is invertible in ${\rm End}_{\m Q}(A)$,
i.e. an endomorphism $\mu$ whose kernel $K_\mu\subset A$ is a finite subgroup.
To each  $\m Q$-endo\-morphism $\nu\in {\rm End}_{\m Q}(A)$ is associated\\
$\star$ a {\it tangent map} $T_\nu\in {\rm End}_{\m C}(V)\simeq \mc M(g,\m C)$,\\
$\star$ a {\it holonomy map} $h_\mu\in {\rm End}_{\m Q}(\La_{\m Q})\simeq \mc M(2g,\m Q)$,
$\La_\m Q:=\La\otimes_\m Z\m Q$.\\
%By extension, to each $\m Q$-endomorphism $\nu$ of $A$ is also associated 
%a tangent map $T_\nu\in {\rm End}_{\m C}(V)$ and
% a holonomy map $h_\nu\in {\rm End}_{\m Q}(\La_\m Q)\simeq \mc M(2g,\m Q)$  where .
The map $h_\nu$ is the restriction of $T_\nu$ to $\La_\m Q$. 
More precisely, an endomorphism (resp.\!  $\m Q$-endomorphism) $\nu$ of $A$ is nothing but 
a $\m C$-endomorphism of $V$ that preserves $\La$ (resp.\! $\La_\m Q$). 
This is why one sometimes writes abusively $\nu$ instead of $T_\nu$ or $h_\nu$.
But it is useful to keep the two notations
because, in coordinates,  $T_\nu$ is a $g\times g$ complex matrix while $h_\nu$ is a 
$2g\!\times\! 2g$ rational matrix.

The {\it Rosati anti-involution} $\nu\mapsto \nu^*$ is the antiinvolution of the $\m Q$-algebra
${\rm End}_{\m Q}(A)$ defined by one of the two equivalent properties:\\
$\star$ $T_{\nu^{*}}$ is the adjoint of $T_\nu$ for the hermitian form $H$ on $V$.\\
$\star$ $h_{\nu^{*}}$ is the adjoint of $h_\nu$ for the symplectic form $\om$ on $\La_\m Q$.
 
\bd 
\label{defisosym} A {\it similarity} of ratio $k$ is an isogeny $\mu$ such that $\mu\mu^*=k^2\mathds{1}$.\\
An isogeny is {\it primitive} if it is not an integral multiple of an isogeny.\\
A {\it unitary} $\m Q$-endomorphism
is a $\m Q$-endomorphism $\nu$ such that $\nu\nu^*\!=\!\mathds{1}$.
\ed

The following lemma is nothing but a useful remark.

\bl
\label{lemisosym}
Let $(A,\om)$ be a principally polarized abelian variety, and  $\mu$ be an isogeny of $A$.
When $k\geq 1$ is an integer, the following are equivalent:\\
$\star$ the isogeny $\mu$ is a similarity of ratio $k$.\\
$\star$ the $\m Q$-endomorphism $\nu:=\frac{1}{k}\mu$ of $A$ is unitary.\\ 
$\star$ the $\m Q$-linear map $h_\nu=\frac{1}{k}h_\mu$ belongs to ${\rm Sp}(\La_\m Q,\om)$.
\el

\bd 
\label{defthestr2} We say that the {\it unitary} $\m Q$-endomorphism
$\nu$  of $A$ {\it preserves a theta structure of level $2$}
if the holonomy $h_\nu$ belongs to the rational symplectic theta subgroup ${\rm Sp}_{g,\m Q}^{\theta,2}$ of level $2$ 
in a symplectic basis of $(\La,\om)$.
%As we have seen in Section \ref{secconcri}, 
This condition does not depend
on the choice of the symplectic basis of $\La$.
\ed

%43
\subsection{Critical values and abelian varieties}
\label{seccriabe}

We can now explain our construction of critical values 
from the point of view of abelian varieties.
In Chapter \ref{seccommul}, 
we will specialize this theorem to the case of 
abelian varieties with complex multiplication.

\bt
\label{thmmaiabe}
Let $(A=V/\La,\om)$ be a principally polarized abelian variety, 
$\nu$ be  a unitary $\m Q$-endomorphism of $A$
preserving a theta structure of level $2$,\;
$T_\nu$ be its tangent map, 
$
G_\nu:=\La/(\La\cap \nu\La)
$
and $d_\nu:=|G_\nu|$.
Then there exists a critical value
$\la_\nu=\ka_\nu\, d_\nu^{1/2}\,{\rm det}_{\m C}(T_\nu)^{1/2}$
on  the group
$G_\nu$ with $\ka_\nu^4=1$.
\et

\br The group $G_\nu$ is naturally isomorphic to the torsion subgroup 
$G_\nu\simeq (\nu^{-1}\La+\La)/\La\subset A$ 
and  the order $d_\nu$ of $G_\nu$ is a divisor of $k^g$ where $k>0$ is chosen 
so that $\mu:=k\nu$ is a $\m Z$-endomorphism of $A$.
In particular, since by assumption, one can choose $k$ odd,
the order $d_\nu$ is odd.

Note also that, since $\nu$ is unitary, the order $d_\nu$ of the group $G_\nu$ and the  absolute value 
of the critical value are related by the equality
$$
|\la_\nu|=d_\nu^{1/2}.
$$ 

Note finally that, once $\nu$ is fixed, except for a fourth root of unity, the ratio $\la_\nu/|\la_\nu|$
does not depend on $A$ in its isogeny class. 
Indeed, changing $A$ in its isogeny class may change $d_\nu$ and the group $G_\nu$ but it does not change the tangent map $T_\nu$.
We will see examples already in Section \ref{secimaqua}. 
\er

The easiest way to determine the square $\ka_\nu^2=\pm 1$ 
is to remember that by Proposition \ref{procrival} 
one has $\la_\nu\equiv 1$ mod $2$. We will see in Section \ref{secimaqua}
that both signs $\pm$ can occur.

\begin{proof}[Proof of Theorem \ref{thmmaiabe}]
The key point is the interrelation between the tangent map $T_\nu$ 
and the holonomy $h_\nu$, together with the use of Proposition \ref{prohsidsi}.
We fix a symplectic $\m Z$-basis 
$(f_1,\ldots,f_g,e_1,\ldots,e_g)$ of $\La$ so that $\om=\sum f_j^*\wedge e_j^*$. 
\begin{align}
\label{eqnmnuhnu}
&\mbox{Let}\;\; m_\nu\in {\rm Sp}_{g,\m Q}^{\theta,2}\;\; \mbox{ be the matrix of $h_\nu^{-1}$ in}\\
&\mbox{the symplectic basis}\;
(e_1,\ldots,e_g,-f_1,\ldots,-f_g)
\nonumber
\end{align}
so that, by \eqref{eqnspgr}, one has $Jm_\nu^{-1}J^{-1}={}^tm_\nu$ and hence the equality in $V^{2g}$
\begin{equation}
\label{eqntftefe}
(T_\nu f_1,\ldots,T_\nu f_g,T_\nu e_1,\ldots,T_\nu e_g)=(f_1,\ldots,f_g,e_1,\ldots,e_g)\,{}^t m_\nu.
\end{equation}
By the adapted symplectic basis in Proposition \ref{prohsidsi},
there exist $\si_1$, $\si_2$ in ${\rm Sp}(g,\m Z)$ and 
an integral matrix ${\bf d}$ with ${\rm det}({\bf d})\neq 0$ such that
\begin{equation}
\label{eqnhsidsi}
m_\nu=\si_1\;D\si_2
\;\;{\rm with}\;\; 
D:=\mbox{\small\scriptsize
$\left(\!
\begin{array}{cc} 
\!{}^t{\bf d}^{-1}\!&{\bf 0}\\  
{\bf 0}&{\bf d}
\end{array}
\!\right)$}.
\end{equation}
The matrix  ${\bf d}$ can be chosen to be a diagonal matrix ${\rm diag}(d_1,\ldots,d_g)$ 
with positive integer coefficients $d_1|d_2|\ldots |d_g$, but we will not use that fact.

We introduce two new symplectic $\m Z$ basis of $\La$.
\begin{align}
\label{eqnffeeffee}
(F_1,\ldots,F_g,E_1,\ldots,E_g)&:=(f_1,\ldots,f_g,e_1,\ldots,e_g)\,{}^t\si_1^{-1},\\
\nonumber
(F'_1,\ldots,F'_g,E'_1,\ldots,E'_g)&:=(f_1,\ldots,f_g,e_1,\ldots,e_g)\,{}^t\si_2.
\end{align}
Combining \eqref{eqntftefe}, \eqref{eqnhsidsi} and \eqref{eqnffeeffee}, one gets
the equalities in $V^g$,
\begin{eqnarray}
\label{eqnbascha0}
(T_\nu F_1,\ldots,T_\nu F_g)&=&(F_1',\ldots,F'_g)\, {\bf d}^{-1}\, ,\\
\nonumber
(T_\nu E_1,\ldots,T_\nu E_g)&=&(E_1',\ldots,E'_g)\, {}^t{\bf d}\, .
\end{eqnarray}
Note that  the group 
$G_\nu\simeq \La/(\La\cap T_\nu\La) \simeq \m Z^g/ {\bf d}\m Z^g$
has order $d_\nu=|{\rm det}_{\m C}({\bf d})|.$

Now we go on our analysis of the unitary $\m Q$-endomorphism $\nu$. We denote by  $\tau$, $\rho$ and $\rho'$ the $g\!\times\! g$ complex matrices 
that give the basis changes  defined by the equalities in $V^g$,
\begin{eqnarray}
\label{eqnbascha}
(F_1,\ldots,F_g)&=&(E_1,\ldots,E_g)\, \tau,\\
\nonumber
(E'_1,\ldots,E'_g)&=&(E_1,\ldots,E_g)\, \rho,\\
\nonumber
(F'_1,\ldots,F'_g)&=&(E_1,\ldots,E_g)\, \rho'.
\end{eqnarray}
Remember that  the matrix $\tau$ is a Riemann matrix, that is  $\tau\in \mc H_g$.

Let $M_\nu$ be the  $g\!\times\! g$ complex matrix that expresses 
$T_\nu$ in the 
basis $(E_1,\ldots,E_g)$ of $\m C^g$ so that one has  the equalities
\begin{eqnarray*}
\label{eqnbascha2}
(T_\nu F_1,\ldots,T_\nu F_g)&=&(F_1,\ldots,F_g)\, \tau^{-1}M_\nu\tau,\\
(T_\nu E_1,\ldots,T_\nu E_g)&=&(E_1,\ldots,E_g)\, M_\nu,
\end{eqnarray*}
and Equalities \eqref{eqnbascha0} can be rewritten as equalities in $\mc M(g,\m C)$:
\begin{eqnarray*}
\label{eqnmbetaud}
M_\nu^{-1}\rho'&=&\tau{\bf d},\\
M_\nu^{-1}\rho&=&{}^t{\bf d}^{-1}.
\end{eqnarray*}

Let
 $
\si:=\si_2\si_1=\mbox{\scriptsize 
$\left(\!\begin{array}{cc} \al&\be\\  
\ga&\de\end{array}\!\right)$} \in {\rm Sp}(g,\m Z) 
$ 
so that by  \eqref{eqnffeeffee}, 
one has 
\begin{align*}
(F'_1,\ldots,F'_g,E'_1,\ldots,E'_g)
&=
(F_1,\ldots,F_g,E_1,\ldots,E_g)\, {}^t\si,
\end{align*}
or, equivalently,
\begin{eqnarray*}
\label{eqnepfpef}
\rho'&=&\tau\, {}^t\al+\,{ }^t\be,\\
\rho&=&\tau\, {}^t\ga+\,{ }^t\de.
\end{eqnarray*}
From these four equalities, one gets
\begin{eqnarray}
\label{eqntaubfd1}
\tau{\bf d}&=&M_\nu^{-1}(\tau\, {}^t\al+{}^t\be),\\
\label{eqntaubfd2}
{}^t{\bf d}^{-1}&=&M_\nu^{-1}(\tau\, {}^t\ga+{}^t\de).
\end{eqnarray}
Hence one has, taking into account that $\tau$ is a symmetric matrix, 
\begin{eqnarray*}
\label{eqntaubfd}
{}^t{\bf d}\tau{\bf d}&=&(\tau\, {}^t\ga+{ }^t\de)^{-1}\,(\tau\, {}^t\al+{ }^t\be),\\
&=&
(\al\tau+\be)\,(\ga\tau+\de)^{-1} .
\end{eqnarray*}
This can be rewritten as
\begin{equation}
\label{eqnsitdtd}
\si\, \tau  \; =\: {}^t{\bf d}\tau{\bf d}.
\end{equation}

Since $h_\nu$ preserves a theta structure of level $2$, by Lemma \ref{lemthestr}, the symplectic matrix $\si$ belongs to ${\rm Sp}_{g,\m Z}^{\theta,2}$.
Therefore 
by Theorem \ref{thmmaithe},
the value
\begin{equation}
\label{eqnlamude}
\la_\nu=j(\si',2\tau) |{\rm det}({\bf d})|= \ka\, {\rm det}_{\m C}(\ga\tau+\de)^{1/2} |{\rm det}({\bf d})|
\end{equation}
is  critical   on the 
group $G_{\bf d}$, 
where $\si':=\!\mbox{\scriptsize 
$\left(\!\begin{array}{cc} \al&\!2\be\!\\  
\!\ga/2\!&\de\end{array}\!\right)$}$ 
and where $\ka$ is  a $8^{\rm th}$ root of unity.
Using \eqref{eqntaubfd2}
one computes,
\begin{align}
\label{eqngataudedm}
{\rm det}_{\m C}(\ga\tau+\de)
& = {\rm det}_{\m C}({\bf d})^{-1}\,{\rm det}_{\m C}( M_\nu).
\end{align}
Plugging this   
into \eqref{eqnlamude}, 
we obtain  the  equality 
$\la_\nu=\ka'\, d_\nu^{1/2}\,{\rm det}_{\m C}(T_\nu)^{1/2}.
$

It remains to explain why the $8^{\rm th}$ root of unity $\ka'$ is a $4^{\rm th}$ root of unity.
This follows from Proposition \ref{procrival} and the claim below.
\vs 

{\bf Claim }  
{\it There exists an odd integer $k$ such that $k^g\, {\rm det}_{\m C}(T_\nu)^{-1}\equiv 1$ mod $2$.}\vs 

\noindent
This claim is true  because there exists an odd integer $k$ such that the complex  matrix $(k T_\nu-1)/2$ preserves a lattice in $\m C^g$. 
Hence all the eigenvalues of $kT_\nu$ are algebraic integers which are equal to $1$ mod $2$. Their product too.
\end{proof}

\br 
\label{remthetatau} When we will deal with the sign issue in Chapter \ref{secsigcri}, 
we will not only need the statement of Theorem \ref{thmmaiabe}
but also the precise way $m_\nu$, $\si$ and $\tau$ are constructed in its proof.  
Note  that, thanks to our choices 
\eqref{eqntftefe}, \eqref{eqnffeeffee} and \eqref{eqnbascha},
the Riemann matrix  $\theta\in \mc H_g$  defined by
$
(f_1,\ldots,f_g)=(e_1,\ldots e_g)\,\theta 
$ 
satisfies the equalities $m_\nu\theta=\theta$ and $\theta=\si_1\tau$.
\er

\br
 As we have seen, the group $G_\nu$ is not cyclic in general but,
even when $g>1$, it may be cyclic of order $d$.
In the next chapters, we will construct many examples of such endomorphisms $\nu$ 
with $G_\nu$ cyclic 
when $A$ has complex multiplication by a CM field.
See for instance Remark \ref{remgnucyc}.
\er

%50
\section{CM number fields}
\label{seccommul}

The aim of this chapter is to specialize our general construction of critical values to the case of  CM abelian varieties and to express it
(Theorem \ref{thmmaimul}) from the point of view of 
symplectic CM algebras 
$(K,\om)$ and their autodual lattices $\La$.

We will first recall in Sections \ref{seccmfield},
and \ref{seccmabelian} 
a few classical definitions and facts.
See \cite{Shimura98}, \cite[Chapter 1]{Milne06} and 
\cite[Chapter 1]{Streng10} for details.

In order to state this theorem,
we will introduce in Section \ref{seccmfthe} the
theta subgroup $U^{\theta,2}_{K,\La}$ of level $2$ of the unitary subgroup $U_K$ of $K$,
and we will describe it thanks to the Cayley transform
(Lemma \ref{lemimauni}).

%51
\subsection{CM algebras, CM types and reflex norms}
\label{seccmfield}

A number field $K\subset \m C$  is said 
to have complex multiplication or to be a {\it CM number field}
if $K$ is a totally imaginary quadratic extension of a totally real number field $K_0$. 
More generally a {\it CM algebra} $K$ is a product of CM fields $K=\prod_jK_j$. It has  even dimension 
$2g:=dim_{\m Q}K$. We denote by $x\mapsto \ol{x} $ 
the {\it complex conjugation} of $K$, that is the automorphism of $K$ equal to the complex conjugation on each $K_j$.

An algebraic number $x$, or an element $x\in K$,  is said to be {\it totally real} (resp.  {\it totally imaginary}, resp. {\it totally unitary}) if all its Galois conjugates in $\m C$ are real (resp. imaginary, resp. unitary). For instance 
$x=\sqrt{2}$ (resp. $x=i\sqrt{2}$, resp. $x=\frac{3+4i}{5}$).
We denote by $K_0$ (resp. $\mc I_K$, resp. $U_K$) the subset of totally real 
(resp. totally imaginary, resp. totally unitary) elements of $K$.

Let $m\geq 1$ be an integer.
An algebraic integer $\mu$  is said to be a 
{\it $m$-Weil number} if all its complex Galois conjugates have modulus equal to $\sqrt{m}$.
\vs 

Those notions are strongly related:\\
$\star$ In a CM field $K\subset \m C$, every real (resp. imaginary or unitary) element is 
totally real (resp. totally imaginary or totally unitary).\\
$\star$ When $t$ is a nonzero totally imaginary algebraic number,
the field $\m Q[t]$ is  CM.
Conversely, when $x$ is in a CM algebra, the difference
$t:=x-\ol x$ is totally imaginary.\\
$\star$ When $\nu\neq \pm 1$ is a totally unitary algebraic number,
the field $\m Q[\nu]$ is  CM.
Conversely, when $x$ is invertible in a CM algebra, the ratio
$\nu:=x/\ol x$ is totally unitary.\\
$\star$ Let $k>0$ be an integer. When $\mu$ is a $k$-Weil number,
the ratio  
$\nu:=\mu/k^{1/2}$ is totally unitary.
Conversely, when $\nu\in \m C$ is totally unitary,
if the  multiple $\mu:=k^{1/2}\nu$ is integral, it is a $k$-Weil number.
\vs 

The following fact implies that the totally unitary 
algebraic numbers $\nu\neq -1$ are exactly the 
images of totally imaginary numbers by the Cayley transform.
For a CM algebra $K$, we set
\begin{eqnarray*}
\mc I^\bullet_K&:=&\{ t\in \mc I_K\mid 1-t\in K^*\},\\
U^\bullet_K&:=&\{\nu\in U_K\mid \nu+1\in K^*\},
\end{eqnarray*}
where $K^*$ is the group of 
invertible elements of $K$.

\bfa
\label{facimauni} 
Let $K$ be a CM algebra.
The maps $t\mapsto \nu= \frac{1+t}{1-t}$ and $\nu\mapsto t=\frac{1-\nu}{1+\nu}$ are inverse  bijections between $\mc I^\bullet_K$ and $U^\bullet_K$.
\efa

\begin{proof} 
Indeed the condition $\ol{\nu}=\nu^{-1}$ is equivalent to $\ol{t}=-t$. \end{proof}

Let $K$ be a CM algebra of dimension $2g$.
A $CM$-type $\Ph$ of $K$
is a  tuple $\Ph=(\rho_1,\ldots, \rho_g)$ of distinct embeddings of $K$ into $\m C$
no two of which are complex conjugate. We will call $(K,\Phi)$ a CM pair
and we will think of $\Ph$ as an algebra morphism $\Phi:K\mapsto \m C^g$.

The CM pair $(K,\Ph)$ is said to be {\it primitive}, if 
$K$ is a field and the restriction of $\Ph$ to any proper
CM subfield is not a CM type.

For an element $\mu$ in  $K$,
we denote by 
$$
N_\Ph(\mu):=\rho_1(\mu)\cdots\rho_g(\mu)
$$
its {\it reflex norm} or {\it type norm}. 
It lives in the reflex field $K^r\subset \m C$. More precisely, 
the {\it reflex field} $K^r$ is the subfield of $\m C$ spanned by the reflex norms
of the elements of $K$. 
The field $K^r$ is always a CM field.

%52
\subsection{CM abelian varieties}
\label{seccmabelian}

An abelian variety is {\it simple} 
if it does not contain proper abelian subvarieties. 
Every abelian variety is isogenous 
to a product of simple abelian varieties.

A {\it CM abelian variety} is an abelian variety $A$
such that ${\rm End}_{\m Q}(A)$ contains a
CM algebra $K_A$ of $\m Q$-dimension $2\dim(A)$.
An abelian variety that is 
isogenous to a CM abelian variety is also CM.
The action of $K_A$ on the tangent space $TA$
of $A$ gives a CM type $\Phi_A:K_A\ra\m C^g$.

A CM abelian variety is {\it simple} if and only if the CM pair $(K_A,\Ph_A)$ is primitive, and 
the map 
\begin{equation} 
\label{eqnakapha}
A\mapsto (K_A,\Ph_A) 
\end{equation}
is a bijection 
between the set of isogeny classes of simple CM abelian varieties 
and the set of isomorphism classes of 
primitive CM pairs. 
\vs

We now recall the inverse map to \eqref{eqnakapha}, 
i.e. the construction of the abelian variety $A$
starting from a CM pair $(K,\Phi)$. 
The most interesting case is when $K$ is a CM number field
and $\Phi$ is primitive.
But, in some examples (as Propositions \ref{proparacri}
and \ref{prorarbrc1})
it will be useful to deal 
with a more general CM pair $(K,\Phi)$.

Let $K=\prod K_j$ be a CM algebra.
Let $\mc O_K=\prod \mc O_{K_j}$, where $\mc O_{K_j}$ is  
the {\it ring of integers of $K_j$}. 
A lattice $\La$ of $K$ is an additive  subgroup of 
$K$
that is commensurable to $\mc O_K$.
It can be seen, via any CM type  $\Ph:K\ra \m C^g$, as a lattice in $\m C^g$.
When there is no possible confusion on the choice of $\Phi$, we will also write  $\La$ for $\Phi(\La)$. 
The following construction due to Shimura is very useful, since it tells us that the complex torus 
$A=\m C^g/\Phi(\La)$ is an abelian variety.
To state it we need more notation.

Given  $t_0\in K$ totally imaginary  and  invertible,
we introduce the $\m Q$-symplectic form $\om_0=\om_{t_0}$ of 
the $2g$-dimensional $\m Q$-vector space $K$
\begin{equation} 
\label{eqnomotrk}
\om_{t_0}(x,x')= {\rm Tr}_{K/\m Q}
(\tfrac{x\ol{x'}}{t_0})
\;\; \mbox{\rm for all $x,x'$ in $K$}.
\end{equation}
Such a pair $(K,\om_{t_0})$ will be called a {\it symplectic CM algebra}.
Replacing $\La$ by a multiple, we can assume
that $\La\subset \mc O_K$ 
and $\om_{t_0}(\La,\La)\subset \m Z$. 
We also denote by $\om_{t_0}$ the $\m R$-linear extension 
of $\om_{t_0}$ to the $2g$-dimensional $\m R$-vector space $\m C^g$.
We introduce then  the CM type $\Ph_0=\Ph_{\om_0}=\Ph_{t_0}$ of $K$ given by
\begin{equation}
\label{eqnimrhot}
\Ph_{t_0}=\{\rho\in {\rm Hom}_{\rm alg}(K,\m C)\mid {\rm Im}( \rho(t_0))>0\}.
\end{equation}
Note that, for any CM type $\Ph$ of $K$, one can find an invertible imaginary element $t_0$ with $\Phi_{t_0}=\Phi$.
By Definition \eqref{eqnimrhot}, this symplectic form
$\om_{t_0}$ is a polarization on $A=\m C^g/\Ph_{t_0}(\La)$, and hence $(A,\om_{t_0})$ 
is a polarized abelian variety.
For our construction  we will need $\om_{t_0}$ to have determinant $1$ on $\La$
so that $\La$ is autodual with respect to $\om_{t_0}$,
which means that $(A,\om_{t_0})$ is principally polarized.

%53
\subsection{Multiplicative theta subgroup of level $2$}
\label{seccmfthe}
We define now the unitary theta subgroup of $K^*$ of level $2$.
\vs 
 
We recall that $\m Z_{(2)}$ is the ring of rational numbers with odd denominator.
For $\La$ lattice in $K$ we set
$\La_{(2)}=\m Z_{(2)}\La$. For instance $\mc O_{K,(2)}$ is the set of ratios of elements of $\mc O_K$ with a denominator prime to $2$.

Let $\La$ be an autodual lattice in a symplectic  CM algebra $(K,\om_0)$ as above.
We introduce two subgroups of the unitary  group 
$U_K=\{\nu\in K\mid \nu\ol{\nu}=1\}$.

\bd 
\label{defkthlal}
The congruence subgroup $U_{K,\La}^2$ of level $2$ and the theta subgroup $U^{\theta,2}_{K,\La}$ of level $2$ are
\begin{eqnarray}
\nonumber
U_{K,\La}^2 &:=&\{\nu\in U_K\mid 
(\nu-1)(\La_{(2)})\subset 2\La_{(2)}\, \}\\
\label{eqnuklath2}
U^{\theta,2}_{K,\La}&:=&\{\nu\in U_{K,\La}^2\mid 
\om_0(\nu x,x)\in 4\,\m Z_{(2)}
\;\mbox{for all}\; x\in \La_{(2)}\}
\end{eqnarray}
\ed

One can also define $U_{K,\La}^2$ as the set of $\nu\in U_K$ for which
there exists an odd integer $d$
such that $d\nu(\La)\subset \La$ and  $d\nu$ acts trivially on the 
quotient $\La/2\La$. 
Note that the elements $\nu$ of $U_{K,\La}^2$ preserve $\La_{(2)}$ and act trivially 
on $\La_{(2)}/2\La_{(2)}$. 
\vs 

Let us fix a symplectic basis of  $\La$. 
For $\nu$ in $K$, set  $m_\nu\in {\mc M}(2g,\m Q)$
for the transpose of the matrix of the multiplication by $\nu$ in this basis. 
When the element $\nu$ is in $U_K$ this matrix is symplectic: 
$m_\nu\in {\rm Sp}(g,\m Q)$. 
Here is an equivalent definition 
\begin{eqnarray*}
U_{K,\La}^2 &=&\{\nu\in U_K\mid 
m_\nu\in {\rm Sp}_{g,\m Q}^2\, \}\\
U^{\theta,2}_{K,\La}&=&\{\nu\in U_K\mid 
m_\nu\in {\rm Sp}^{\theta,2}_{g,\m Q}\}
\end{eqnarray*}
This definition does not depend on the choice of the symplectic basis
of $\La$.

We also define the analog subsets in the imaginary elements of $K$:
\begin{eqnarray}
\nonumber
\mc I_{K,\La}^2 &:=&\{t\in \mc I_K\mid 
t\La_{(2)}\subset \La_{(2)}\, \}\\
\label{eqniklath2}
\mc I^{\theta,2}_{K,\La}&:=&\{t\in \mc I_{K,\La}(2)\mid 
\om_0(t x,x)\in 2\,\m Z_{(2)}
\;\mbox{for all}\; x\in \La_{(2)}\}
\end{eqnarray}
or, equivalently,
\begin{eqnarray*}
\mc I_{K,\La}^2&=&\{t\in \mc I_K\mid 
m_t\in {\g s\g p}_{g,\m Q}^2\, \}\\
\mc I^{\theta,2}_{K,\La}&=&\{t\in \mc I_K\mid 
m_t\in {\g s\g p}^{\theta,2}_{g,\m Q}\}
\end{eqnarray*}

The following lemma is an analog of Lemma \ref{lemcaythe} 
using the Cayley transform of Fact \ref{facimauni}. 
It will give a bijection between the following two subsets.
\begin{eqnarray*}
\mc I^{\bullet,2}_{K,\La}&:=&\{ t\in \mc I_{K,\La}^2\mid N_{K/\m Q}(1\!-\!t)\in \m Z^*_{(2)}\},\\
U^{\bullet,2}_{K,\La}&:=&\{\nu\in U_{K,\La}^2\mid N_{K/\m Q}((\nu\!+\!1)/2)\in \m Z^*_{(2)}\}.
\end{eqnarray*}
For an element  $x\in K$ preserving $\La_{(2)}$, 
the condition $N_{K/\m Q}(x)\in \m Z^*_{(2)}$,  tells us that 
the inverse $x^{-1}$
exists and also preserves $\La_{(2)}$.

When $K$ is a CM field, this condition  $N_{K/\m Q}(x)\in \m Z^*_{(2)}$
means that $x$ is  a ratio of two elements of $\mc O_K$ which are prime to $2$.

\bl
\label{lemimauni} Let $(K,\om_0)$ be a symplectic CM algebra, 
$\La\subset K$ an autodual lattice.
$(a)$ The Cayley transform $t\mapsto \frac{1+t}{1-t}$
is a  bijection from $\mc I^{\bullet,2}_{K,\La}$ onto $U^{\bullet,2}_{K,\La}$.\\
$(b)$ For $t\in \mc I^{\bullet,2}_{K,\La}$ and $\nu=\frac{1+t}{1-t}$, one has the equivalence
%\begin{eqnarray*}
\label{eqniltult}
$t\in \mc I^{\theta,2}_{K,\La}
\Longleftrightarrow
\nu\in U^{\theta,2}_{K,\La} \, .$
%\end{eqnarray*}
\el

\begin{proof}  The proof is the same as for Lemma \ref{lemcaythe}.

$(a)$ The equality $\nu=\frac{1+t}{1-t}$ can be rewritten as 
$
(1-t)(1+\nu)/2=1.
$
For such 
a pair $(t,\nu)$, one has the equivalence\\
$N_{K/\m Q}(1-t)\in \m Z^*_{(2)}$ $\Longleftrightarrow$ 
$N_{K/\m Q}((1+\nu)/2)\in \m Z^*_{(2)}$ and, in that case,\\ 
$t\La_{(2)}\subset \La_{(2)}$ $\Longleftrightarrow$  $(1+\nu)\La_{(2)}\subset 2\La_{(2)}$.\\
Hence one has the equivalence:\;
$t\in \mc I^{\bullet,2}_{K,\La}$ $\Longleftrightarrow$  $\nu\in U^{\bullet,2}_{K,\La}$.

$(b)$ In that case, one has the equivalences:\\ 
$t\in \mc I^{\theta,2}_{K,\La}$  $\Longleftrightarrow$ 
$\om_0(tx,x)\equiv 0\;{\rm mod}\; 2$ for all $x$ in $\La_{(2)}$ $\Longleftrightarrow$\\ 
$\om_0((1+t)x,(1-t)x)\equiv 0\;{\rm mod}\; 4$ for all $x$ in $\La_{(2)}$
 $\Longleftrightarrow$\\ 
$\om_0(\nu x',x')\equiv 0\;{\rm mod}\; 4$ for all $x'$ in $\La_{(2)}$
 $\Longleftrightarrow$
 $\nu\in U^{\theta,2}_{K,\La}$. 
\end{proof}

%54
\subsection{Critical values and CM field}
\label{seccricom}

All the examples of critical values found in the lists of \cite[Section 1.5]{CSAGI} can be explained thanks to the following theorem. 
This theorem  is a special case of our main theorem \ref{thmmaiabe} 
for CM abelian varieties and is expressed
from the point of view of CM number fields.

\bt 
\label{thmmaimul} 
Let $K$ be a CM algebra of degree $2g$, 
$t_0\in \mc I_K$ invertible,
$\om_{t_0}$ be the symplectic form as in \eqref{eqnomotrk} and $\Phi=\Phi_{t_0}$ be the CM type  as in \eqref{eqnimrhot}.

Let $\La\subset K$ be an autodual lattice and $\nu\in K$ be a unitary element that belongs to  the unitary theta subgroup  $U^{\theta,2}_{K,\La}\subset U_K$ of level $2$.

Let $G_\nu:=\La/(\La\cap \nu \La)$  and $d_\nu:=|G_\nu|$.
Then there exists a critical value
$\la_\nu=\ka_\nu\, d_\nu^{1/2}\,N_\Ph(\nu)^{1/2}$  on $G_\nu$ with $\ka_\nu^4=1$.
\et

Note that this critical value $\la_\nu$ 
on $G_\nu$ is a $d_\nu$-Weil number since 
$N_\Phi(\nu)$ is a totally unitary algebraic number.

This gives a critical value $\la_\nu$  modulo a fourth root of unity $\ka_\nu$.
We can  determine this root of unity by using Corollary \ref{cormaisym} but
note that replacing $\nu$ by $-\nu$ changes the sign of 
$N_\Ph(\nu)$ when $g$ is odd.
It is easier to determine $\ka_\nu$ up to sign by using the fact that 
$(\la_\nu-1)/2$ is an algebraic integer.

\begin{proof}[Proof of Theorem \ref{thmmaimul}]
As seen in Section \ref{seccmabelian}, 
the pair $(A=\m C^g/\Phi(\La),\om_{t_0})$ is a principally polarized abelian variety.
The element $\nu\in K$ induces a unitary $\m Q$-endomorphism of $A$  that preserves a theta structure of level $2$ on $A$.
Theorem \ref{thmmaiabe} gives us a critical value  
$\la_\nu:=\ka_\nu\, d_\nu^{1/2}\,\det_\m C(T_\nu)^{1/2}$
on $G_\nu$. We conclude thanks to the equality
$N_\Phi(\nu)=\det_\m C(T_\nu)$.
\end{proof}

%60
\section{Examples}
\label{seclisexa}

In this Chapter we show on three examples how to 
compute explicit critical values by using 
Theorem \ref{thmmaimul}. 
In each of these examples, we 
follow always the same strategy, keeping the notation \eqref{eqnmnuhnu} as in the proof of Theorem \ref{thmmaiabe}.\vs 

\noindent
$(a)$ We choose the CM pair $(K,\Phi)$,  and its symplectic form $\om_0$.\\
$(b)$ We choose a first lattice $\La=\La_0$ in $K$ that is autodual for $\om_0$.\\
$(c)$ We choose $\nu$ unitary in $K$ and  compute the group $G_{0,\nu}:=\La_0/(\La_0\cap\nu\La_0)$.\\
$(d)$ We modify slightly $\La$ and $\om_0$ so that $\nu$ belongs to $U_{K,\La}^{\theta,2}$.\\
$(e)$ We compute the square root of the reflex norm $N_\Phi(\nu)$.\\
$(f)$ We comment on the $4^{\rm th}$-root of unity involved in the formula for $\la_\nu$.
In this chapter, we will determine this $4^{\rm th}$-root of unity 
only up to sign. 
The determination of this sign will be given in 
Chapter \ref{secsigcri}.
\vs 

In these examples, the CM algebra $K$ is not always a number field
and the lattice $\La$ is not always
a fractional ideal of $K$.

%61
\subsection{Imaginary quadratic fields}
\label{secimaqua}

We  begin by the case where $K$ is an imaginary quadratic field
and hence $A$ is a CM elliptic curve. 
This case, which  has already been worked out in \cite{CSAGI},
will help the reader to understand the above strategy.

\bp
\label{prorapirb0}
Let $d\!=\!a\!+\!b$  with $a,b$ positive integers 
 and  
$a\!\equiv\!\frac{(d+1)^2}{4}\;{\rm  mod}\; 4$.
Then the complex number $\la:=\sqrt{a}+i\sqrt{b}$
is a $d$-critical value.
\ep

\begin{proof} 
Note that the congruence relation on $a$ 
is equivalent to 
\begin{equation}
\label{eqnconrel} 
a-b\equiv 1\;\; {\rm mod}\; 4
\;\;\;{\rm and}\;\;\;
ab\equiv 0\;\; {\rm mod}\; 4.
\end{equation}

\noindent $(a)$ We choose for CM algebra the field $K=\m Q[\al]$ with $\al=i\sqrt{ab}$. 
We introduce the $\m Q$-basis $(e_{0,1},f_{0,1})$ of $K$ where $e_{0,1}=1$ and $f_{0,1}=\al$.

We choose the symplectic form $\om_0$ on the $\m Q$-vector space $K$ to be 
$$\textstyle
\om_0(x,x')={\rm Tr}_{K/\m Q}(\frac{x\ol{x}'}{2\al})
$$ 
so that $\om_0=f_{0,1}^*\wedge e_{0,1}^*$, i.e.  $\om_0(f_{0,1},e_{0,1})=1$ . 
Since ${\rm Im}(\al)>0$, the CM type defined by 
\eqref{eqnimrhot} is the singleton $\Phi=\{\rho\}$ where $\rho$ is the injection of $K$ in $\m C$. 

\vs 

\noindent $(b)$  We first introduce the lattice $\La_0=\m Z[\al]=\m Z e_{0,1}\oplus \m Z f_{0,1}$. 
This lattice is autodual for the symplectic form $\om_0$.
\vs 

\noindent $(c)$  
We choose the totally unitary element $\nu$ of $K$ to be 
\begin{equation} 
\label{eqnhuabab}
\nu:=
\frac{\sqrt{a}+i\sqrt{b}}{\sqrt{a}-i\sqrt{b}}
=\frac{1}{d}(a-b+2i\sqrt{ab}).
\end{equation}
The matrix $m_{0,\nu}$ of multiplication by $\nu^{-1}$ in the basis $(e_{0,1},-f_{0,1})$ is
\begin{equation*} 
m_{0,\nu}=\frac{1}{d}
\mbox{
$\left(\!\begin{array}{cc} a\!-\!b&-2ab\\  
2&a\!-\!b\end{array}\!\right)$}. 
\end{equation*}
Since the coefficients of the matrix $dm_{0,\nu}$ are integral with gcd equal to $1$,
the group $G_{0,\nu}=\La_0/(\La_0\cap \nu\La_0)$ is isomorphic to $\m Z/d\m Z$. 
The fixed point of $m_{0,\nu}$ in $\mc H_1$ 
is $\theta_0:= i\sqrt{ab}$.\vs 

\noindent $(d)$
The matrix $m_{0,\nu}$ is symplectic 
but does not  belong to the symplectic theta subgroup  ${\rm Sp}^{\theta,2}_{2,\m Q}$ 
because the  lower-left coefficient
of  $m_{0,\nu}$
is not equal to $0$ mod $4$.
This is why we introduce the basis $e_1$, $f_1$ with 
$e_1=2e_{0,1}$ and $f_1=f_{0,1}$ and the lattice $\La:=\m Z e_1\oplus\m Zf_1$.
This lattice $\La$ is a sublattice of $\La_0$ of index $2$
that is autodual for $\tfrac12 \om_0$. 
Since $d$ is odd,
one still has
$$
G_\nu=\La/(\La\cap\nu\La)\simeq \La_0/(\La_0\cap\nu\La_0)=G_{0,\nu}\simeq \m Z/d\m Z,
$$ 
The matrix $m_{\nu}$ of multiplication by $\nu^{-1}$ in the basis $(e_1,-f_1)$ is
\begin{equation*} 
m_{\nu}=\frac{1}{d}
\mbox{
$\left(\!\begin{array}{cc} a\!-\!b&-ab\\  
4&a\!-\!b\end{array}\!\right)$}. 
\end{equation*}
Since $ab\equiv 0$ mod $4$, the element 
$\nu$  belongs to the unitary theta subgroup  $U^{\theta,2}_{K,\La}$ of level $2$.
The fixed point of $m_{\nu}$ in $\mc H_1$ 
is $\theta:= i\sqrt{ab}/2$.
\vs 

\noindent $(e)$  Hence by Theorem \ref{thmmaimul}, one has a $d$-critical value\\
$\la_\nu= \ka_\nu d^{1/2}N_\Phi(\nu)^{1/2}= \ka_\nu(\sqrt{a}+i\sqrt{b})$
with $\ka_\nu^4=1$.
\vs 

\noindent $(f)$ 
Since by Proposition \ref{procrival}, the ratio $(\la_\nu -1)/2$ is an algebraic integer,
one  has $\ka_\nu=\pm 1$. 
The precise sign $\pm$ in $\ka_\nu$  will be computed in 
Section \ref{secsigimaqua}. 
\end{proof}

\br 
The extension of Proposition  \ref{prorapirb0} when  
$a-\frac{(d+1)^2}{4}\equiv 2$ mod $4$ is not correct.
This explains why,  in the assumptions of Theorem \ref{thmmaimul},  the unitary theta group 
$U^{\theta,2}_{K,\La}$ of level $2$
can not be replaced by the congruence unitary subgroup 
$U_{K,\La}^2$ of level $2$.
\er

\br
Note that in this proof, when $a$ and $b$ are not relatively prime, 
we can not choose $\La$ to be the whole ring of integers $\mc O_K$.

Note also that if we fix $d_0= a_0 +b_0$ with $a_0\equiv \frac{(d_0+1)^2}{4}$ mod $4$
and let $d=md_0$, $a=ma_0$, $b=mb_0$ where $m$ varies among the positive integers
$m\equiv 1$ mod $4$, the value of the totally unitary element $\nu$ in \eqref{eqnhuabab}
does not depend on $m$, but the field $K$, the lattice $\La$, 
the group $G_\nu\simeq\m Z/d\m Z$ and the critical value $\la_\nu$ do depend on $m$.
\er

%62
\subsection{Products of imaginary quadratic  fields}
\label{secproima}

The new critical values that we obtain in this section 
(Proposition \ref{proparacri})
are obtained by using products of imaginary fields.
\vs 

Let $d$ be the product $d=d_1d_2$ of two odd  integers.
When $\la_1$ is a $d_1$-critical value and $\la_2$ is a $d_2$-critical value,
then the product $\la=\la_1\la_2$ is a critical value on the product 
$\m Z/d_1\m Z\times \m Z/d_2\m Z$ (see \cite[Section 1.4]{CSAGI}). 
In particular, when $d_1$ and $d_2$ are coprime, this value  
$\la$ is  $d$-critical.

In the list of \cite[Section 1.5]{CSAGI}, 
we found a $d$-critical value  with $d=15$,
$$
\la= (\sqrt{3}+ i\sqrt{2})(\sqrt{2}+ i).
$$  
This $d$-critical value is of the form 
$\la=\la_1\la_2$ but is not obtained by the above process.
Indeed, $\la_1:=\sqrt{3}+ i\sqrt{2}$ and $\la_2:=\sqrt{2}+ i$ are not critical values because
the ratios $(\la_j-1)/2$ are not  algebraic integers.

The following proposition explains why this value $\la$ is $d$-critical.
The proof uses an abelian surface that is isogenous but not isomorphic 
to a pro\-duct of two elliptic curves.

\bp
\label{proparacri}
Let $d=d_1d_2$ be a product of two coprime odd numbers. 
For $j=1,2$, let $d_j=a_j\!+\!b_j$ all positive,
$a_j\!-\! \frac{(d_j+1)^2}{4}\equiv 2$ mod $4$,
and $\eps_j\!=\!\pm 1$.
Then 
$\la:=(\sqrt{a_1}\! +\! i\eps_1\sqrt{b_1})(\sqrt{a_2}\! +\! i\eps_2\sqrt{b_2})
$ is $d$-critical.
\ep

For instance~:\\ 
$\star$ $\la= (\sqrt{2}+ i)(\sqrt{3}+ i\sqrt{2})$ and $\la= (\sqrt{2}+ i)(\sqrt{3}- i\sqrt{2})$  are $15$-critical.\\
$\star$ $\la= (\sqrt{2}+ i)(\sqrt{6}+ i)$  and $\la= (\sqrt{2}+ i)(\sqrt{6}- i)$  are $21$-critical.\\
Notice that the first two $\la$'s are Galois conjugate but that the last two are not. 
This last example emphasizes that we have to pay attention to the signs $\eps_j$ occuring in the 
formula for $\la$. 

\begin{proof} 
Note that the congruence relation on $a_j$, $b_j$  is equivalent to 
\begin{equation*}
a_j-b_j\equiv 1\;\; {\rm mod}\; 4
\;\;\;{\rm and}\;\;\;
a_jb_j\equiv 2\;\; {\rm mod}\; 4.
\end{equation*}

\noindent $(a)$ We choose for CM algebra the product $K=K_1\times K_2$ 
of the quadratic fields $K_j=\m Q[\al_j]$ with $\al_j=i\sqrt{a_jb_j}$. 
We introduce the $\m Q$-basis $e_{0,1}$, $e_{0,2}$, $f_{0,1}$, $f_{0,2}$ of $K$ where 
$e_{0,1}=(1,0)$, $e_{0,2}=(0,1)$, $f_{0,1}=(\al_1,0)$ and $f_{0,2}=(0,\al_2)$.
Note that the element $\al:=(\al_1,\al_2)\in K$ is totally imaginary and invertible.

We define the symplectic form $\om_0$ on $K$ by 
$$\textstyle
\om_0(x,x')={\rm Tr}_{K/\m Q}(\frac{x\ol{x'}}{4\al})
$$ 
so that $\om_0=f_{0,1}^*\wedge e_{0,1}^*+f_{0,2}^*\wedge e_{0,2}^*$. 
Since both ${\rm Im}(\al_j)>0$, the CM type defined by 
\eqref{eqnimrhot} is  $\Ph=\{\rho_1,\rho_2\}$, where $\rho_j:K\ra \m C$ is the projection on the factor $K_j\subset\m C$. 
\vs 

\noindent $(b)$  We choose for lattice $\La_0=\m Z e_{0,1} \oplus \m Ze_{0,2}\oplus\m Z f_{0,1}\oplus \m Z f_{0,2}$. 
This lattice is autodual for the symplectic form $\om_0$.
\vs 

\noindent $(c)$ We set $d'_j=\pm d_j$ so that $d'_j\equiv 1$ mod $4$. 
We choose the  unitary element 
\begin{equation*} 
\nu:=(\nu_1,\nu_2)\in K\;\;\; \mbox{\rm with}\;\;\;
\nu_j= 
\frac{1}{d'_j}(a_j-b_j+2i\eps_j\sqrt{a_jb_j})\in K_j.
\end{equation*}

The same computation as the one given in Section \ref{secimaqua},
proves that the group $G_{0,\nu}$ is isomorphic to $\m Z/d\m Z$, Indeed,
$$
G_{0,\nu}=\La_0/(\La_0\cap\nu\La_0)\simeq (\m Z/d_1\m Z)\times (\m Z/d_2\m Z)\simeq \m Z/d\m Z
$$
The matrix  of multiplication by $\nu^{-1}$  in the basis $e_{0,1}$, $e_{0,2}$, $-f_{0,1}$, $-f_{0,2}$ is
\begin{equation} 
\label{eqnmonumod}
m_{0,\nu}= \mbox{\scriptsize 
$\left(\!\begin{array}{cccc} u_1&0&\eps_1v_1&0\\  
0&u_2&0&\eps_2v_2\\
\eps_1w_1&0&u_1&0\\
0&\eps_2w_2&0&u_2\end{array}\!\right)$},
\;\; {\rm with}
\end{equation}
\begin{equation}
\label{eqnujvjwj}
u_j\!:=\!\tfrac{a_j-b_j}{d'_j}\equiv 1\;{\rm mod}\;4,\;\;
v_j\!:=\!\tfrac{-2a_jb_j}{d'_j}\equiv 4\;{\rm mod}\;8,\;\; 
w_j\!:=\!\tfrac{2}{d'_j}\equiv 2\;{\rm mod}\;8.
\end{equation}
For $j\!=\!1,2$ set $\theta_j\!:=\!i\sqrt{a_jb_j}$.
The fixed point of $m_{0,\nu}$ in $\mc H_2$ is
$\theta_0:=\mbox{\scriptsize 
$\left(\!\begin{array}{cc} \!\theta_{1}\!&0\\  
0&\!\theta_{2}\!\end{array}\!\right)$}.$
\vs 

\noindent $(d)$  Since $\eps_1$ and $\eps_2$ are odd the element $\nu$ is not in $\mc U^{\theta,2}_{K,\La_0}$.
This is why we have to choose another basis 
\begin{align}
\label{eqnnewbas}
(e_1,e_2,-f_1,-f_2)&=
(e_{0,1},e_{0,2},-f_{0,1},-f_{0,2}) \, P
\end{align}
and another lattice $\La=\m Z e_1\oplus \m Ze_2\oplus \m Z f_1 \oplus \m Z f_2$.
This lattice $\La$ will be a sublattice of index $4$ in $\La_0$
that will be autodual for $\tfrac12 \om_0$, 
and such that $\nu$ belongs to $U^{\theta,2}_{K,\La} $.
Since $d$ is odd,
one will still have
$$
G_\nu=\La/(\La\cap\nu\La)\simeq \La_0/(\La_0\cap\nu\La_0)=G_{0,\nu}\simeq \m Z/d\m Z.
$$ 

To explain precisely  this choice of $\La$ 
we need to distinguish two subcases that we call $(d1)$ and $(d2)$. With no loss of generality, we assume $\eps_1=1$.\vs 

\noindent $(d1)$ We assume that $\eps_1=-\eps_2=1$. 

We choose 
the new basis \eqref{eqnnewbas} to be given by the basis change matrix 
$
P=
\mbox{
$\left(\!\begin{array}{cc} {\bf p}&{\bf 0}\\  
{\bf 0}&{\bf p}
\end{array}\!\right)$}
$ 
with 
$
{\bf p}=
\mbox{\scriptsize 
$\left(\!\begin{array}{cc} 1\!&1\\  
1\!&\!-1\!
\end{array}\!\right)$}
$ 
so that
\begin{eqnarray*}
\label{eqnbaschaid1}
(e_1,e_2,-f_1,-f_2)
&=&(e_{0,1}\!+\!e_{0,2},\; e_{0,1}\!-\!e_{0,2},\;-f_{0,1}\!-\!f_{0,2},\;-f_{0,1}\!+f_{0,2})
\end{eqnarray*}
This basis is symplectic for the form $\tfrac12\om_0$.

The matrix $m_{\nu}$ of multiplication by $\nu^{-1}$ in 
the basis $e_1$, $e_2$, $-f_1$, $-f_2$  
is 
\begin{equation} 
\label{eqnmnumod}
m_{\nu}=P^{-1}m_{0,\nu}P=\mbox{\scriptsize 
$\left(\!\begin{array}{cccc} u_{+}&u_{-}&v_{-}&v_{+}\\  
u_{-}&u_{+}&v_{+}&v_{-}\\
w_{-}&w_{+}&u_{+}&u_{-}\\  
w_{+}&w_{-}&u_{-}&u_{+}\end{array}\!\right)$},
\end{equation} 
with $u_{\pm}=(u_1\pm u_2)/2$, \;
$v_{\pm}=(v_1\pm v_2)/2$ 
and $w_{\pm}=(w_1\pm w_2)/2$.
 
The congruence conditions \eqref{eqnujvjwj} imply that 
$m_\nu\equiv {\bf 1}$ mod $2$ and that the diagonal coefficients
of both the upper right and lower left blocks satisfy
$$
v_-\equiv w_-\equiv 0\; {\rm mod}\; 4.
$$
This implies that the element $\nu$ belongs to $U^{\theta,2}_{K,\La}$.
Set $\theta_\pm:=(\theta_1\pm\theta_2)/2$.
The fixed point of $m_\nu$ in $\mc H_2$ is
$\theta=P^{-1}\theta_0=\mbox{\scriptsize 
$\left(\!\begin{array}{cc} \theta_{+}&\theta_{-}\\  
\theta_{-}&\theta_{+}\end{array}\!\right)$}.$
\vs

\noindent $(d2)$ We assume that $\eps_1=\eps_2=1$. 

We choose 
the new basis \eqref{eqnnewbas} to be  given by the basis change matrix 
$
P=
\mbox{
$\left(\!\begin{array}{cc} {\bf 2}&{\bf p}\\  
{\bf 0}&{\bf 1}
\end{array}\!\right)$}
$ 
with 
$
{\bf p}=
\mbox{\scriptsize 
$\left(\!\begin{array}{cc} 1\!&1\\  
1\!&\!-1\!
\end{array}\!\right)$}
$ 
so that
\begin{eqnarray*}
\label{eqnbaschaid2}
(e_1,e_2,-f_1,-f_2)&=&(2e_{0,1},\; 2e_{0,2},\;e_{0,1}\!+\!e_{0,2}\!-\!f_{0,1},\;e_{0,1}\!-e_{0,2}\!-\!f_{0,2}).
\end{eqnarray*}
This basis is symplectic for the form $\tfrac12\om_0$.

The matrix $m_{\nu}$ of multiplication by $\nu^{-1}$ in 
the basis $e_1$, $e_2$, $-f_1$, $-f_2$  
is 
\begin{equation} 
\label{eqnmnumod2}
m_{\nu}\!=\!P^{-1}m_{0,\nu}P=\mbox{\scriptsize 
$\left(\!\begin{array}{cccc} \!u_1\!-\!w_1\!&-w_2&(v_1\!-\!w_1\!-\!w_2)/2&\!(u_1\!-\!u_2\!-\!w_1\!+\!w_2)/2\!\\  
-w_1&\!u_2\!+\!w_2\!&\!(u_2\!-\!u_1\!+\!w_2\!-\!w_1)/2\!&(v_2\!-\!w_1-\!w_2)/2\\
2w_1&0&u_1+w_1&w_1\\  
0&2 w_2&w_2&u_2-w_2\end{array}\!\right)$}.
\end{equation} 
The congruence conditions \eqref{eqnujvjwj} imply that 
$m_\nu\equiv {\bf 1}$ mod $2$ and that the diagonal coefficients
of both the upper right and lower left blocks satisfy
$$
(v_1\!-\!w_1\!-\!w_2)/2\equiv (v_2\!-\!w_1-\!w_2)/2
\equiv 2w_1 
\equiv 2 w_2
\equiv 0\; {\rm mod}\; 4.
$$
This implies that the element $\nu$ belongs to $U^{\theta,2}_{K,\La}$.
The fixed point of $m_\nu$ in $\mc H_2$ is
$\theta=P^{-1}\theta_0=\frac12\mbox{\scriptsize 
$\left(\!\begin{array}{cc} \!\theta_{1}\!-\! 1\!&-1\\  
-1&\!\theta_{2}\! +\! 1\!\end{array}\!\right)$},$\vs 

\noindent $(e)$  Hence, in all cases, by Theorem \ref{thmmaimul}, one has a $d$-critical value
\begin{eqnarray*}
\la_\nu&=& \ka_\nu d^{1/2}N_\Phi(\nu)^{1/2}= 
\ka_\nu (d_1\nu_1)^{1/2}(d_2\nu_2)^{1/2}\\
&=&\ka'_\nu(\sqrt{a_1}+i\eps_1\sqrt{b_1})(\sqrt{a_2}+i\eps_2\sqrt{b_2}),
\end{eqnarray*}
with ${\ka_\nu'}^4=1$.
\vs

\noindent $(f)$ 
Since by Proposition \ref{procrival}, the ratio $(\la_\nu -1)/2$ is an algebraic integer,
one must have $\ka'_\nu=\pm 1$.
The precise sign $\ka'_\nu$ 
will be computed in Section \ref{secsigproima}.
\end{proof}

%63
\subsection{Quartic CM fields}
\label{secquafie}

In this Section we find new critical values by using quartic CM algebras (Proposition \ref{prorarbrc1}). 
Among them we will find the last four critical values
that were pointed out in Section \ref{secintexp}. 

\bp
\label{prorarbrc1}
Let $d=a\!+\! b\!+\!c$ be positive integers with $b^2\! -\! 4ac>0$. \\
${\bf (A)}$ Assume that $a\equiv 1$ {\rm mod} $4$ and
$b\!\equiv\! c\!\equiv\! (0\, {\rm or}\, 1)$ {\rm mod} $4$. 

Then the sum
$
\la= \sqrt{a}\!+\!\sqrt{c}\! +\! i\sqrt{b\!-\! 2 \sqrt{ac}} 
$ 
is  $d$-critical.\\
${\bf (B)}$ Assume that $a\equiv 3$ {\rm mod} $4$ and
$b\!\equiv\! c\!\equiv\! (0\, {\rm or}\, 3)$ {\rm mod} $4$. 

Then the sum
$
\la= \sqrt{b\!-\! 2 \sqrt{ac}} \! +\! i\sqrt{a}\!+\!i\sqrt{c}
$ 
is  $d$-critical.
\ep

For instance\\ 
$\star$ $\la=1\!+\sqrt{5} + i\sqrt{5\!-\!2 \sqrt{5}}$ is $11$-critical.\\
$\star$  $\la=1\!+\sqrt{5} + i\sqrt{9\!-\!2 \sqrt{5}}$
is $15$-critical.\\
$\star$ $\la=2\sqrt{2\!-\! \sqrt{3}}+ i\sqrt{3}\!+2i $ is $15$-critical.\\
$\star$ $\la=1\!+2\sqrt{2}+2i\sqrt{2\!-\!\sqrt{2}}$ 
is $17$-critical.\\ 
$\star$ $\la= \sqrt{11\!-\!2 \sqrt{21}}+ i\sqrt{3}\!+i\sqrt{7} $ is $21$-critical.

\br 
\label{remlamod2}
Note that $\la$ satisfies the condition $\la\equiv 1$ mod $2$ 
from Proposition \ref{procrival}.
Indeed, checking for instance the case when $a\equiv b\equiv c\equiv 1$ mod $4$,
one writes $\la=\la_1+\la_2+\la_3$ with $\la_1=\sqrt{a}$, $\la_2=\sqrt{c}$
and $\la_3=i\sqrt{b-2\sqrt{ac}}$, and one has
$\la_1\equiv\la_2\equiv\la_3\equiv 1$ mod $2$
because 
$\la_1^2\equiv\la_2^2\equiv\la_3^2\equiv 1$ mod $4$.
\er

The key remark that relates these values $\la$ and the reflex norm as in 
Theorem \ref{thmmaimul} is the following old-fashioned factorization.

\bl
\label{lemfaclam}
Let $d=a\!+\! b\!+\!c$ be positive integers with positive discriminant $\De:=b^2\! -\! 4ac>0$.
Then one has the equality
\begin{eqnarray*} 
\sqrt{a}\!+\!\sqrt{c}\! +\! i\sqrt{b\!-\! 2 \sqrt{ac}} 
&=& \sqrt{a}\;\left( 1+i\sqrt{\tfrac{b+\sqrt{\De}}{2a}}\right)
\;\left( 1-i\sqrt{\tfrac{b-\sqrt{\De}}{2a}}\right).
\end{eqnarray*} 
\el  

The proof of this factorization is left to the reader. For instance, one has\\
\centerline{$
1\!+\!\sqrt{5}\! +\! i\sqrt{9\!-\! 2 \sqrt{5}} 
= \left( 1+i\sqrt{\frac{9+\sqrt{61}}{2}}\,\right)
\;\left( 1-i\sqrt{\frac{9-\sqrt{61}}{2}}\,\right).$}

\begin{proof}[Proof of Proposition \ref{prorarbrc1}] We will prove simultaneously {\bf (A)} and {\bf (B)}.
The only difference between the two proofs will occur in the last step $(f)$ 
when computing the fourth root of unity $\ka_\nu$.
\vs 

\noindent $(a)$  We choose for CM algebra  the quartic $\m Q$-algebra $K=\m Q[\al]/(F(\al))$ 
where 
\begin{equation*}
\label{eqnpalabc}
F(\al)=a\al^4+b\al^2+c,
\end{equation*}
is a polynomial with  positive integral coefficients $a,b,c$ and with positive discriminant $\De=b^2-4ac$.
This CM algebra is a CM field if and only if the polynomial $F$ is irreducible.
We set 
$$
\mbox{
$\de:=2a\al^2+b\in K$} 
$$
so that $\de^2=\De$ and let $K_0$
be the  totally real subalgebra
$K_0:=\m Q[\de]$. One has $K_0=\m Q 1\oplus \m Q\de$ where $1$ is the unity element of $K$.
We introduce the $\m Q$-basis of $K$ 
$$
e_{0,1}:=1,\;\; e_{0,2}:=\de,\;\; f_{0,1}:=\al\de,\;\; f_{0,2}:=\al.
$$

We define the symplectic form $\om_0$ on $K$ by 
$$\textstyle
\om_0(x,x')={\rm Tr}_{K/\m Q}(\frac{x\ol{x}'}{4\al\de}).
$$ 
so that, by direct computation, one has $\om_0=f_{0,1}^*\wedge e_{0,1}^*+f_{0,2}^* \wedge e_{0,2}^*$. 
We introduce the  morphisms 
$\rho_\pm:K\rightarrow \m C$
defined by 
$$
\rho_+(\al)=\al_+:=i\sqrt{(b-\sqrt{\De})/2a}
\;\;{\rm and}\;\;
\rho_-(\al)=\al_-:=-i\sqrt{(b+\sqrt{\De})/2a}.
$$
The signs have been chosen carefully so that  
both ${\rm Im}(\rho_\pm(\al\de))>0$. Therefore the CM type defined by 
\eqref{eqnimrhot}  is $\Ph=\{\rho_+,\rho_-\}$.
\vs 

\noindent $(b)$ We first introduce the lattice $\La_0=\m Z e_{0,1}\oplus \m Ze_{0,2}\oplus\m Z f_{0,1}\oplus \m Z f_{0,2}$. 
This lattice is autodual for the symplectic form $\om_0$.
\vs 

\noindent $(c)$ We set $t=\al$.
We choose the totally unitary element 
\begin{equation} 
\nu:=\frac{1+t}{1-t}=
1+\frac{2}{d}(a\al^3+a\al^2+(a+b)\al-c)\in K.
\end{equation}
In the basis $e_{0,1},e_{0,2},-f_{0,1},-f_{0,2}$, the matrix $m_{0,\nu}$ of multiplication by $\nu^{-1}$  is
\begin{equation} 
\label{eqnmulnuquar}
m_{0,\nu}=\frac{1}{d}\mbox{\scriptsize 
$\left(\!\begin{array}{cccc}a\!-\!c&\De&\De&-b\!-\!2c\\  
1&a\!-\!c&-b\!-\!2c&1\\
1&2a\!+\!b&a\!-\!c&1\\
2a\!+\!b&\De&\De&a\!-\!c\end{array}\!\right)$}.
\end{equation}
By Proposition \ref{prohsidsi}, there exists $\si_{0,1}$ and $\si_{0,2}$ in 
${\rm Sp}(g,\m Z)$ and a diagonal matrix ${\bf d}_0={\rm diag}(d_1,d_2)$ 
with $d_1|d_2$ such that the symplectic matrix $m_{0,\nu}$ 
can be written as 
%\begin{eqnarray*}
$m_{0,\nu}=\si_{0,1}\;
\mbox{\small\scriptsize
$\left(\!
\begin{array}{cc} 
\!{}^t{\bf d}^{-1}_0\!&{\bf 0}\\  
{\bf 0}&{\bf d}_0
\end{array}
\!\right)$}
\;\si_{0,2}.$
%\end{eqnarray*} 

Note that, by the above computation, the matrix $dm_{0,\nu}$ has integer coefficients, the gcd of its coefficients is equal to $1$
and the gcd of all its  $2\times 2$ minor determinants is equal to $d$.
Therefore one can choose $d_1=1$ and $d_2=d$.
This proves that the group 
$G_{0,\nu}=\La_0/(\La_0\cap\nu\La_0)$ is isomorphic to $\m Z/d\m Z$. 
\vs 

\noindent $(d)$ The matrix 
$m_{0,\nu}$ is symplectic 
but does not  belong to the symplectic theta subgroup  ${\rm Sp}^{\theta,2}_{2,\m Q}$ 
because the diagonal of the lower-left block
of  $dm_{0,\nu}$, which is $(1,\De)$  
is not equal to $0$ mod $4$.
This is why we again have to choose another basis 
\begin{align}
\label{eqnnewbasbis}
(e_1,e_2,-f_1,-f_2)&=
(e_{0,1},e_{0,2},-f_{0,1},-f_{0,2}) \, P
\end{align}
and another lattice $\La=\m Z e_1\oplus \m Ze_2\oplus \m Z f_1 \oplus \m Z f_2$.
This lattice $\La$ will be sublattice of $\La_0$ of index $16$
that will be autodual for $\tfrac14 \om_0$, 
and such that $\nu$ belongs to $U^{\theta,2}_{K,\La} $,
Therefore, since $d$ is odd,
one will still have
$$
G_\nu=\La/(\La\cap\nu\La)\simeq \La_0/(\La_0\cap\nu\La_0)=G_{0,\nu}\simeq \m Z/d\m Z,
$$ 
and $\nu$ will belong to the unitary theta subgroup  $U^{\theta,2}_{K,\La}$ of level $2$.

To explain precisely  this choice of $\La$ 
we again need to distinguish two subcases 
that we call $(d1)$ and $(d2)$.
\vs 

\noindent $(d1)$ We assume that $a$ is odd and $b\equiv c\equiv 0$ mod $4$. 

We choose 
the new  basis \eqref{eqnnewbasbis} to be  given by the  basis change  matrix
$
P=
\mbox{
$\left(\!\begin{array}{cc} \!2{\bf p}\!&{\bf 0}\\  
{\bf 0}&2\, {}^t{\bf p}^{-1}\!\!
\end{array}\!\right)$}
$ 
with 
$
{\bf p}=
\mbox{\scriptsize 
$\left(\!\begin{array}{cc} 2\!&b/2\\  
0\!&\!1/2\!
\end{array}\!\right)$},
$ so that
$$
(e_1,e_2,-f_1,-f_2)=(4e_{0,1},\; b\, e_{0,1}+e_{0,2},\;-f_{0,1}+b\,f_{0,2},\;-4f_{0,2}).
$$
The matrix $m_\nu$ of multiplication by $\nu^{-1}$ 
in the basis $e_1$, $e_2$, $-f_1$, $-f_2$ is
\begin{equation}
\label{eqnmnuqd1}
m_{\nu}=P^{-1}m_{0,\nu}P=\frac1d\mbox{\scriptsize 
$\left(\!\begin{array}{cccc} 
2a\!-\!d&-ac&\!b^2\!+\!bc\!-\!ac\!&\!-2b\!-\!2c\!\\  
4&d\!-\! 2c&-2b\!-\!2c&4\\
4&2a\!+\!2b&2a\!-\!d&4\\
2a\!+\!2b&\!b^2\!+\!ab\!-\!ac\!&-ac&d\!-\! 2c\end{array}\!\right)$}.
\end{equation}
The element $\nu$ belongs to $U^{\theta,2}_{K,\La}$
because  $m_\nu\equiv{\bf 1}$ mod $2$, and  
both the diagonal of the upper-right block
and of the  lower-left block of  the matrix $d\, m_{\nu}$ satisfy
$$
b^2\!+\!bc\!-\!ac\equiv 4\equiv 4\equiv  b^2\!+\!ab\!-\!ac\equiv 0\;\;
{\rm mod}\;\; 4.
$$ 

\noindent $(d2)$ We assume that 
$a$ is odd and 
$b\equiv c\equiv a$ mod $4$. 

We choose 
the new basis \eqref{eqnnewbasbis} to be  given by the
basis change matrix  
$
P=
\mbox{
$\left(\!\begin{array}{cc} \!2{\bf p}\!&{\bf p}\\  
{\bf 0}&\!2\,{}^t{\bf p}^{-1}\!\!
\end{array}\!\right)$}
$ 
with 
$
{\bf p}=
\mbox{\scriptsize 
$\left(\!\begin{array}{cc} 2\!&b/2\\  
0\!&\!1/2\!
\end{array}\!\right)$},
$ 
so that
$$
(e_1,e_2,-f_1,-f_2)=(4e_{0,1},\;
2b\,e_{0,1}+2e_{0,2},\;2e_{0,1}-f_{0,1}+b\,f_{0,2},\; b\,e_{0,1}+e_2-2f_{0,2}).
$$
The matrix $m_\nu$ of multiplication by $\nu^{-1}$ in the basis $e_1$, $e_2$, $-f_1$, $-f_2$ is
\begin{equation}
\label{eqnmnuqd2}
m_\nu =P^{-1}m_{0,\nu}P= \frac1d\mbox{\scriptsize 
$\left(\!\begin{array}{cccc} 
2a\!-\!d\!-2\!&\!-2ac\!-\!2a\!-\!2b&
\!b^2\!+\!bc\!-\!ac\!-\!1\!&\!-ac\!\!-d\!-\!b\!-\!1\!\\  
\!2\!-\!2a\!-\!2b&\!d\!-\! 2c\!-\!2f\!&
\! ac\!-d\!-\!b\!+\!1\!&1\!-\! f\\
4&4a\!+\!4b&
2a\!-\!d\!+\! 2&2a\!+\!2b\!+\!2\\
4a\!+\!4b&\!4f\!&
2a\!+\!2b\!-\!2ac&\!d\!-\! 2c\!+\!2f\!\end{array}\!\right)$},
\end{equation}
where  $f\!:=\!b^2\!+\!ab\!-\!ac$.
The element $\nu$ belongs to $U^{\theta,2}_{K,\La}$
because $m_\nu\equiv{\bf 1}$ mod $2$, and  both the diagonal of the upper-right 
and  lower-left block of  the matrix $d\, m_{\nu}$ satisfy
$$
b^2\!+\!bc\!-\!ac\! -1\equiv 4\equiv 4\equiv 4f\equiv 0\;\;
{\rm mod}\;\; 4.
$$

\noindent $(e)$ We compute the critical value $\la_\nu$ 
given by Theorem \ref{thmmaimul}.
We first notice that $N_{K/\m Q}(1+\al)=\frac{d}{a}$ 
because this norm is the sum of the coefficients of the polynomial $F/a$, and hence that
\begin{eqnarray*}
N_\Phi(\nu)&=&
\frac{N_\Phi(1+\al)}{N_\Phi(1-\al)}=\frac{N_\Phi(1+\al)^2}{N_{K/\m Q}(1+\al)}
=\frac{a}{d}\, N_\Phi(1+\al)^2.
\end{eqnarray*}

Therefore by Theorem \ref{thmmaimul}  combined with Lemma \ref{lemfaclam}, one gets 
\begin{eqnarray*}
\la_\nu&=& 
\ka_\nu\, d^{1/2} N_\Phi(\nu)^{1/2}=
\ka_\nu\, a^{1/2} N_\Phi(1+\al)\\
&=&
\ka_\nu\, a^{1/2} (1+\al_+)(1+\al_-)= 
\ka_\nu\,(\sqrt{a}\!+\!\sqrt{c}\! -\! i\sqrt{b\!-\! 2 \sqrt{ac}}  ).
\end{eqnarray*}

\noindent $(f)$ 
Since by Proposition \ref{procrival},  
the ratio $(\la-1)/2$ is an algebraic integer,
arguing as in Remark \ref{remlamod2},
one must have
\begin{eqnarray}
\label{eqnlanupm1}
\la_\nu=\pm (\sqrt{a}\!+\!\sqrt{c}\! -\! i\sqrt{b\!-\! 2 \sqrt{ac}} )
&{\rm when}&
a\equiv 1 \;{\rm mod}\: 4,\\
\label{eqnlanupm3}
\la_\nu=\pm (\sqrt{b\!-\! 2 \sqrt{ac}}\! +\! i\sqrt{a}\!+\!i\sqrt{c} )
&{\rm when}&
a\equiv 3 \;{\rm mod}\: 4.
\end{eqnarray}
The precise sign $\pm$ will be computed  in Section \ref{secsigquafie}.
\end{proof}

Note that the complex conjugation 
will give another critical value with opposite imaginary part. 
However, it is important to remember for latter use that the critical values
$\la_\nu$ we have constructed in this section  satisfy:

\mbox{ }\hfill 
$\frac{{\rm Im}(\la_\nu)}{{\rm Re}(\la_\nu)}<0$ 
in Case  \eqref{eqnlanupm1} and $\frac{{\rm Im}(\la_\nu)}{{\rm Re}(\la_\nu)}>0$  in Case \eqref{eqnlanupm3}.\hfill{ }

%70
\section{Ideals in CM fields}
\label{secabemok}

In this chapter, we come back to the general CM abelian varieties, but we  specialize our 
Theorem \ref{thmmaimul} 
to the case when the lattice $\La$ is a fractional ideal $\g m$ of $K$, or, equivalently, 
the abelian variety $A=\m C^g/\Phi(\La)$
has multiplication by $\mc O_K$ (Theorem \ref{thmmaimul2}).
In this case 
a result of Taniyama-Shimura 
tells us when
this abelian variety $A$ is principally polarized
(Fact \ref{facshitan}).

An important case is when $K$ is the cyclotomic field $\m Q[\zeta_n]$ (Proposition \ref{promaimul3}). In this case, 
for all CM type $\Phi$ of $K$ 
one can find an ideal $\g m$ of $K$ 
such that the abelian variety $\m C^g/\Phi(\g m)$ is principally polarized (Corollary \ref{corcycfie}). 

In the last Section \ref{secfibluc}, we give another concrete
family of $d$-critical values where $d=L_{n}$ is the $n^{\rm th}$
Lucas number with $n\geq 5$ prime (Corollary \ref{corfibonacci}).

%71
\subsection{Using class field theory}
\label{secshitan}

When $K$ is a CM number field, we would  like to find $t_0\in \mc I_K$ and an ideal $\g m$ of $K$  
such that the symplectic form
$\om_{t_0}$
restricted to $\g m$
takes values in $\m Z$ and has determinant equal to $1$,
because this means that the polarized  abelian variety 
$(\m C^g/\Phi_{t_0}(\g m),\om_{t_0})$
is principally polarized. 
Shimura and Taniyama explained when this is possible: in concrete words, they show that,\\  
$(A)$  such a construction is  possible, 
if we allow the CM type $\Phi$ to vary,\\ 
$(B)$ if the extension  $K/K_0$ is ramified, 
one can  prescribe the CM type $\Phi$,\\
$(C)$ if the extension $K/K_0$ is unramified, they describe the possible CM types.\\ 
Here is their precise result in which 
$\mc D_K$ and $\mc D_{K/K_0}$ are the differents.
\bfa
\label{facshitan}
{\bf (Shimura, Taniyama)} Let $K$ be a CM field.
\\
{\bf (A)} , There exist
an ideal $\g m$ of $\mc O_K$
and an imaginary element $t_0\in \mc I_K$
such that
the polarized abelian variety
$(\m C^g/\Phi_{t_0}(\g m),\om_{t_0})$ is principal.
\\
{\bf (B)} If
$\mc D_{K/K_0}\neq \mc O_K$, for all CM type $\Phi$
of $K$, there exist an ideal $\g m$ of $\mc O_K$
and $t_0\in \mc I_K$ such that $\Phi_{t_0}=\Phi$ 
and
the polarized abelian variety
$(\m C^g/\Phi_{t_0}(\g m),\om_{t_0})$ is principal.
\\
{\bf (C)} If $\mc D_{K/K_0}= \mc O_K$, 
let $t_{00}\in \mc I_K$  and 
$\g n:=t_{00}\mc D_K\cap {K_0}$. 
One has 
$(i)\Leftrightarrow(ii)$:\\
$(i)$ There exist an ideal $\g m$ of $\mc O_K$
and $t_0\in \mc I_K$ such that 
$\Phi_{t_0}=\Phi_{t_{00}}$ 
and
the polarized abelian variety
$(\m C^g/\Phi_{t_0}(\g m),\om_{t_0})$ is principal.\\
$(ii)$ The element ${\rm Art}_{K/K_0}(\g n)\in {\rm Gal}(K/K_0)$ is trivial.
\efa

Condition $(ii)$ involves the  image 
of the ideal $\g n$ of $K_0$ by the Artin map and its restriction
${\rm Art}_{K/K_0}(\g n)$ 
to the unramified extension $K$ of $K_0$. This restriction belongs to
the Galois group ${\rm Gal}(K/K_0)$ which has order $2$.

Condition $(ii)$ is automatically satisfied 
when the ideal $\g n$ is principal. 
We will see in Corollary \ref{corcycfie} that this is always the case for cyclotomic fields.

Note that for the cyclotomic fields $K=\m Q[\zeta_n]$, both cases $(B)$ and $(C)$ can occur (see Lemma \ref{lemcycfie}),

\begin{proof}[Proof of Fact \ref{facshitan}]
This is \cite[Proposition 1]{Shimura77}, but it will be useful to point out how class field theory is used in the proof. 
We will discuss successively  the points $(C)$,  $(B)$ and $(A)$.

We first note that 
the relative different $\mc D_{K/K_0}$ is generated by imaginary elements. Therefore,
given a nonzero $t_{00}$ in $\mc I_K$, there exists a fractional ideal $\g n$
of ${K_0}$ such that
$t_{00}^{-1}\mc D_K=\g n^{-1}\mc O_K$.
We want to know if there exist a totally positive 
element $s_0$ in $K_0$
and an ideal $\g m$ of $\mc O_K$ which is autodual for $\om_{t_0}$ where 
$t_0:=s_0t_{00}$. This means that
\begin{eqnarray*}
\label{eqntoodk}
t_0^{-1} \mc D_K\g m\ol{\g m}=\mc O_K
\;\;\mbox{\rm or, equivalently,}\;\;
\g m\ol{\g m}=s_0\g n\mc O_K.
\end{eqnarray*}
This exactly means that the class $[\g n]$ in the narrow class group 
 $\mc C\ell^+(K_0)$ is 
the image by the norm map $N_{K/K_0}$  of a class 
$[\g m]$ of the narrow class group 
$\mc C\ell^+(K)$.
Let $L^+_0$ be the maximal unramified (at all finite places) abelian extension of $K_0$
and $L^+$ the maximal unramified abelian extension of $K$.
By class field theory, the Artin maps ${\rm Art}_{K}$
and ${\rm Art}_{K_0}$ give a commuta\-tive diagram
with horizontal isomorphisms as in \cite[p. 400]{Washington97},
\begin{equation}
 \begin{CD}
\mc C\ell^+(K) @>{\rm Art}_{K}>> {\rm Gal}(L^+/K)\\
@VN_{K/K_0}VV @VV{\rm Res}_{L^+_0}V\\
\mc C\ell^+(K_0) @>{\rm Art}_{K_0}>> {\rm Gal}(L^+_0/K_0)
\end{CD}
\end{equation}
The question is whether the element 
$\si:={\rm Art}_{K_0}([\g n])\in {\rm Gal}(L^+_0/K_0)$ is in the image of the restriction map
${\rm Res}_{L^+_0}$.
Equivalently the question is whether
the restriction of $\si$ to $L_0^+\cap K$ is trivial.\vs 

In case $(C)$, the 
field $K$ is included in $L_0^+$
and the criterion is the triviality of $\si$ 
in $K$.
\vs 

In case $(B)$, 
the intersection $L_0^+\cap K$ is equal to $K_0$. Therefore the criterion is always true.
\vs 

In case $(A)$, we do not worry about the 
total positivity of $s_0$ therefore we 
introduce the Hilbert class field $L_0$ of $K_0$ which is also the maximal totally real subfield 
of $L_0^+$, and the Hilbert class field $L$ of $K$ which is equal to $L^+$. 
Class field theory gives a similar commutative diagram, 
with horizontal isomorphisms ${\rm Art}_{K}$
and ${\rm Art}_{K_0}$,
\begin{equation}
 \begin{CD}
\mc C\ell(K) @>{\rm Art}_{K}>> {\rm Gal}(L/K)\\
@VN_{K/K_0}VV @VV{\rm Res}_{L_0}V\\
\mc C\ell(K_0) @>{\rm Art}_{K_0}>> {\rm Gal}(L_0/K_0)
\end{CD}
\end{equation}
Since $L_0$ is totally real, one has $K\cap L_0=K_0$ and 
the restriction map
${\rm Res}_{L_0}$ is onto. The norm map $N_{K/K_0}$ is also onto.
\end{proof}

%72
\subsection{Critical values and ideals in CM field}
\label{seccricom2}

We now give a consequence  of Theorem \ref{thmmaimul} 
applied to the case where the lattice $\La$ 
is an ideal $\g m$ of $\mc O_K$, combined with 
Shimura Taniyama's Fact \ref{facshitan}.
\vs 

We begin by an explicit construction of elements of 
the group 
$U^{\theta,2}_{K,\g m}$ 
when $\g m$ is an  ideal.
This construction does not depend on $\g m$.

\bl
\label{lemitkl2} Let $K$ be a CM number field,
$\om_{t_0}$ a symplectic form on $K$ as in 
\eqref{eqnomotrk}
and $\g m$ be an autodual fractional ideal of $K$.\\
$(a)$ 
For all $s$ in $K$ with denominator prime to $2$,
the difference 
$t:=s-\ol{s}$ belongs to $\mc I^{\theta,2}_{K,\g m}$.\\
$(b)$ If moreover $N_{K/\m Q}(1\!+\!t)$
has an odd numerator, then
$\nu:=\frac{1+t}{1-t}$ is in $U^{\theta,2}_{K,\g m}$.
\el

\begin{proof} $(a)$ 
Indeed, one has 
$\om_{t_0}(tx,x)=2\,\om_{t_0}(sx,x)\in 2\m Z_{(2)}$, 
for all $x$  in $\g m$.

$(b)$ Since $N_{K/\m Q}(1-t)=N_{K/\m Q}(1+t)$, 
this follows from Point $(a)$ and  from Lemma \ref{lemimauni}.
\end{proof}

The three parts of the following Theorem \ref{thmmaimul2} correspond to the three parts of Shimura Taniyama's Fact \ref{facshitan}.

\bt 
\label{thmmaimul2} 
Let $K$ be a CM field and $K_0:=K\cap \m R$.
Let 
$\nu=\mu/\ol{\mu}$ with $\mu=1+s-\ol{s}$
where $s\in K$ has denominator prime to $2$ and 
$N_{K/\m Q}(\mu)$ has odd numerator.
Set $G_\nu:=\mc O_K/(\mc O_K\cap \nu\mc O_K)$, $d_\nu:=|G_\nu|$ and $\la_{0,\nu}=d_\nu^{1/2} N_\Ph(\nu)^{1/2}$.
\\
{\bf (A)} There exist a CM type $\Phi$ of $K$
and a critical value $\la_\nu=\ka_\nu\la_{0,\nu}$
on the group $G_\nu$ with $\ka_\nu^4=1$.\\
{\bf (B)} If $\mc D_{K/K_0}\neq\mc O_K$, then 
for all CM types $\Phi$ of $K$, there exists
a critical value $\la_\nu=\ka_\nu\la_{0,\nu}$
on the group $G_\nu$ with $\ka_\nu^4=1$.\\
{\bf (C)} If $\mc D_{K/K_0}=\mc O_K$, 
let $t_{00}\in \mc I_K$ and set $\Phi=\Phi_{t_{00}}$ and  $\g n:=t_{00}\mc D_K\cap {K_0}$. Assume that the element ${\rm Art}_{K/K_0}(\g n)\in {\rm Gal}(K/K_0)$ is trivial, then  
there exists
a critical value $\la_\nu=\ka_\nu\la_{0,\nu}$
on the group $G_\nu$ with $\ka_\nu^4=1$.
\et

\bc 
\label{cormaimul2}
Moreover, in all these cases, if $s$ belongs to $\mc O_K$,
one has
\begin{equation}
\la_\nu=\ka_\nu N_\Phi(\mu)
\;\;{\it and}\;\; 
G_\nu=\mc O_K/\mu\mc O_K.
\end{equation}
\ec

We recall that the expression {\it $s$ has denominator prime to $2$},
means that the prime divisors of the ideal $2\mc O_K$
never occur in the prime decomposition of the fractional ideal $s\mc O_K$ 
with a negative power.

\begin{proof}[Proof of Theorem \ref{thmmaimul2}]
By Shimura Taniyama Fact \ref{facshitan}, we know that in these three cases, there exists
$t_0\in \mc I_K$ with $\Ph_{t_0}=\Ph$ and a fractional ideal $\g m$ of $K$ such 
that the polarized abelian variety
$(\m C^g/\Phi(\g m),\om_{t_0})$ is principal.
We will choose $\La=\g m$. 
By Lemma \ref{lemitkl2}, our 
unitary element  $\nu\in K$ 
belongs to the unitary theta subgroup
$U^{\theta,2}_{K,\g m} $ of level $2$.

Therefore Theorem \ref{thmmaimul} tells us
that $\la_\nu$ is a critical value on the group 
$G_\nu=\g m/(\g m\cap \nu\g m)$.
Since we are dealing with a finite quotient of $\g m$, 
we can replace $\g m$ by $\mc O_K$ so that
$
G_\nu\simeq\mc O_K/(\mc O_K\cap \nu\mc O_K).
$
\end{proof}

\begin{proof}[Proof of  Corollary \ref{cormaimul2}]
When $\mu$ is in $\mc O_K$, since $\mu+\ol\mu=2$ and $N_{K/\m Q}(\mu)$ is odd, 
the integers $\mu$ and  $\ol\mu$ are relatively prime and one has $\mu\mc O_K\cap \ol{\mu}\mc O_K=\mu\ol{\mu}\mc O_K$ so that
$G_\nu\simeq
\mc O_K/\ol{\mu}^{-1}(\mu\mc O_K\cap \ol{\mu}\mc O_K)=
\mc O_K/\mu\mc O_K
$ 
as required.
\end{proof}

\br 
\label{remgnucyc}
Note that in Theorem \ref{thmmaimul2},
it is easy to determine the structure of the  abelian group 
$G_\nu$ by decomposing the ideal 
$\mc J_\nu:=\mc O_K\cap \nu\mc O_K$ as a product 
$\prod_j\g p_j^{n_j}$ of prime ideals 
$\g p_j$ of $\mc O_K$ so that 
$G_\nu\simeq \prod_j\mc O_K/\g p_j^{n_j}$.
In particular, the group $G_\nu$ is cyclic 
if and only if the prime ideals 
$\g p_j$ that divide $\mc J_\nu$ are over different primes $p_j$ of $\m Z$
and have inertia degree equal to $1$.
\er

%73
\subsection{Critical values and cyclotomic fields}
\label{seccycfie}

In this section we apply Theorem \ref{thmmaimul2} to the cyclotomic fields $\m Q[\zeta_n]$. 
\vs 

We first recall  notation related to cyclotomic fields.
Let $\zeta_n$ be the primitive $n^{\rm th}$-root of unity $\zeta_n:=e^{2i\pi/n}$, let $\Ph_n$ be the $n^{\rm th}$
cyclotomic polynomial, let $K_n:=\m Q[\zeta_n]$ be the cyclotomic field, let $\ph(n)$ be the Euler number 
$\ph(n):=[K_n:\m Q]=d^\circ \Phi_n$, and let $K_{n,0}:=\m Q[\zeta_n+\zeta_n^{-1}]$ be the real cyclotomic field.
Their rings of integers are $\mc O_{K_{n}}=\m Z[\zeta_n]$
and $\mc O_{K_{n,0}}=\m Z[\zeta_n+\zeta_n ^{-1}]$.

\bl
\label{lemcycfie} Let $n\not\equiv 2$ mod $4$ and $K_n=\m Q[\zeta_n]$.\\ 
$(a)$ The different $\mc D_{K_n}$ is a principal ideal. 
More precisely, there exists an imaginary 
element $t_n\in \mc O_{K_n}$ 
such that $\mc D_{K_n}=t_n\mc O_{K_n}$.\\
$(b)$ The extension $K_n/K_{n,0}$ is ramified if and only if $n$ is a prime power.
\el

Since $\m Q[\zeta_{m}]=\m Q[\zeta_{2m}]$ for $m$
odd, the assumption on $n$ is innocuous.

\begin{proof} This is classical. We just sketch the argument.

$(a)$ Let  $\xi_n=\prod_{p|n}(1-\zeta_p)$
where the product is over the prime divisors of $n$.
We first claim that $\mc D_{K_n}=n\xi_n^{-1}\mc O_{K_n}$.
We write $n=\prod p^{r_p}$. 
For all integers $k_p$ one has 
$$\textstyle
Tr_{K_n/\m Q}(\xi_n\prod_{p}\zeta_{p^{r_p}}^{k_p})
=\prod_pTr_{K_{p^{r_p}}/\m Q}((1-\zeta_p)\zeta_{p^{r_p}}^{k_p})
\equiv 0
\;\;\mbox{\rm mod}\;\; n
$$
By definition of the different, this proves the inclusion  $\mc D_{K_n}\subset n\xi_n^{-1}\mc O_{K_n}$.
By \cite[Proposition 2.7]{Washington97} the absolute value of the discriminant $d_{K_n}$
is equal to $n^{\ph(n)}/\prod_{p|n}p^{\ph(n)/(p-1)}$.
Since this number is also equal to $N_{K_n/\m Q}(n\xi_n^{-1})$,
this proves the equality $\mc D_{K_n}=n\xi_n^{-1}\mc O_{K_n}$.

To conclude, we need to find a unit $u_n\in \mc O_{K_n}^*$ 
such that $u_n\xi_n$ is imaginary. 

We first notice that for $p$ odd the 
ratio $u_p:=\frac{\zeta_p-\zeta_p^{-1}}{1-\zeta_p}$ is a unit
and that for $p=2$ the ratio $u_p:=\frac{2\zeta_4}{1-\zeta_2}=i$
is also a unit. Therefore we set $u'_n:=\prod_{p|n}u_p$,
so that the product $u'_n\xi_n$ is imaginary 
when $n$ has an odd number of prime factors,
and is real when $n$ has an even number of prime factors.

In the first case, we are done with $u_n=u'_n$.

In the second case, $n$ is not a prime power and by Lemma 
\ref{lemcycfie2} below, there exists an imaginary unit $u''_n$ in 
$\m Q[\zeta_n]$ and we are done with $u_n:=u'_nu''_n$.

$(b)$ We will not use this fact which is proven in 
\cite[Proposition 2.15]{Washington97}. It 
explains why we have to pay attention to the case $(C)$ in 
Theorem \ref{thmmaimul2}. 
\end{proof}

In the proof we have used the following lemma.

\bl
\label{lemcycfie2} Let $n\not\equiv 2$ mod $4$. If $n$ is not a prime power, then the ring $\m Z[\zeta_n]$
contains an imaginary unit.
\el

\begin{proof}
If $n$ is even, the ring $\m Z[\zeta_n]$ contains the imaginary unit $i=\zeta_4$. 

If $n$ is odd, it admits two distinct odd prime factors $p,q$
and $\m Q[\zeta_n]$ contains the imaginary unit $\zeta_{pq}-\zeta_{pq}^{-1}$.  
\end{proof}

Combining Lemma \ref{lemcycfie} and 
Shimura-Taniyama Fact \ref{facshitan}, we now deduce that
there exists a principally polarized abelian variety with multiplication by $\mc O_{K_n}$ and with 
any prescribed CM type.

\bc
\label{corcycfie}
 Let $K_n=\m Q[\zeta_n]$, let $g=\ph(n)/2$
 and let $t_{00}\in \mc I_{K_n}$ nonzero.\\ 
$(a)$ The fractional ideal $\g n:=t_{00}\mc D_{K_n}\cap {K_{n,0}}$ is principal.\\
$(b)$
There exist an ideal $\g m$ of $\mc O_{K_n}$
and $t_0\in \mc I_{K_n}$ such that $\Phi_{t_0}=\Phi_{t_{00}}$ 
and such that 
$(\m C^g/\Phi_{t_0}(\g m),\om_{t_0})$ is a principally polarized abelian variety.
\ec

\begin{proof}
$(a)$ According to Lemma \ref{lemcycfie}.$a$, the different $\mc D_{K_n}$ 
is a principal ideal of $\mc O_{K_n}$ generated by an imaginary element $t_n\in \mc I_{K_n}$.
Therefore, the element $t_nt_{00}$ is in $K_{n,0}$ and one has 
$$
\g n=t_nt_{00}\mc O_{K_{n,0}}.
$$
This proves that $\g n$ is a principal ideal of $K_{n,0}$.

$(b)$ This follows from   Fact \ref{facshitan}.C, 
since, by $(a)$, the ideal $\g n$ is principal.
\end{proof}

Thanks to Corollary \ref{corcycfie}, the statement of Theorem \ref{thmmaimul2} is cleaner
for the cyclotomic field. Here it is:

\bp
\label{promaimul3} 
Let $K_n:=\m Q[\zeta_n]$ and $\nu=\mu/\ol{\mu}$ with $\mu=1+s-\ol{s}$
where $s\in K_n$ has denominator prime to $2$ and 
$N_{K_n/\m Q}(\mu)$ has odd numerator.\\
Let  $G_\nu:=\mc O_{K_n}/(\mc O_{K_n}\cap \nu\mc O_{K_n})$
and $d_\nu=|G_\nu|$.\\
$(i)$ Then for all CM types $\Phi$ of $K_n$, there exists
a critical value\\ 
\hspace*{8em}$\la_\nu=\ka_\nu\,d_\nu^{1/2} N_\Ph(\nu)^{1/2}$ on $G_\nu$
with $\ka_\nu^4=1$.\\
$(ii)$ 
Moreover, if $s$ is in $\mc O_{K_n}$,
one has
\hspace*{1em}$\la_\nu=\ka_\nu N_\Phi(\mu)$
and $G_\nu=\mc O_{K_n}/\mu\mc O_{K_n}.$
\ep

\br 
Point $(ii)$  does not apply to $\mu=1+\ze_3-\ol{\ze_3}=1+i\sqrt{3}$ 
because the norm of this $\mu$ is even. Indeed in this case $G_\nu=\{1\}$.
\er

The group $G_\nu$ has order $d_\nu$ but is not always 
cyclic. In the next section, we will give 
a very simple example
where one can check that $G_\nu$ is cyclic.

\begin{proof}[Proof of Proposition \ref{promaimul3}]
$(i)$ We distinguish two cases. 

In case the extension $K_n/K_{n,0}$ is ramified, we apply Theorem \ref{thmmaimul2}.B.

In case the extension $K_n/K_{n,0}$ is unramified, we
apply
Theorem \ref{thmmaimul2}.C, remembering 
that, by Corollary \ref{corcycfie}, 
 the ideal 
$\g n:=t_{00}\mc D_{K_n}\cap K_{n,0}$ is principal,
for any $t_{00}$ in $\mc I_{K_n}$. 

$(ii)$ Apply Corollary \ref{cormaimul2}.
\end{proof}

%74
\subsection{Fibonacci and Lucas numbers}
\label{secfibluc}

As an application of Proposition \ref{promaimul3}, 
we construct in Corollary \ref{corfibonacci} new $d$-critical values 
where the integers $d$ 
are  Lucas numbers $L_{n}$ with $n$ prime.\vs 

These Lucas numbers are intimately related with the Fibonacci numbers.
We recall that the Fibonacci numbers $(F_n)_{n\geq 0}= (0,1,1,2,3,5,8,\ldots)$ are defined by their first two terms and the recurrence relation $F_{n+1}=F_n+F_{n-1}$.
They can also be defined as the coefficients of a matrix $n^{\rm th}$-power
\begin{eqnarray}
\label{eqnfnfnfn}
\mbox{\scriptsize 
$\left(\!\begin{array}{cc} 1&1\\  
1&0\end{array}\!\right)^n$}=\mbox{\scriptsize 
$\left(\!\begin{array}{cc} F_{n+1}&F_n\\  
F_n&F_{n-1}\end{array}\!\right)$}. 
\end{eqnarray}
The Lucas numbers
$(L_n)_{n\geq 0}= (2,1,3,4,7,11,18,\ldots)$ 
are also defined by their first two terms and the same recurrence relation $L_{n+1}=L_n+L_{n-1}$. 
Hence they are equal to the trace of the matrix
\eqref{eqnfnfnfn}: 
$L_n=F_{n+1}+F_{n-1}$.
Together they can  be defined by
\begin{eqnarray}
\label{eqnlnfnr5}
 (\tfrac{1+\sqrt{5}}{2})^n
&=& 
\tfrac{L_n+F_n\sqrt{5}}{2}.
\end{eqnarray}

We will need a few properties of Fibonacci and Lucas numbers.

\bl 
\label{lemfibonacci}
$(a)$ One has $F_n=\prod_{1\leq k< n/2}(1+4 \cos^2(\frac{k\pi}{n}))$, and \\
$(b)$ 
$F_n=\prod_{d|n} F'_{d}$ where $F'_n$ is the integer
$F'_n:=\prod_{\substack{1\leq k< n/2\\k\wedge n=1}}(1\!+\!4 \cos^2(\frac{k\pi}{n}))$. 
\el

\begin{proof}[Proof of Lemma \ref{lemfibonacci}] All this is classical. Here are a few hints.

$(a)$ The Fibonacci polynomials $P_n(X)$ are defined by 
$P_1=1$, $P_2=X$ and $P_{n+1}=XP_n+P_{n-1}$. They 
satisfy the equality $F_n=P_n(1)$ and
are related to the Chebychev polynomial so that the zeros of 
$P_n$ are the $2i\cos(\frac{k\pi}{n})$ for $1\leq k \leq n$.
 
$(b)$ These numbers $F'_n$ are both algebraic integers
and invariant by Galois. Hence they are integers.
\end{proof}

\bl 
\label{lemlucas} Assume that $m$ and $n$ are odd.\\
$(a)$ One has $L_n=\prod_{1\leq k< n/2}(1+4 \sin^2(\frac{k\pi}{n}))$, and \\
$(b)$ 
$L_n=\prod_{d|n} L'_{d}$ where $L'_n$ is the integer
$L'_n:=\prod_{\substack{1\leq k< n/2\\k\wedge n=1}}(1\!+\!4 \sin^2(\frac{k\pi}{n}))$. \\
$(c)$ These integers $L'_n$ are odd except  $L'_3=4$.
\el
For instance, one has 
$L'_5\!=\! 11$, $L'_7\!=\! 29$, $L'_9\!=\! 19$,$\ldots$.

\begin{proof}[Proof of Lemma \ref{lemlucas}]
$(a)$ and $(b)$ From Formula \eqref{eqnlnfnr5}, one gets the equality 
$
F_{2n}=F_nL_n\, .
$
Our claims  together with $L'_n=F'_{2n}$ follow from Lemma \ref{lemfibonacci}.

$(c)$ Just compute the sequence $(L_n\; {\rm mod}\; 8)_{n\geq 0}$
and notice that it has period $12$ with pattern 
$[2,1,3,4,7,3,2,5,7,4,3,7]\; {\rm mod}\; 8$. 
\end{proof}

The following lemma will be useful.

\bl
\label{lemquocyc}
For $n\geq 5$ odd, set $\mu_n:=1+\zeta_n-\zeta_n^{-1}$. 
Then the quotient ring 
$\m Z[\zeta_n]/\mu_n\m Z[\zeta_n]$ is cyclic and isomorphic to $\m Z/L'_n\m Z$.
\el 

\begin{proof}[Proof of Lemma \ref{lemquocyc}]
Denote by $x_n$  the image of $\zeta_n^{-1}$ in the quotient ring 
$A_n=\m Z[\zeta_n]/\mu_n\m Z[\zeta_n]$. 
By definition $L'_n$ is equal to the norm of $\mu_n$, that is $L'_n:=N_{K_n/\m Q}(\mu_n)$. Therefore the ring $A_n$ has order $L'_n$. 
The element $x_n$ is invertible in $A_n$ and satisfies the equality $x_n^2=x_n+1$.
Therefore by \eqref{eqnlnfnr5}, one has for all integer $k\geq 1$,
$$
2\, x_n^k=L_k+(2x_n-1)F_k.
$$ 
We apply this equality with $k=n$, remembering that,
since $L'_n$ divides $L_n$, one has $L_n\equiv 0$ in $A_n$. 
One gets
$$
(2x_n-1)F_n=2\;\;\mbox{\rm  in }\;\; A_n.
$$
Since, $n\neq 3$, by Lemma \ref{lemlucas}, $L'_n$ is odd, and the element $2$ is invertible in $A_n$.
Therefore $x_n$ is a multiple of $1$  in $A_n$ and the group 
$A_n$ is cyclic.
\end{proof}

\bp
\label{profibonacci}
For all $n\geq 5$ odd, and all signs 
$\eps_k=\pm 1$, the product
$\la'_n=\prod\limits_{\substack{ k<n/2\\k\wedge n=1}}
(1+2i \eps_k \sin(\frac{k\pi}{n}))$ 
or $-\la'_n$
is  $L'_n$-critical.
\ep

\br 
Note that $\la'_n$ satisfies the condition $\la'_n\equiv 1$ mod $2$ 
from Proposition \ref{procrival}.
Indeed, one computes, 
\begin{eqnarray*}
\label{eqnlanmod2}
\la'_n&\equiv&
\textstyle
\prod_{\substack{k<n/2\\k\wedge n=1}}(1+2\cos(\tfrac{k\pi}{n}))\equiv 
\prod_{\substack{k<n\\k\wedge 2n=1}}(1+2\cos(\tfrac{k\pi}{n}))\;\; {\rm mod}\;\; 2\,\mc O_{K_n}\\
\la'_n&\equiv & 
\textstyle
\prod_{\substack{k<n\\k\wedge 2n=1}}(3+4\cos^2(\tfrac{k\pi}{2n}))
\equiv F'_{2n}\equiv L'_{n}
\;\equiv\; 
1\;\; {\rm mod}\;\; 2\,\mc O_{K_n}.
\end{eqnarray*}
\er

\begin{proof}[Proof of proposition \ref{profibonacci}]
This is a consequence of Proposition \ref{promaimul3} combined with Lemma \ref{lemquocyc}.
We work with the cyclotomic field $K=K_n$ and choose the CM type $\Phi:=\{\rho_k\mid k\in (\m Z/n\m Z)^*\}$ where 
$\rho_k(\zeta_n)=\zeta_n^{\eps_k k}$, 
and choose $\mu=\mu_n:=1+\zeta_n-\zeta_n^{-1}$
so that one has 
$\la'_n:=N_\Phi(\mu_n)$.
According to Proposition \ref{promaimul3}, $\ka_n \la'_n$ is a critical value on the group
 $\mc O_{K_n}/\mu_n\mc O_{K_n}$ for a $4^{\rm th}$ root of unity $\ka_n$.
According to Lemma \ref{lemquocyc}, this group is cyclic of order 
$L'_n$.
Moreover, since $\la'_n\equiv 1$ mod $2$, by Proposition \ref{procrival}, one has $\ka_n=\pm 1$.
\end{proof}

The precise sign $\ka_n$ involved with $\la'_n$ 
could be computed thanks to the exact transformation
formula for the Riemann theta function, and 
the conclusion should be that the value $\la'_n$
-with a plus sign- is always $L'_n$-critical.

The following is an immediate corollary.

\bc
\label{corfibonacci}
For all $n\geq 5$ prime, 
and all  signs
$\eps_k=\pm 1$, the product\\
$\la_n=\prod\limits_{1\leq k< n/2}
(1+2i\eps_k \sin(\frac{k\pi}{n}))$ 
or $-\la_n$
is  $L_n$-critical.
\ec

\br
Assume now that the integer $n$ is coprime to $6$ but is not a prime number.
Since $\la_n$ is the product of $\la'_d$ 
for $d$ divisors of $n$, 
the value $\pm \la_n$ is critical on the group $G_n$ product of the cyclic groups $\m Z/L'_d\m Z$. 
The product of these integers $L'_d$ is equal to $L_n$.
but, since these integers $L'_d$ are not always coprime, 
this group $G_n$ is not always the cyclic group $\m Z/L_n\m Z$.
\er

To be concrete, we give three examples, 
the first ones being the conjugates of 
$\la_5=1\!+\!\sqrt{5}+i\sqrt{5\!-\!2\sqrt{5}}$ 
that already occured in \cite[Section 1.5]{CSAGI}, 
and that was also part of Proposition \ref{prorarbrc1}~:
\vs 

\noindent
$\star$ $n=5$ : $\la_5:=(1\pm 2i\sin(\frac{\pi}{5}))(1\pm 2i\sin(\frac{2\pi}{5}))$ : $11$-critical.\\
$\star$ $n=7$ : $\la_7:=(1\pm 2i\sin(\frac{\pi}{7}))(1\pm 2i\sin(\frac{2\pi}{7}))(1\pm 2i\sin(\frac{3\pi}{7}))$ : $29$-critical.\\
$\star$ $n=9$ : $\la'_9:=(1\pm 2i\sin(\frac{\pi}{9}))(1\pm 2i\sin(\frac{2\pi}{9}))(1\pm 2i\sin(\frac{4\pi}{9}))$ : $19$-critical.

%80
\section{On the sign of the critical values}
\label{secsigcri}

This aim of this chapter is to explain how to prove that the values 
$\la$ given in the examples of Chapter \ref{seclisexa} are indeed 
$d$-critical. 
The arguments given in Chapter \ref{seclisexa} prove that either $\la$ or $-\la$ is $d$-critical. We will complete here these arguments and prove that $\la$ is indeed  $d$-critical.
This will rely on extra computations together with a formula
due to Stark for the theta cocycle $j(\si,\tau)$
that we give in Section \ref{secsigtrafor}.
Since this sign is  often equal to the plus sign, these computations 
give a
cleaner statement for our examples in Chapter \ref{seclisexa}.

%81
\subsection{Square root of determinant of Riemann matrices}
\label{secsqadetrie}

We begin by two simple lemmas on the determinant of a Riemann matrix.
For $\tau$ in $\mc H_g$, we define the function 
\begin{equation}
\label{eqnhtau}
h(\tau)=\det(\tfrac{\tau}{i})^{1/2}
\end{equation} 
as the continuous square root of 
$\det(\frac{\tau}{i})$ such that $h(i{\bf 1}_g)=1$.

\bl 
\label{lemsqadetrie} 
Let $\tau$ in $\mc H_g$.\\ 
$(a)$ One has $\det_\m C(\tau/i)\neq 0$.\\
$(b)$ All the eigenvalues $\la$ of $\tau/i$ have positive real part:\; 
${\rm Re}(\la)>0$.\\
$(c)$ The function $h(\tau)$ is the product of the square roots $\la^{1/2}$ of the eigenvalues $\la$ of $\tau/i$
with $\arg(\la^{1/2})\in (-\frac{\pi}{4},\frac{\pi}{4})$.\\
$(d)$ If $g=2$, the function $h(\tau)$ has positive real part:\; ${\rm Re}(h(\tau))>0$.
\el

\begin{proof} 
$(a)$ 
If $z\in \m C^g$ is in the kernel of $\tau$, 
one has ${\rm Im}({}^t\ol{z}\tau z)=0$. Since $\tau$ is symmetric, this can be rewritten as
${}^t\ol{z}{\rm Im}(\tau) z=0$. Since ${\rm Im}(\tau)$ is positive this gives $z=0$. And hence $\tau$ is invertible.

$(b)$ Indeed, if ${\rm Re}(\la)\leq 0$, the matrix $\tau-i\la{\bf 1}_g$ 
is in $\mc H_g$ and $\det(\tau-i\la{\bf 1}_g)\neq 0$.

$(c)$ This follows from $(b)$ by analytic continuation.

$(d)$ This follows from $(c)$.
\end{proof}

When $S$ is a real $g\times g$ symmetric matrix,  its signature $s$
is the number of positive eigenvalues minus the number of negative eigenvalues. 
%We denote then by ${\rm det}'(S)$ the product of its non zero eigenvalues and by $S^{-1}$
%the Moore-Penrose inverse of $S$. Remember that when $S$ is invertible,
%${\rm det}'(S)={\rm det}(S)$ and the Moore-Penrose inverse of $S$ 
%is the usual inverse of $S$.

\bl 
\label{lemsqadetrie2} 
Let $\tau$ in $\mc H_g$ and $S$ be an invertible  real $g\times g$ symmetric matrix and $s$ be its signature. Then one has
\begin{equation}
\label{eqnhtauhtau}
h(\tau)\, h(-\tau^{-1}\!-\!S)
\;=\; e^{is\pi/4}\,h(\tau\!+\!S^{-1})\,|{\rm det}(S)|^{1/2}.
\end{equation} 
\el

\begin{proof}
The squares are equal and one checks the sign by looking at the limit when $\tau=it{\bf 1}_g$ with $t$ going to $0$.
\end{proof}

%82
\subsection{The sign in the transformation formula}
\label{secsigtrafor}

The following lemma give precise formulas for the theta cocycle $j(\si,\tau)$ that we introduced in Lemma \ref{lemtrafor}. We recall that, 
for $\si\!=\!\mbox{\scriptsize 
$\left(\!\begin{array}{cc} \al&\be\\  
\ga&\de\end{array}\!\right)$} \in {\rm Sp}_{g,\m Z}^\theta$
and $\tau\in \mc H_g$, the theta cocycle is defined by the equality 
\begin{eqnarray}
\label{eqnsigtrafor}
\theta
(0,\si\tau)
&=& j(\si,\tau)\;\,
\theta
(0,\tau).
\end{eqnarray}

To compute this cocycle we need to introduce a function
$
J_{1/2}(\si,\tau)
$
which is a continuous square root of 
the tangent cocycle $J(\si,\tau):={\rm det}_{\m C}(\ga\tau+\de)$.
The function $J_{1/2}(\si,\tau)$, for $\det(\de)\neq 0$ 
is given by the formula 
\begin{eqnarray}
\label{eqnj12sitau}
J_{1/2}(\si,\tau)
&:=&
|\det(\de)|^{1/2}\,
h(\tau)\, h(-\tau^{-1}\!-\!\de^{-1}\ga),
\end{eqnarray}
where $h(\tau)$ is the function defined in \eqref{eqnhtau}.
 
\bl 
\label{lemstasty} 
{\bf (Stark)} Let   $\si\in {\rm Sp}_{g,\m Z}^\theta$ 
with $\det(\de)\neq 0$ and $\tau\in \mc H_g$.\\ 
$(a)$ Then,  one has
\begin{eqnarray}
\label{eqnthecoc1}
j(\si,\tau)
&=&
\ka(\si)\;
J_{1/2}(\si,\tau)
\end{eqnarray}
where $\ka(\si)$ is the $8^{\rm th}$ root of unity given by the normalized Gauss sum
\begin{eqnarray}
\label{eqnkapgausum}
\ka(\si)
&=&
\frac{1}{|\det(\de)|^{1/2}}
\sum_{x\in \de^{-1}\m Z^g/\m Z^g}
e^{i\pi{}^tx{}^t\de\be x}
\end{eqnarray}
$(b)$ If  $\det(\de)=\pm d_0$ with $d_0>0$  odd, square-free
and $\det(\ga)\neq 0$, one has 
\begin{eqnarray}
\label{eqnthecoc2}
j(\si,\tau)
&=&
e^{is\pi/4}\,\eps_{d_0} \; (\!\tfrac{c}{d_0}\!)\;
|{\rm det}(\de^{-1}\ga)|^{1/2}\,
h(\tau\!+\!\ga^{-1}\de),
\mbox{ where}
\end{eqnarray}
- $s$ is the signature of the real symmetric matrix $S:=\de^{-1}\ga$,\\
- $\eps_{d_0}=1$ when $d_0\equiv 1$ mod $4$ and 
$\eps_{d_0}=-i$ when $d_0\equiv 3$ mod $4$,\\
- $c=2\,d_0\,{}^tm\de^{-1}\ga m$ where $m$ is any primitive vector of $\m Z^g$,\\
- $(\tfrac{c}{d_0})$ is a Jacobi symbol,
%- $\ga^{-1}\de$ is the Moore-Penrose inverse of the real symmetric matrix $\de^{-1}\ga$.
\el

This precise determination of the theta cocycle
has a long history. 

For $g=1$, these formulas which are due to Hecke were useful  in \cite{CSAGI}. 

When $g>1$, the first formulas \eqref{eqnthecoc1}  and \eqref{eqnkapgausum}
were kown  to Siegel and Igusa, 
see \cite[p.228]{Igusa64}. They  can be explicitely  found 
in \cite[2.2.26 p.169]{LiVe80}, \cite[p.7]{Sta82} or \cite[p. 26-27]{Freitag91}. 
There also exists a very simple formula for $\ka(\si)^2$ due to Igusa
in \cite[p.182]{Igusa72}. It says that the map $\si\mapsto \ka(\si)^2$
is a character of ${\rm Sp}_{g,\m Z}^\theta$ and when 
$\si$ is in ${\rm Sp}_{g,\m Z}^{2}$, one has 
$\ka(\si)^2= e^{i\frac{\pi}{2}\,{\rm tr}(\de-\mathds{1}_g)}$. 

The second formula \eqref{eqnthecoc2} is due to Stark in \cite{Sta82} 
for an odd prime $d_0$ and was extended by   Styer 
in \cite[Thm 2 p.660]{styer84a} to  a square-free $d_0$, even allowing a degenerate block $\ga$. 
The proof of \eqref{eqnthecoc2} relies on 
the first Formulas, on  \eqref{eqnhtauhtau} 
and on the classical formula for Gauss sums on cyclic groups. 
An extension of Formula \eqref{eqnthecoc2} in the case where 
$d_0$ is not square-free can be found in \cite{styer84b}.

The strong relationship between the transformation formula 
\eqref{eqnsigtrafor} and the Weil representation was discovered by A. Weil.
See \cite{LiVe80} and \cite{Friedberg85}. 
%Variations of these transformation formulas for number fields have been given by Eichler in  \cite{Eichler77} (see also \cite[section 5.1]{Garrett90}) by A. Schwartz in \cite{Schwartz93}, by Rhodes in \cite{Rhodes97} and by Richter in \cite{Richter00} or \cite{Richter04}.
\vs 

We will only need Lemma \ref{lemstasty} when 
$d_0=1$. This case is easier because  the 
Gauss sum 	and the Jacobi symbol are equal to $1$.

\bc
\label{corsigtrafor1} 
Let   $\si\in {\rm Sp}_{g,\m Z}^\theta$, $\tau\in \mc H_g$.
If $\det(\de)=\pm 1$, $\det(\ga)\neq 0$,
one has
\begin{eqnarray}
\label{eqnthecoc3}
j(\si,\tau)
&=& e^{is\pi/4}\,|{\rm det}(\de^{-1}\ga)|^{1/2}\,
h(\tau\!+\!\ga^{-1}\de).
\end{eqnarray}
\ec

%83
\subsection{Finding the sign of  critical values}
\label{secsigcrival}

We now explain the strategy  to compute sign of  the critical value $\la_\nu$ 
in the examples of Chapter \ref{seclisexa}. 
This strategy follows  the proof of Theorem \ref{thmmaiabe}.\vs 

We recall that  $(A=\m C^g/\La,\om)$ is a principally polarized abelian variety 
and that $\nu$ is  a unitary $\m Q$-endomorphism of $A$
preserving a theta structure of level $2$,
that $G_\nu$ is the finite abelian group
$
G_\nu:=\La/(\La\cap \nu\La).
$
Remember that $\nu$ can be thought both
as a unitary transformation $T_\nu$ of $\m C^n$
or as a symplectic transformation $h_\nu$ of $(\La_\m Q,\om)$, and that,
as in \eqref{eqnmnuhnu},  
in a symplectic basis
$e_1,\ldots,e_g,-f_1,\ldots,-f_g$ of $\La$, the multiplication by $\nu^{-1}$  is given by a symplectic rational matrix $m_\nu\in {\rm Sp}_{g,\m Q}^{\theta,2}$ 
preserving a theta structure of level $2$.
Since it is elliptic, this matrix $m_\nu$ has a fixed point $\theta\in \mc H_g$:
$$
m_\nu \theta=\theta
$$
By Proposition \ref{prohsidsi}, there exists $\si_{1}$ and $\si_{2}$ in 
${\rm Sp}(g,\m Z)$ and an integral matrix ${\bf d}$ 
such that the symplectic matrix $m_{\nu}$ 
can be written as 
\begin{eqnarray}
\label{eqnmnusdsd}
m_{\nu}&=&\si_{1}\;
D
\;\si_{2}\;\;
{\rm with}\;\;
\mbox{\small
$D=\left(\!
\begin{array}{cc} 
\!{}^t{\bf d}^{-1}\!&{\bf 0}\\  
{\bf 0}&{\bf d}
\end{array}
\!\right)$}.
\end{eqnarray} 
When  $g=2$, we will often but not always choose
${\bf d}=\mbox{\scriptsize 
$\left(\!\begin{array}{cc} 1&0\\  
0&d\end{array}\!\right)$}.$

By Lemma \ref{lemthestr}  the matrix $\si:=\si_{2}\si_{1}$ belongs to $ {\rm Sp}_{g,\m Z}^{\theta,2}$. 
This matrix 
$\si=\mbox{\scriptsize 
$\left(\!\begin{array}{cc} \al&\be\\  
\ga&\de\end{array}\!\right)$},$
and the element  $\tau=:\si_{1}^{-1}\theta\in \mc H_g$ 
are related by the equality $\si\tau=D\tau$ that is: 
$$
\si\tau = { }^t{\bf d} \tau {\bf d}.
$$
Therefore one can apply  Lemma \ref{lemstasty} to compute the critical value $\la_\nu$ given by Corollary \ref{cormaisym}.
In all the examples below,
since we have already computed $\la_\nu^2$ in Chapter \ref{seclisexa}, 
we will be able to avoid a few computations 
thanks to the following corollary:

\bc 
\label{corsqadetrie}
Assume $\det(\ga)\neq 0$. Let $s$ be the signature of $\de^{-1}\ga$.\\
$(a)$ If $g=1$ and $|\de|= 1$, then one has ${\rm Re}(\la_\nu)>0$.\\
$(b)$ If $g=2$ and $|\det(\de)|= 1$, then one has ${\rm Re}(e^{-is\pi/4}\la_\nu)>0$.
\ec

\begin{proof}
By Corollary \ref{cormaisym}, one has  
$
\la_\nu=j(\si',2\tau)\, |G_\nu|
$.
%We first assume $\det(\ga)\neq 0$ so that, 
By Formula \eqref{eqnthecoc3}, one has
$
\la_\nu=e^{is\pi/4}\,
h(\tau+\ga^{-1}\de) \, |det(\ga)|^{1/2}\, |G_\nu|.
$
We now use  Lemma \ref{lemsqadetrie}.\\
When $g=1$, one has $s=\pm 1$, and $-\frac{\pi}{4}< {\rm arg}(h(\tau\!+\!\ga^{-1}\de))<\frac{\pi}{4}$.\\ 
When $g=2$, one has $Re(h(\tau\!+\!\ga^{-1}\de))>0$. 
\end{proof}

Hence, for computing precisely $\la_\nu$,  the remaining two steps
beyond the steps $(a)$ to $(f)$  in  Chapter \ref{seclisexa} are :
\vs 

\noindent
$(g)$ 
We compute a decomposition 
$m_{\nu}=\si_{1} D\si_{2}$ 
with $\si_1$, $\si_2$ in ${\rm Sp}(g,\m Z)$.\\
$(h)$ We write $\si\!=\!\si_2\si_1=\mbox{\scriptsize 
$\left(\!\begin{array}{cc} \al&\be\\  
\ga&\de\end{array}\!\right)$}$ and 
conclude thanks to Corollary \ref{corsqadetrie}.

%84
\subsection{The sign for imaginary quadratic fields}
\label{secsigimaqua}

In this section we check the sign of the $d$-critical values 
from Proposition \ref{prorapirb0}. 
We will use freely the notation of Sections \ref{secimaqua} and \ref{secsigcrival}.
Remember that, from Part $(f)$ of Section \ref{secimaqua}, we know that   $\la_\nu=\pm(\sqrt{a}+i\sqrt{b})$.
\vs

\noindent $(g)$ 
We recall that in the basis $e_{1},\, - f_{1}$, one has
$m_{\nu}=\frac{1}{d}
\mbox{\scriptsize 
$\left(\!\begin{array}{cc} a\!-\!b&-ab\\  
4&a\!-\!b\end{array}\!\right)$} $
and that $\theta:= i\sqrt{ab}/2\in \mc H_1$ is fixed by $m_{\nu}$.
One checks that $m_{\nu}=\si_{1}D\si_{2}$ with 
$\si_{1}=
\mbox{\scriptsize 
$\left(\!\begin{array}{cc} \!a\!-\!b\!&\frac{a-b-1}{4}\\  
4&1\end{array}\!\right)$}$, 
\;\; 
$\mbox{\scriptsize 
$D=\left(\!
\begin{array}{cc} 
d^{-1}&0\\  
0&d
\end{array}
\!\right)$},$
and 
$\si_{2}=
\mbox{\scriptsize 
$\left(\!\begin{array}{cc} 1&\frac{a-b-d^2}{4}\\  
0&1\end{array}\!\right)$}.
$ 

Hence one has
$\si\tau=d^2\si$ 
with 
$\si=\si_{2}\si_{1}=
\mbox{\scriptsize 
$\left(\!\begin{array}{cc} \!2a\!-\!2b\!-\! d^2\!&a\!-\!(d\!+\!1)^2/4\\  
4&1\end{array}\!\right)$}$
and with
$\tau=\si_{1}^{-1}\,\theta =
\frac{(a-b-d^2)+2i\sqrt{ab}}{4d^2}.
$
By \eqref{eqnconrel}  this matrix $\si$ belongs to the  theta subgroup ${\rm Sp}^{\theta,2}_{1,\m Z}$ of level $2$.
%One also computes$\tau=P^{-1}\tau_0=\frac{a-b-d^2+2i\sqrt{ab}}{4d^2}$.
\vs

\noindent $(h)$ One has $g=1$ and  $\de=1$. Hence by 
Corollary \ref{corsqadetrie}, the real part of the  
critical value $\la_\nu$ is positive.
Hence, using Part $(f)$ of Section \ref{secimaqua},  we conclude that $\la_\nu=\sqrt{a}+i\sqrt{b}$.
And this value is always $d$-critical.

%85
\subsection{The sign for products of imaginary quadratic fields}
\label{secsigproima}

In this section we check the sign of the $d$-critical values 
from Proposition \ref{proparacri}. 
We will use freely the notation of Sections \ref{secproima} and \ref{secsigcrival}.
Remember that, from Part $(f)$ of Section \ref{secproima}, we know that   $\la_\nu=\pm(\sqrt{a_1}\!+\! i\eps_1\sqrt{b_1})(\sqrt{a_2}\!+\! i\eps_2\sqrt{b_2})$.
\vs

\noindent $(g)$ 
To simplify the calculation we will begin by the simpler symplectic basis
$e_{0,1},\ldots,e_{0,g},-f_{0,1},\ldots,-f_{0,g}$ 
for which the matrix $m_{0,\nu}\in {\rm Sp}(g,\m Q)$ of multiplication by $\nu^{-1}$ has been computed in \eqref{eqnmonumod}
\begin{equation*} 
m_{0,\nu}= \mbox{\scriptsize 
$\left(\!\begin{array}{cccc} u_1&0&\eps_1v_1&0\\  
0&u_2&0&\eps_2v_2\\
\eps_1w_1&0&u_1&0\\
0&\eps_2w_2&0&u_2\end{array}\!\right)$},
\;\; {\rm with}
\end{equation*}
\begin{equation*}
u_j\!:=\!\tfrac{a_j-b_j}{d'_j}\equiv 1\;{\rm mod}\;4,\;\;
v_j\!:=\!\tfrac{-2a_jb_j}{d'_j}\equiv 4\;{\rm mod}\;8,\;\; 
w_j\!:=\!\tfrac{2}{d'_j}\equiv 2\;{\rm mod}\;8.
\end{equation*}

One finds an explicit decomposition $m_{0,\nu}=\si_{0,1}D_0\si_{0,2}$ with 
\begin{eqnarray*}
\si_{0,1}&=&
\mbox{\scriptsize 
$\left(\!\begin{array}{cccc} \!a_1\!-\!b_1\!&0&\!\eps_1(a_1\!-\!b_1\!-\!1)/2\!&0\\  
0&\!a_2\!-\!b_2\!&0&\!\eps_2(a_2\!-\!b_2\!-\!1)/2\!\\
2\eps_1&0&1&0\\
0&2\eps_2&0&1
\end{array}\!\right)$}\\
\si_{0,2}&=&
\mbox{\scriptsize 
$\left(\!\begin{array}{cccc} 1&0&\!\eps_1(a_1\!-\!b_1\!-\!d_1^2)/2\!&0\\  
0&1&0&\!\eps_2(a_2\!-\!b_2\!-\!d_2^2)/2\!\\
0&0&1&0\\
0&0&0&1
\end{array}\!\right)$},
\end{eqnarray*} 
and $\;\; 
D_0=\mbox{\scriptsize 
$\left(\!
\begin{array}{cc} 
\!{}^t{\bf d}_0^{-1}\!&{\bf 0}\\  
{\bf 0}&{\bf d}_0
\end{array}
\!\right)$}$ where 
${\bf d}_0=
\mbox{\scriptsize 
$\left(\!\begin{array}{cc} 
d'_1&0\\  
0&d'_2
\end{array}\!\right)$}$.

The product
$\si_0=\si_{0,2}\si_{0,1}$ can be written as $\si_0=
\mbox{\scriptsize 
$\left(\!\begin{array}{cc} 
\al_0&\be_0\\  
\ga_0&\de_0
\end{array}\!\right)$}$ where
\begin{eqnarray}
\label{eqnsialbegade}
\nonumber
\al_0=\mbox{\scriptsize 
$\left(\!\begin{array}{cc} 
\! \!2a_1\!-\!2b_1\!-\! d_1^2\!\!\!&0\\  
0& \!\!\!2a_2\!-\!2b_2\!-\! d_2^2\!\!
\end{array}\!\right)$}\!
&,&\!
\be_0=\mbox{\scriptsize 
$\left(\!\begin{array}{cc} 
\!\!\eps_1(4a_1\!-\!(d_1\!+\!1)^2)/2\!\!\!&0\\  
0&\!\!\!\eps_2(4a_2\!-\!(d_2\!+\!1)^2)/2\!\!
\end{array}\!\right)$}\\
\ga_0=\mbox{\scriptsize 
$\left(\!\begin{array}{cc} 
2\eps_1&0\\  
0& 2\eps_2
\end{array}\!\right)$}
&,&
\de_0=\mbox{\scriptsize 
$\left(\!\begin{array}{cc} 
1&0\\  
0& 1
\end{array}\!\right)$}.
\end{eqnarray}

As we have seen in section \ref{secproima} we need 
to introduce a new lattice 
$\La$ with basis 
$(e_1,e_2,-f_1,-f_2)=(e_{0,1},\,e_{0,2},\,-f_{0,1},\,-f_{0,2})\,P$
where $P$ is a basis change matrix. 
This matrix $P$ is a scalar multiple of an element of
${\rm Sp}(g,\m R)$.
We need to distinguish two cases. 
A key difference between these two cases will be the value of the parameter $s$, 
when we will apply Corollary  \ref{corsqadetrie}.
\vs

{\bf First case:}
We assume that $\eps_1=-\eps_2=1$.

As in Section \ref{secproima}.$d1$, we choose
$
P=
\mbox{
$\left(\!\begin{array}{cc} {\bf p}&{\bf 0}\\  
{\bf 0}&{\bf p}
\end{array}\!\right)$}
$ 
with 
$
{\bf p}=
\mbox{\scriptsize 
$\left(\!\begin{array}{cc} 1\!&1\\  
1\!&\!-1\!
\end{array}\!\right)$}.
$ 
\vs 

\noindent $(g1)$  In this new basis, the matrix $m_\nu=P^{-1}m_{0,\nu}P$ of multiplication by $\nu^{-1}$ belongs to $ {\rm Sp}^{\theta,2}_{2,\m Q}$. One can write $m_\nu=\si_1D\si_2$
with 
$$
\si_1:=P^{-1}\si_{0,1}P\; ,\;\;
D=P^{-1}D_0P
\; ,\;\;
\si_2=P^{-1}\si_{0,2}P.
$$
Since the matrix $P$ is block diagonal, 
the matrix $D$ too, and one has 
$$
D=
\mbox{\scriptsize 
$\left(\!\begin{array}{cc} \!{}^t{\bf d}^{-1}\!&{\bf 0}\\  
{\bf 0}&{\bf d}
\end{array}\!\right)$}
\;\;{\rm
with}\;\; 
{\bf d}={\bf p}^{-1}{\bf d}_0{\bf p}=
\mbox{\scriptsize 
$\left(\!\begin{array}{cc} d'_+&d'_-\\  
d'_-&d'_+
\end{array}\!\right)$}
\;\; {\rm where}\;\;
d'_\pm=\frac{d'_1\pm d'_2}{2}.
$$
Both $\si_1$ and $\si_2$ belong to ${\rm Sp}(g,\m Z)$ 
and one has 
$$\si=\si_1\si_2=P^{-1}\si_0 P=
\mbox{\scriptsize 
$\left(\!\begin{array}{cc} \al&\be\\  
\ga&\de
\end{array}\!\right)$}=
\mbox{\scriptsize 
$\left(\!\begin{array}{cc} 
{\bf p}^{-1}\al_{_0}{\bf p}\!&\! {\bf p}^{-1}\be_{_0}{\bf p}\\  
 {\bf p}^{-1}\ga_{_0}{\bf p}\!&\! {\bf p}^{-1}\de_{_0}{\bf p}
\end{array}\!\right)$}.
$$
By  Lemma \ref{lemthestr},  this matrix $\si$ belongs to the  theta subgroup ${\rm Sp}^{\theta,2}_{2,\m Z}$ of level $2$.\vs

\noindent $(h1)$ 
Using \eqref{eqnsialbegade}, one checks that $\det(\de)=1$
and that the symmetric matrix 
$\ga^{-1}\de$ has signature $s=0$.
Hence by Corollary \ref{corsqadetrie}, the real part of the  
critical value $\la_\nu$ is positive.
Therefore, using Part $(f)$ of Section \ref{secproima},  we conclude that $\la_\nu=(\sqrt{a_1}\!+\! i\sqrt{b_1})(\sqrt{a_2}\!-\! i\sqrt{b_2})$.
And this value is always $d$-critical.
\vs

{\bf Second case:}
We assume that $\eps_1=\eps_2=1$.

As in Section \ref{secproima}.$d2$,  we choose
$
P=
\mbox{
$\left(\!\begin{array}{cc} {\bf 2}&{\bf p}\\  
{\bf 0}&{\bf 1}
\end{array}\!\right)$}
$ 
with 
$
{\bf p}=
\mbox{\scriptsize 
$\left(\!\begin{array}{cc} 1\!&1\\  
1\!&\!-1\!
\end{array}\!\right)$}.
$ 
\vs 

\noindent $(g2)$  
In this new basis, the matrix $m_\nu=P^{-1}m_{0,\nu}P$ of multiplication by $\nu^{-1}$ also belongs to $ {\rm Sp}^{\theta,2}_{2,\m Q}$.
Since the matrix $P$ is not block diagonal, 
neither is the matrix $P^{-1}DP$. 
This is why we introduce the commutator
$$
\si_{0,3}=PD_0^{-1}P^{-1}D_0=
\mbox{
$\left(\!\begin{array}{cc} {\bf 1}&{\bf x}\\  
{\bf 0}&{\bf 1}
\end{array}\!\right)$}
\;\; {\rm where}\;\;
{\bf x}={\bf p}-{\bf d}_0{\bf p}{\bf d}_0.
$$
One can write $m_\nu=\si_1D\si_2$
with 
$$
\si_1:=P^{-1}\si_{0,1}P\; ,\;\;
D:=D_0
\; ,\;\;
\si_2:=P^{-1}\si_{0,3}\si_{0,2}P.
$$
Both $\si_1$ and $\si_2$ belong to ${\rm Sp}(g,\m Z)$ and one has 
$\si=\si_2\si_1=P^{-1}\widetilde{\si}_0 P$ where 
$$
\widetilde{\si}_0:=\si_{0,3}\si_0=
\mbox{
$\left(\!\begin{array}{cc} \al_0\!+\!{\bf x}\ga_0&\be_0\!+\!{\bf x}\de_0\\  
\ga_0&\de_0
\end{array}\!\right)$}.
$$ 
By Lemma \ref{lemthestr},  this matrix $\si$ belongs to the  theta subgroup ${\rm Sp}^{\theta,2}_{2,\m Z}$ of level $2$.
This matrix
$\si=\mbox{\scriptsize 
$\left(\!\begin{array}{cc} \al&\be\\  
\ga&\de
\end{array}\!\right)$}$ has lower blocks 
$\ga= 2\ga_{_0}$ and 
$\de=\de_{_0}$
\vs

\noindent $(h2)$ The end of the argument is as in $(i1)$.
One checks that $\det(\de)=1$.
This time the symmetric matrix 
$\ga^{-1}\de$ is positive and hence has signature $s=2$.
Therefore by 
Corollary \ref{corsqadetrie}, the imaginary part of the  
critical value $\la_\nu$ is positive.
Hence, using again Part $(f)$ of Section \ref{secproima},  we conclude that $\la_\nu=(\sqrt{a_1}\!+\! i\sqrt{b_1})(\sqrt{a_2}\!+\! i\sqrt{b_2})$.
And this value is always $d$-critical.

%86
\subsection{The sign for quartic CM fields} 
\label{secsigquafie}

In this section we check the sign of the $d$-critical values 
from Proposition \ref{prorarbrc1}. 
We will use freely the notation of Sections \ref{secquafie} and \ref{secsigcrival}.
Remember from Part $(f)$ of Section \ref{secquafie}, we know that  $\la_\nu=\pm i^{(a-1)/2}
(\sqrt{a}\!+\!\sqrt{c}\! -\! i\sqrt{b\!-\! 2 \sqrt{ac}})$.\vs

As in Section \ref{secquafie}, we distinguish two cases, where $\eps=\pm 1$.
\vs

{\bf First case:}
We assume that $a\equiv \eps$ mod $4$,\; and $b\equiv c\equiv 0 $ mod $4$.

As in  Section \ref{secquafie}.$d1$,  we choose
$
P=
\mbox{
$\left(\!\begin{array}{cc} 
\!2{\bf p}\!&{\bf 0}\\  
{\bf 0}&\!2\,{}^t{\bf p}^{-1}\!
\end{array}\!\right)$}
$ 
with 
$
{\bf p}=
\mbox{\scriptsize 
$\left(\!\begin{array}{cc} 2\!&b/2\\  
0\!&\!1/2\!
\end{array}\!\right)$}.
$ 
\vs  

\noindent
$(g1)$  We have computed in \eqref{eqnmnuqd1}  the matrix $m_\nu$ 
of multiplication by $\nu^{-1}$ in the
corresponding basis $e_1$, $e_2$, $-f_1$, $-f_2$.
This matrix $m_\nu$ belongs to ${\rm Sp}^{\theta,2}_{2,\m Q}$.
Using the method of elimination by symplectic matrices,
we find explicit matrices ${\bf d}$, $\si_1$ and $\si_2$  
satisfying \eqref{eqnmnusdsd}. They are given by 
$
{\bf d}=\mbox{\scriptsize 
$\left(\!\begin{array}{cc} 1&0\\  
0&d\end{array}\!\right)$} 
$,
\begin{eqnarray*} 
\si_{1}
&=&
\mbox{\scriptsize 
$\left(\!\begin{array}{cccc}
1&
(\eps d\!-\! 1)(b\!+\!c)/2&
0&0\\  
0&1&0&0\\
0&1&1&0\\
1&\!(\eps d (d\!-\!2a)\!+\!d\!-\!2c)/4\!&\!(1\!-\!\eps d)(b\!+\!c)/2\!&1\end{array}\!\right)$},\\
\si_{2}
&=&
\mbox{\scriptsize 
$\left(\!\begin{array}{cccc}
\!1\!-\!2\eps (b\!+\!c)\!&
\!\eps (b\!+\!c)(2c\!-\!d\!+\!\eps)/2\!-\!c\!&
\!\eps(b\!+\!c)^2\!-\!4c\!&\!-2\eps(b\!+\!c)\!\\  
4&d\!-\!2c&-2b\!-\!2c&4\\
0&1&1&0\\
\eps&(\eps (d\!+\!2b)\!-\!1)/4&0&\eps\end{array}\!\right)$}.
\end{eqnarray*}
Since $d\equiv\eps$ mod $4$ and $b\equiv c\equiv 0$ mod $4$, 
these two matrices have integer coefficients.
Hence, by Lemma \ref{lemthestr},  the matrix $\si:=\si_{2}\si_{1}$ belongs to $ {\rm Sp}_{g,\m Z}^{\theta,2}$. \vs 

\noindent $(h1)$ If we write
$\si=\mbox{\scriptsize 
$\left(\!\begin{array}{cc} \al&\be\\  
\ga&\de\end{array}\!\right)$},$
the $2\times 2$ lower blocks are given by
$$
\ga=\mbox{\scriptsize 
$\left(\!\begin{array}{cc} 0&2\\  
2\eps&x\end{array}\!\right)$}
\;\;{\rm and}\;\; 
\de=\mbox{\scriptsize 
$\left(\!\begin{array}{cc} 1&0\\  
y&\eps\end{array}\!\right)$},
$$
where $x$ and $y$ are integers.
This proves that $\det(\de)=\eps=\pm 1$ and the symmetric matrix 
$\de^{-1}\ga=\mbox{\scriptsize 
$\left(\!\begin{array}{cc} 0&2\\  
2&\eps\, (x\!-\!2y)\end{array}\!\right)$},$
has signature $s=0$. 

Therefore by 
Corollary \ref{corsqadetrie}, the real part of the  
critical value $\la_\nu$ is positive.
Hence, using again Part $(f)$ of Section \ref{secquafie},  we conclude that\\ $\la_\nu=\sqrt{a}\!+\!\sqrt{c}\! -\! i\sqrt{b\!-\! 2 \sqrt{ac}}$\;
when $\eps =1$,\\
$\la_\nu=\sqrt{b\!-\! 2 \sqrt{ac}}\! +\!i\sqrt{a}\!+\!i\sqrt{c}$\;
when $\eps =-1$.\\
And these values and their complex conjugate  are always $d$-critical.
\vs

{\bf Second case:}
We assume that $a\equiv b\equiv c\equiv \eps $ mod $4$.

As in Section \ref{secquafie}.$d2$,  we choose
$
P=
\mbox{
$\left(\!\begin{array}{cc} \!2{\bf p}\!&{\bf p}\\  
{\bf 0}&\! 2\,{}^t{\bf p}^{-1}\!
\end{array}\!\right)$}
$ 
with 
$
{\bf p}=
\mbox{\scriptsize 
$\left(\!\begin{array}{cc} 2\!&b/2\\  
0\!&\!1/2\!
\end{array}\!\right)$}.
$ 
\vs 

\noindent $(g2)$ We have computed in \eqref{eqnmnuqd2}  the matrix $m_\nu$ of multiplication by $\nu^{-1}$ in the
corresponding basis $e_1$, $e_2$, $-f_1$, $-f_2$.
This matrix $m_\nu$ belongs to ${\rm Sp}^{\theta,2}_{2,\m Q}$.
Using the method of elimination by symplectic matrices,
we find explicit matrices ${\bf d}$, $\si_1$ and $\si_2$  
satisfying \eqref{eqnmnusdsd}. They are given by 
$
{\bf d}=\mbox{\scriptsize 
$\left(\!\begin{array}{cc} 1&0\\  
0&d\end{array}\!\right)$} 
$,
\begin{eqnarray*} 
\si_{1}
&=&
\mbox{\scriptsize 
$\left(\!\begin{array}{cccc}
1&
2b\!+\!2c\!+\! 2&
0&0\\  
-1&d\!-\!2c\!-\!2&
(b\!+\!c\!+\!1)h/2&-h/4\\
0&-4&\!1\!-\!2\eps b\!-\!2\eps c\!-\! 2\eps&\eps\\
2&-2d\!+\!4c&
\!-(b\!+\!c\!+\!1)(h\!+\!2\eps )\!&
\!(h\!+\!2\eps)/2\!\end{array}\!\right)$},\\
\si_{2}
&=&
\mbox{\scriptsize 
$\left(\!\begin{array}{cccc}
1&
\!-F\!&
\!\eps (s\!+\!1)^2\!-\!c\!&
\!-F/2\!\\  
-1&\!c\!+\!df/2\!&
\!(b\! +\! c\!+\!1)(1\!-\!\eps d)/2-1\!&\!(2c\!+\!d f\!-\!2)/4\!\\
0&2&1&1\\
0&2&0&1
\end{array}\!\right)$},
\end{eqnarray*}
where $h:=\eps d+2\eps b-1$, $f:=(2a\!-\!d\!-\!2)\eps\!-\!1$ and 
$F=(b\!+\!c\!+\!1)f\!+\!2c\!+\!2$.

Since $a\equiv b\equiv c\equiv \eps$ mod $4$, 
these two matrices have integer coefficients.
Hence, by Lemma \ref{lemthestr},  the matrix $\si:=\si_{2}\si_{1}$ belongs to $ {\rm Sp}_{g,\m Z}^{\theta,2}$. \vs 

\noindent $(h2)$ If we write
$\si=\mbox{\scriptsize 
$\left(\!\begin{array}{cc} \al&\be\\  
\ga&\de\end{array}\!\right)$},$
the $2\times 2$ lower blocks are given by
$$
\ga=\mbox{\scriptsize 
$\left(\!\begin{array}{cc} 0&-8\\  
0&-4\end{array}\!\right)$}
\;\;{\rm and}\;\; 
\de=\mbox{\scriptsize 
$\left(\!\begin{array}{cc} 1-4\eps (b\! +\!c\!+\!1)&2\eps\\  
-2\eps (b\! +\!c\!+\!1)&\eps\end{array}\!\right)$},
$$
where $x$ and $y$ are integers.
This proves that $\det(\de)=\eps=\pm 1$.
We can not apply our Corollary \ref{corsqadetrie} because the
symmetric matrix 
$\de^{-1}\ga=\mbox{\scriptsize 
$\left(\!\begin{array}{cc} 0&0\\  
0&-4\eps\end{array}\!\right)$}$
%has signature $s=-\eps$. 
is degenerate.
%Corollary \ref{corsqadetrie}, the real part of $e^{i\eps\pi/4}\la_\nu$ is positive.
Instead we use Formulas \eqref{eqnthecoc1} and compute
$$
j(\si',2\tau)=h(\tau)\,h(-\tau^{-1}-\de^{-1}\ga)=
\det({\bf 1}+\de^{-1}\ga\tau)^{ 1/2}=(1-4\eps\tau_{2,2})^{1/2}
$$
which always  has positive real part. Hence, $\la_\nu$ also has positive real part
and, using again Part $(f)$ of Section \ref{secquafie},  we conclude that\\ $\la_\nu=\sqrt{a}\!+\!\sqrt{c}\! -\! i\sqrt{b\!-\! 2 \sqrt{ac}}$\;
when $\eps =1$,\\
$\la_\nu=\sqrt{b\!-\! 2 \sqrt{ac}}\! +\!i\sqrt{a}\!+\!i\sqrt{c}$\;
when $\eps =-1$.\\
And these values and their complex conjugate  are always $d$-critical.

%87
\subsection{Conclusion}
\label{secconclusion}

I would like to end this paper by pointing out why this sign 
issue that we discussed in Chapter \ref{secsigcri} 
is so  delicate.
\vs 

It follows from the sign discussion 
in Section \ref{secsigimaqua}
that, for  $d= 5$,  the
value
$\la:=1\!+\!2i$ is a $d$-critical value.
But one can  check, using the Buchberger's
algorithm, that  the opposite value
$-\la$ is not a $d$-critical value.
\vs

Similarly, it follows from the sign discussion 
in Section \ref{secsigquafie}
that, for  $d= 15$,  the
value
$\la:=1\!+\!\sqrt{5}+i\sqrt{9\!-\!2\sqrt{5}}$ 
is an even $d$-critical value,
that is a value admitting an even $\la$-critical function.
As we have noticed in Remark \ref{remevencri}, all the critical values 
$\la$ that we construct in this article are even. 
But one can also check, using again the Buchberger's
algorithm,  that  the opposite value
$-\la$ is not an even $d$-critical value.

%9
{\small
\bibliography{theta2}
}
\vs 

{\small
\noindent
Y. \textsc{Benoist}: CNRS, 
Universit\'e Paris-Saclay,\hfill
\texttt{yves.benoist@u-psud.fr}}

\end{document}